\documentclass[a4paper, reqno]{amsart}

\AtBeginDocument{%
   \def\MR#1{}
}

\usepackage{mathtools}
\usepackage{amsthm}
\usepackage{amsfonts}
\usepackage{mathrsfs}
\usepackage[usenames,dvipsnames]{color}
\usepackage[all]{xy}
\usepackage{xr-hyper}
\usepackage[colorlinks=
   citecolor=Black,
   linkcolor=Red,
   urlcolor=Blue]{hyperref}
\usepackage{amssymb}
\usepackage{graphicx}
\usepackage[normalem]{ulem}
\usepackage{cancel}

\theoremstyle{plain}
\newtheorem{theorem}{Theorem}
\newtheorem{introthm}{Theorem}
\newtheorem{lemma}[theorem]{Lemma}
\newtheorem{prop}[theorem]{Proposition}
\newtheorem{cor}[theorem]{Corollary}

\theoremstyle{definition}
\newtheorem{definition}[theorem]{Definition}
\newtheorem{remark}[theorem]{Remark}
\newtheorem{intro-remark}{Remark}

\newtheorem{eg}[theorem]{Example}

\numberwithin{equation}{section}
\numberwithin{theorem}{section}

\DeclareMathOperator{\Nilp}{Nilp}
\DeclareMathOperator{\Gal}{Gal}
\DeclareMathOperator{\Norm}{N}
\DeclareMathOperator{\GL}{GL}
\DeclareMathOperator{\GU}{GU}
\DeclareMathOperator{\SU}{SU}
\DeclareMathOperator{\id}{id}
\DeclareMathOperator{\diag}{diag}
\DeclareMathOperator{\val}{val}
\DeclareMathOperator{\Stab}{Stab}
\DeclareMathOperator{\Adm}{Adm}
\DeclareMathOperator{\End}{End}
\DeclareMathOperator{\Lie}{Lie}
\DeclareMathOperator{\Spf}{Spf}
\DeclareMathOperator{\Tr}{Tr}
\DeclareMathOperator{\height}{ht}

\newcommand{\ad}{\mathrm{ad}}
\newcommand*{\LargerCdot}{\raisebox{-0.25ex}{\scalebox{1.2}{$\cdot$}}}

\title[Supersingular Locus of Unitary Shimura Varieties]{The
  Supersingular Locus of Unitary Shimura Varieties with Exotic Good
  Reduction}

\author{Haifeng Wu}

\address{Universit\"{a}t Duisburg-Essen} \email{haifeng.wu@uni-due.de}

\begin{document}
\maketitle

\begin{abstract}
  In this paper, we use a group-theoretic approach to give a concrete
  description of the geometric structure of the supersingular locus of
  unitary Shimura varieties with exotic good reduction. This approach
  also is a more uniform way to prove results of this form obtained
  previously by, for example, Vollaard-Wedhorn \cite{MR2800696} and
  Rapoport-Terstiege-Wilson \cite{MR3175176}. 
\end{abstract}

\section{Introduction}
\label{sec:intro}

We are interested in the geometry of the basic loci of Shimura
varieties, which may have important applications in the Langlands
program and Kudla program, for example, see the work of M.\,Harris \&
R.\,Taylor in \cite{MR1876802} and S.\,Kudla \& M.\,Rapoport in
\cite{MR2800697} \& \cite{MR3281653}. The basic locus of a Shimura
variety is the unique closed and, in some sense, the most interesting
Newton stratum. However, usually the geometric structure of the basic
locus cannot be described explicitly, for example, we even do not know
the dimension of the basic locus in the Siegel moduli spaces with
Iwahori level structure in the odd case (cf.~\cite[Theorem
1.1]{MR2915544}). Many mathematicians contributed to the general
problem of giving a concrete description of basic loci of Shimura
varieties. For the work in this area before 2005, we refer to the
introduction of \cite{MR2666394}. Let us review the work after 2005.

\begin{itemize}
\item I.\,Vollaard \& T.\,Wedhorn study the supersingular locus of the
  reduction of the Shimura variety for $\mathrm{GU}(1, n-1)$ at an
  inert prime $p$ in \cite{MR2800696}.
\item U.\,G\"{o}rtz \& C.-F.\,Yu study the supersingular locus of the
  Siegel modular varieties with Iwahori level structure
  $\mathcal{A}_{g, I}$ in \cite{MR2915544}. They show that if $g$ is
  even, the dimension of the supersingular locus is $g^{2}/2$. If $g$
  is odd, they give an estimate of the dimension of the supersingular
  locus. And in any case, the supersingular locus is not
  equidimensional if $g \geq 2$.
\item M.\,Rapoport, U.\,Terstiege \& S.\,Wilson study the supersingular
  locus of the Shimura variety for $\mathrm{GU}(1, n-1)$ over a
  ramified prime with the parahoric level structure given by a
  selfdual lattice in \cite{MR3175176}.
\item B.\,Howard \& G.\,Pappas study the supersingular locus of the
  Shimura variety for $\mathrm{GU}(2,2)$ at an inert prime in
  \cite{MR3272278}.
\item U.\,G\"{o}rtz \& X.\,He in \cite{GH} claim that the supersingular
  locus of the Shimura variety for $\mathrm{GU}(2,2)$ at a split prime
  can be written down similarly to \cite{MR3272278}.
\item In \cite{TianYichaoXiaoLiang} Tian and Xiao describe the basic
  locus in quaternionic Shimura varieties, and in
  \cite{HelmTianYichaoXiaoLiang}, Helm, Tian and Xiao analyze the case
  of Shimura varieties attached to $G(U(r, s) \times U(s, r))$. In
  both cases, they apply their results to study (and prove in certain
  cases) the Tate conjecture for the special fibers of these Shimura
  varieties.
\end{itemize}

In all the above cases except the G\"{o}rtz-Yu case, the supersingular
locus is a union of Ekedahl-Oort strata and admits a stratification by
classical Deligne-Lusztig varieties, and the index set and the closure
relations between strata can be described in terms of the Bruhat-Tits
building of a certain inner form of the underlying group. Such Shimura
varieties are called of Coxeter type in \cite{GH}. U.\,G\"{o}rtz \&
X.\,He study the analogous problem in the equi-characteristic case,
i.e.~the basic affine Deligne-Lusztig varieties of Coxeter type. They
give a complete list of ADLV of Coxeter type (cf.~\cite[Theorem
5.1.2]{GH}). In the mixed characteristic case, the affine
Deligne-Lusztig ``variety'' is a priori only a set, X.\,Zhu shows that
the ADLV has a perfect scheme structure, and, when the underlying
group is unramified, it is canonically isomorphic to the perfection of
the special fiber of its associated Rapoport-Zink space in his mixed
affine Grassmannian in \cite[Proposition 0.4]{zhu} as perfect schemes.

Recently, M.\,Chen \& E.\,Viehmann claim that they can give a complete
description of the Shimura variety for $\mathrm{GU}(2, n-2)$ at an
inert prime in \cite{ChenViehmann}.

This paper is a contribution to the program of giving a concrete
description of the basic loci of the Shimura varieties of Coxeter
type. We use a group-theoretic approach to study the basic locus of
unitary Shimura varieties with exotic good reduction. Although we
focus on a concrete case, our method gives a uniform way to study all
the basic loci of Coxeter type Shimura varieties of PEL type. Let us
talk about the main steps and ingredients for the general strategy.
\begin{enumerate}
\item A priori, we need a suitable Rapoport-Zink space $\mathcal{M}$.
  By the local model diagram, this is equivalent to requiring a
  suitable local model $\mathbf{M}$. In general, the naive local model
  $\mathbf{M}^{\mathrm{naive}}$ is not flat and the honest local model
  $\mathbf{M}^{\mathrm{loc}}$ may not have a moduli description. For
  the purposes of studying basic loci, a topologically flat local
  model $\mathbf{M}$ is enough for us because, by definition, the
  basic locus is reduced. The key point of proving the topological
  flatness of $\mathbf{M}$ is to show that the set of
  Kottwitz-Rapoport strata is the same as the $\mu$-admissible set.
\item Using Dieudonn\'{e} theory, establish the bijection from the
  affine Deligne-Lusztig variety to the
  $\overline{\mathbb{F}}_{p}$-valued points of $\mathcal{M}$. The
  surjectivity would follow from the topological flatness of
  $\mathbf{M}$. This essentially means that one proves
  \cite[Proposition 0.4]{zhu} in the setting at hand.
\item Using Lusztig's partition, get some version of the Crucial Lemma
  (cf. \cite[Lemma 2.1]{MR2666394}, see also
  Section~\ref{sec:essential-lemma}).
\item Identify the Deligne-Lusztig varieties on the group side with
  the open BT strata on the RZ side as schemes using Lusztig's
  partition (or equivalently, the Crucial Lemma).
\item Identify closures on both sides using the normality of the
  closure of Deligne-Lusztig varieties (cf. Remark~\ref{rmk:15}). The
  closure relation would be described in terms of the Bruhat-Tits
  building.
\end{enumerate}

Now, let us return to our concrete case. We study the supersingular
locus of the unitary Shimura varieties for $\mathrm{GU}(1, n-1)$ at a
ramified prime with special parahoric level structure. More precisely,
let $E$ be an imaginary quadratic field extension of $\mathbb{Q}$
together with a ramified rational prime $p \geq 3$. Let $(W, \varphi)$
be a hermitian space of signature $(1, n-1)$, $\mathbb{G}$ the
corresponding unitary similitude group. Let $C_{p}$ be the special
parahoric subgroup corresponding to the $0$-th vertex of the local
Dynkin diagram~\eqref{eq:244} and \eqref{eq:245}, $C^{p}$ a
sufficiently small open compact subgroup of
$\mathbb{G}(\mathbb{A}^{p}_{f})$. Let $\mathcal{A}$ be the integral
model of the Shimura variety $\mathrm{Sh}_{C^{p}}(\mathbb{G},h)$, then
$\mathcal{A}$ is smooth by \cite[Proposition 4.16]{MR2595559}. The
smoothness of $\mathcal{A}$ is unexpected because $p$ is ramified, so
we use the terminology ``exotic good reduction''.

The supersingular locus of the special fiber of $\mathcal{A}$ can be
studied using Rapoport-Zink's $p$-adic uniformization theorem. Now let
us consider the corresponding Rapoport-Zink spaces.

Let $F$ be a ramified quadratic field extension of $\mathbb{Q}_{p}$,
together with the unique non-trivial automorphism $\bar{} \in
\mathrm{Gal}(F/\mathbb{Q}_{p})$ and the uniformizer $\pi$ such that
$\bar{\pi} = -\pi$. We denote $L$ the completion of the maximal
unramified field extension of $\mathbb{Q}$ and let $\breve{F}
\coloneqq F \otimes_{\mathbb{Q}_{p}} L$. Let $\mathbb{F}$ denote the
algebraically closed field $\overline{\mathbb{F}}_{p}$.

For an $\mathbb{F}$-scheme $S$, a unitary $p$-divisible group of
signature $(1, n-1)$ over $S$ (cf.~\cite[3.1]{RSZ}) is a triple
$(X, \iota_{X}, \lambda_{X})$, where $\iota_{X}$ is an
$\mathcal{O}_{F}$-action satisfying the Kottwitz condition, the Wedge
condition and the extra Spin condition if $n$ is even. The
polarization $\lambda_{X}$ satisfies the condition that the Rosati
involution on $\mathrm{End}(X)$ attached to $\lambda_{{X}}$ induces
the non-trivial automorphism on $\mathcal{O}_{F}$ over
$\mathbb{Q}_{p}$. Furthermore, the periodicity condition is assumed:
if $n$ is even, $\ker(\lambda_{X}) = X[\iota_{X}(\pi)]$; if $n$ is
odd, $\ker(\lambda_{X}) \subset X[\iota_{X}(\pi)]$ is of height $n-1$.
 
We fix a supersingular unitary $p$-divisible group
$(\mathbb{X}, \iota_{\mathbb{X}}, \lambda_{\mathbb{X}})$ of signature
$(1, n-1)$ over $\mathbb{F}$, and consider the moduli functor
$\mathcal{N}^{e}$:
\begin{align}
  \label{eq:248}
  (\mathbb{F} \text{-schemes}) & \longrightarrow (Sets), \\
  S & \longmapsto \{(X, \iota_{X}, \lambda_{X}, \rho_{X}) / \cong\},
      \nonumber
\end{align}
where $(X, \iota_{X}, \lambda_{X})$ is a unitary $p$-divisible group
and $\rho_{X}$ is an $\mathcal{O}_{F}$-linear quasi-isogeny such that
$\rho^{*}(\lambda_{\mathbb{X}})$ and $\lambda_{X}$ differ locally on
$\overline{S}$ by a scalar in $\mathbb{Q}_{p}^{\times}$. Then
$\mathcal{N}^{e}$ is of relative dimension $n-1$ and has the same
underlying topological space with the honest Rapoport-Zink space.
(cf.~Proposition~\ref{pr:9}). 

Let $G = \mathrm{GU}(N, \varphi)$ be the unitary similitude group of
signature $(1, n-1)$ where $(N, b \sigma)$ is the isocrystal given by
the framing object
$(\mathbb{X}, \iota_{\mathbb{X}}, \lambda_{\mathbb{X}})$ and $\varphi$
is the hermitian form corresponding to the polarization
$\lambda_{\mathbb{X}}$. Let $K = \mathrm{Stab}(\mathbb{M})$ be the
special parahoric subgroup corresponding to the $0$-th vertex of the
local Dynkin diagram of $G$ (see \eqref{eq:244} and \eqref{eq:245})
and $\mu$ the geometric minuscule cocharacter $(1, 0^{n-1};1)$. Then,
via Dieudonn\'{e} theory, we have a bijection
\begin{align}
  \label{eq:249}
  \Phi \colon
  X(\mu, b)_{K} & \longrightarrow \mathcal{N}^{e}(\mathbb{F}), \\
  g & \longmapsto g \mathbb{M}, \nonumber
\end{align}
where $X(\mu, b)_{K}$ is a union of affine Deligne-Lusztig
varieties. Then the map $\Phi$ induces a scheme structure on the left
hand side. Let $X(\mu, b)_{K}'$ (resp. $\mathcal{S}$) be the connected
component with trivial Kottwitz invariant of $X(\mu, b)_{K}$
(resp. $\mathcal{N}^{e}$). In \cite{GH} G\"{o}rtz-He show that the
affine Deligne-Lusztig variety is a disjoint union of fine affine
Deligne-Lusztig varieties (aka. Ekedahl-Oort strata)
\begin{equation}
  \label{eq:250}
  X(\mu, b)'_{K} = \biguplus_{w \in
    \mathrm{EO}_{\mathrm{cox}}} X^{f}_{w}(b),
\end{equation}
and each Ekedahl-Oort stratum is a disjoint union of classical
Deligne-Lusztig varieties
\begin{equation}
  \label{eq:251}
  X^{f}_{w}(b) \cong \coprod_{j \in \mathbb{J}/ \mathbb{J}\cap
    P_{\tilde{\mathbb{S}} - \Sigma}} j \cdot Y_{\Sigma^{\sharp}}(w), 
\end{equation}
where
\begin{equation}
  \label{eq:252}
      Y_{\Sigma^{\sharp}}(w) = \{ g \in P_{\tilde{\mathbb{S}}-\Sigma} /
    P_{\Sigma^{\sharp}} : g^{-1} b_{\ad} \sigma(g) \in P_{\Sigma^{\sharp}} w P_{\Sigma^{\sharp}} \}.
\end{equation}

For the framing object, we associate to it a hermitian space $C$. A
lattice $\Lambda$ in $C$ is called a vertex lattice if
$ \Lambda \subset \Lambda^{\sharp} \subset \pi^{-1} \Lambda$, where
$\Lambda^{\sharp}$ is the dual of $\Lambda$. The dimension of the
$\mathbb{F}_{p}$-vector space $\Lambda / \pi \Lambda^{\sharp}$ is
called the type of the lattice, denoted by $t(\Lambda)$. Let
$\mathcal{B}$ be the set of vertex lattices. Via the crucial lemma
(cf. Lemma~\ref{lm:3}), each basic EO element
$w \in \mathrm{EO}_{\mathrm{cox}}$ is attached to a vertex lattice
$\Lambda$. And we can show that the map $\Phi$ induces an isomorphism
from the closure of the Deligne-Lusztig variety
$Y_{\Sigma^{\sharp}(w)}$ to a closed subscheme $\mathcal{S}_{\Lambda}$
of $\mathcal{S}$.

Using these group-theoretic results and Smithling's result in
\cite{smithling-moduli}, via the map $\Phi$, we have the main theorem.

\begin{introthm}[see Theorem~\ref{thr:bt}]
  \label{thr:main}
  \leavevmode
  \begin{enumerate}
  \item There is a stratification, which is called the Bruhat-Tits
    stratification, of $\mathcal{S}$ by locally closed subschemes
    \begin{equation}
      \label{eq:22444}
      \mathcal{S} = \biguplus_{\Lambda \in \mathcal{B}} 
      \mathcal{S}_{\Lambda}^{\circ},
    \end{equation}
    and each stratum is isomorphic to the Deligne-Lusztig variety
    associated to the orthogonal group
    $\mathrm{SO}(\mathbb{B}_{\Lambda})$ and a $\sigma$-Coxeter
    element. The closure of each stratum
    $\mathcal{S}_{\Lambda}^{\circ}$ in $\mathcal{S}$ is given by
    \begin{equation}
      \label{eq:225555}
      \overline{\mathcal{S}_{\Lambda}^{\circ}} = \biguplus_{\Lambda'
        \subset \Lambda} \mathcal{S}_{\Lambda'}^{\circ} 
      = \mathcal{S}_{\Lambda}.
    \end{equation}
  \item The scheme $\mathcal{S}$ is geometrically connected of pure
    dimension $[\frac{n-1}{2}]$. The irreducible components of
    $\mathcal{S}$ are those $\mathcal{S}_{\Lambda}$ with
    $t(\Lambda) = n$.
  \end{enumerate}
\end{introthm}

Then, using the $p$-adic uniformization theorem, we have the
description of the supersingular locus of $\mathcal{A} \otimes
\mathbb{F}$.

\begin{introthm}[see Theorem~\ref{thr:conclusion}]
  \label{thr:ss}
  The supersingular locus $\mathcal{A}^{\mathrm{ss}}_{\mathbb{F}}$ is
  of pure dimension $[\frac{n-1}{2}]$. We have natural bijections
  \begin{equation}
    \label{eq:242222}
    \{ \text{irreducible components of } \mathcal{A}^{\mathrm{ss}}_{\mathbb{F}} \}
    \stackrel{1:1}{\longrightarrow} \mathbb{I}(\mathbb{Q}) \backslash
    (J(\mathbb{Q}_{p}) / K_{\mathrm{max}} \times
    \mathbb{G}(\mathbb{A}^{p}_{f}) / C^{p}),
  \end{equation}
  and 
  \begin{equation}
    \label{eq:243333}
    \{ \text{connected components of } \mathcal{A}^{\mathrm{ss}}_{\mathbb{F}} \}
    \stackrel{1:1}{\longrightarrow} \mathbb{I}(\mathbb{Q}) \backslash
    (J(\mathbb{Q}_{p}) / J^{0} \times \mathbb{G}(\mathbb{A}^{p}_{f})/
    C^{p}). 
  \end{equation}
  where $J^{0}$ is the subgroup of $J(\mathbb{Q}_{p})$
  consisting of those $j$ with trivial Kottwitz invariant and
  $K_{\mathrm{max}}$ is the stabilizer of some maximal-type vertex
  lattice in $J(\mathbb{Q}_{p})$.
\end{introthm}

This paper is structured as follows. In Section 2 we collect some
group data from the literature. In Section 3 we establish the
bijection between the Rapoport-Zink space and the affine
Deligne-Lusztig variety. In Section 4 we describe the set-theoretic
structure of the Rapoport-Zink space using G\"{o}rtz-He's
group-theoretic result. In Section 5 we establish the Bruhat-Tits
stratification scheme-theoretically. In Section 6 using the $p$-adic
uniformization theorem we describe the supersingular locus.

\subsection*{Acknowledgment}
\label{sec:acknowledgement}

I would like to give special thanks to my advisor Prof. Ulrich
G\"{o}rtz for introducing me to the area of the reductions of Shimura
varieties. I am very grateful that he spent innumerable hours on
discussion with me and provided infinite patience and encouragement to
me. This work was supported by the SFB/TR 45 ``Periods, Moduli Spaces
and Arithmetic of Algebraic Varieties'' of the DFG.

\section{Group data}

In this section we collect some combinatorial results from the literature, which
will be used later.

\subsection{Notations}

We list some notations which will be used through the whole
paper. Let $p$ be an odd prime number, $F$ a ramified quadratic field
extension of $\mathbb{Q}_{p}$. We denote by
$\bar{} \in \Gal(F/\mathbb{Q}_{p})$ the non-trivial automorphism. Let
$\pi$ be a uniformizer of $F$ such that $\bar{\pi} = - \pi$ and
$\pi^{2} = \varpi$, where $\varpi = \epsilon p$ is a uniformizer of
$\mathbb{Q}_{p}$ and $\epsilon$ is a unit in $\mathbb{Z}_{p}$. We
denote $L$ the completion of the maximal unramified field extension of
$\mathbb{Q}_{p}$ and $\breve{F} = F \otimes_{\mathbb{Q}_{p}} L$. Let
$\sigma$ be the Frobenius automorphism of $L / \mathbb{Q}_{p}$. Let
$\Gamma = \Gal(\overline{L}/L)$.

\subsection{Hermitian forms over local fields}

Let $E/E_{0}$ be a quadratic extension of local fields of mixed
characteristic $(0,p)$, $W$ an $n$-dimensional vector space over $E$
together with a non-degenerate hermitian form
\begin{equation}
  \varphi \colon W \times W \to E
\end{equation}
with respect to the quadratic extension $E/E_{0}$, i.e. $\varphi$ is
$E$-linear in the first factor and $^{*}$-linear in the second factor,
where $^{*} \in \Gal(E/E_{0})$ is the non-trivial automorphism. The
pair $(W, \varphi)$ is called a \emph{hermitian space}.

The isomorphism classes of hermitian forms can be determined by their
discriminants in the group $E_{0}^{\times}/ \mathrm{N}_{E / E_{0}}
E^{\times}$ by~\cite[Theorem 3.1]{Jacobo}. 

\begin{definition}
  A hermitian space $(W,\varphi)$ is called \emph{split} if it has
  trivial discriminant, i.e. the image of $(-1)^{n(n-1)/2} \det W$ in
  the group $E_{0}^{\times} / \Norm_{E/E_{0}} E^{\times}$ is trivial,
  otherwise $(W,\varphi)$ is called \emph{non-split}.
\end{definition}

\begin{remark}
  \label{rmk:17}
  The local class field theory shows the group
  $E_{0}^{\times} / \Norm_{E/E_{0}} E^{\times}$ is of order $2$. If
  the field extension $E/E_{0}$ is ramified, the group
  $E_{0}^{\times} / \Norm_{E/E_{0}} E^{\times}$ is generated by the
  units in $\mathcal{O}_{E_{0}}$. Therefore, in the ramified case,
  when $n$ is odd, there is only one similarity class of hermitian
  forms; when $n$ is even, there are two similarity classes of
  hermitian forms.
\end{remark}

\begin{prop}
  \label{pr:19}
  For an $n$-dimensional hermitian space $(W, \varphi)$ with respect
  to $E/E_{0}$, let $\mathrm{SU}(W, \varphi)$ be the special unitary
  group over $E_{0}$. Then when $n$ is odd, $\mathrm{SU}(W, \varphi)$
  is always quasi-split; when $n$ is even, $\mathrm{SU}(W, \varphi)$
  is quasi-split if and only if the hermitian form $\varphi$ is split.
\end{prop}

\begin{proof}
  We have the \emph{Witt decomposition}
  \begin{equation}
    \label{eq:21}
    W = H_{1} \oplus \cdots \oplus H_{q} \oplus W_{0},
  \end{equation}
  where $H_{i}$ is a hyperbolic plane for all $i$, $W_{0}$ is anisotropic of at
  most dimension $2$ by \cite[63:19]{MR1754311}. When $n$ is odd, $W_{0}$ is a
  line; when $n$ is even, $W$ is split if and only if $W_{0} = 0$. Let $S$ be
  the maximal $E_{0}$-split torus with respect to the
  decomposition~\eqref{eq:21}. Then, by definition, $\mathrm{SU}(W, \varphi)$ is
  quasi-split if and only if the centralizer of $S$ is a maximal torus, which is
  equivalent to the condition that the $E_{0}$-rank of $S$ is $[\frac{n}{2}]$,
  i.e. $W_{0}$ is a line when $n$ is odd and $W_{0} = 0$ when $n$ is even.
\end{proof}

We are interested in lattices in hermitian spaces. A lattice $M$ in
$(W,\varphi)$ is called \emph{$\gamma$-modular} if
$M^{\vee}=\gamma^{-1}M$, where $M^{\vee}$ is the dual lattice of $M$
with respect to $\varphi$ and $\gamma$ is a uniformizer of $E$; $M$ is
called \emph{nearly $\gamma$-modular}\footnote{Here we adopt the
  terminology in \cite{RSZ}.}  if
$M \subset M^{\vee} \stackrel{1}{\subset} \gamma^{-1}M$, where the
symbol $\stackrel{k}{\subset}$ means that the quotient of the
inclusion is of dimension $k$ over the residue field of $E$.

\begin{lemma}
  \label{split-lemma}
  Let $(W,\varphi)$ be an $n$-dimensional hermitian space with respect
  to the ramified field extension $\breve{F} / L$. Then when $n$ is
  odd, $\varphi$ is similar to a split hermitian form; when $n$ is
  even, $\varphi$ is split if and only if it contains a $\pi$-modular
  lattice.
 \end{lemma}

\begin{proof}
  If $n$ is odd, there exists $a \in L^{\times}$ such that $a \varphi$ has
  trivial discriminant because the group $L^{\times} /
  \Norm_{\breve{F}/L}\breve{F}^{\times}$ is generated by the units
  $\mathcal{O}_{L}^{\times}$. If $n=2m$ is even, \cite[Proposition
  8.1(b)]{Jacobo} shows that $\phi$ is split if and only if $W$ contains a
  $\pi$-modular lattice.
\end{proof}

\subsection{Combinatorics}
\label{sec:group-data}

Let $(V,\phi)$ be an $n$-dimensional split hermitian space over $F$,
$(e_{1},\ldots, e_{n})$ a basis such that
$\phi(e_{i},e_{j})=\delta_{i,n+1-j}$. Let $G=\GU(V,\phi)$ be the
general unitary group defined over $\mathbb{Q}_{p}$, i.e. for each
$\mathbb{Q}_{p}$-algebra $R$,
\begin{equation}
\label{eq:107}
G(R) = \left\{ g \in \GL_{F \otimes_{\mathbb{Q}_{p}} R} (V \otimes
_{\mathbb{Q}_{p}} R) \left|
  \begin{array}{ll}
    \phi(gv, gw) = c(g) \phi(v,w) \\
    \text{for some } c(g) \in R^{\times} \\
    \text{and for any } v, w \in V.
  \end{array}
  \right.
\right\}.
\end{equation}
The algebraic group $G$ is a reductive group over $\mathbb{Q}_{p}$,
and its derived group $G_{\mathrm{der}} = \SU(V,\phi)$ is semisimple
and simply connected. Let $D$ be the torus $G/G_{\mathrm{der}}$. We
identify $\pi_{1}(G) = \pi_{1}(D) = X_{*}(D)$.

Let $S \subset G$ be the maximal $L$-split torus consisting of
diagonal matrices defined over $\mathbb{Q}_{p}$, $T$ its centralizer,
$N$ its normalizer. Then $T$ is a maximal torus of $G$ because $G$ is
quasi-split. Over $\breve{F}$, we have the following isomorphism:
\begin{align}
\label{eq-split}
G_{\breve{F}} & \simeq  \GL_{n,\breve{F}} \times
                         \mathbb{G}_{\mathrm{m},\breve{F}},\\
g & \mapsto  (g_{0}, c(g)), \nonumber
\end{align}
where $g_{0} = (g_{i, j}^{0} \cdot g_{i, j}^{1})_{i,j}$, if we write $g =
(g_{i,j}^{(0)} \otimes g_{i, j}^{(1)})_{i, j}$ with $g^{(0)}_{i, j} \in F$ and
$g_{i,j}^{(1)} \in \breve{F}$. Then via the identification ~\eqref{eq-split},
the action of the non-trivial automorphism $\bar{} = \bar{} \otimes \id_{L} \in
\Gal(\breve{F} / L)$ on RHS is given by the map $(g_{0},c) \mapsto (\overline{c
  \phi^{-1} \prescript{t}{}{g_{0}}^{-1} \phi}, \bar{c}) $.

In this section, we collect some group data from \cite{MR546588}
\cite{Haines-Rapoport} \cite{MR2435422} \cite{MR2516305} \cite{MR2764885}
\cite{MR3177281}.

\subsubsection{Affine root systems and Iwahori-Weyl groups}
\label{sec:affine-root-systems}

First of all, we will compute the relative root system
$(X^{*}, X_{*}, \Phi, \Phi^{\vee})$ of $G$ and its Iwahori-Weyl group.

\paragraph{(a) odd case.} We write $n = 2m + 1$. Then
\begin{align*}
  S(L) & =  \{\diag(s_{1}, \ldots, s_{n}) : s_{i} \in L^{\times} \text{
    and } s_{1}s_{n} = \dots = s_{m}s_{m+2} = s_{m+1}^{2} \}, \\
  T(L) & =  \{\diag(t_{1}, \ldots, t_{n}) : t_{i} \in
             \breve{F}^{\times} \text{ and } t_{1} \bar{t}_{n} = \dots
             = t_{m} \bar{t}_{m+2} = t_{m+1} \bar{t}_{m+1}  \}.
\end{align*}
Under the identification ~\eqref{eq-split}, $X_{*}(T)$ can be
identified with $\mathbb{Z}^{n} \times \mathbb{Z}$. And
$X_{*}(T)_{\Gamma}$ is identified with $\mathbb{Z}^{m} \times
\mathbb{Z}$ under the canonical projection $X_{*}(T) \to
X_{*}(T)_{\Gamma}$. Let $X_{*} = X_{*}(T)_{\Gamma} \otimes \mathbb{R}
= \mathbb{R}^{m} \times \mathbb{R}$, then we identify $X_{*}(S)$ with
its image $2X_{*}(T)_{\Gamma}$ in $X_{*}$.

Similarly, $X^{*}(T) = \mathbb{Z}^{n} \times \mathbb{Z}$, so we may
identify $X^{*}(S)$ with $\mathbb{Z}^{m} \times \mathbb{Z}$ under the
canonical projection $X^{*}(T) \to X^{*}(T)_{\Gal(\breve{F}/L)} /
\text{torsion} = X^{*}(S)$. Let $X^{*} = X^{*}(S) \otimes \mathbb{R}$.
Then the set of roots $\Phi$ is just the image of the absolute roots
$\Phi(T,G)$, which is of type $\mathrm{A}_{n-1}$ by \eqref{eq-split},
under the natural map $X^{*}(T) \to X^{*}(S)$. Let $\epsilon_{i} \in
X^{*}$ be the function on $X_{*}$ sending $(x_{1}, \ldots, x_{m}; y)$
to $x_{i}$ for any $i \in \{1,2,\ldots,m\}$. Then
\begin{equation}
\label{eq:10}
\Phi =
\left\{
    \begin{array}{ll}
      \pm \epsilon_{i},  & 1 \leq i \leq m ,\\
      \pm 2 \epsilon_{i}, & 1 \leq i \leq m, \\
      \pm \epsilon_{i} \pm \epsilon_{j}, & 1\leq i < j \leq m
    \end{array}
  \right\}.
\end{equation}
So $\Phi$ is non-reduced. 

Let's look at the set of affine roots $\Phi_{a}$, by \cite[Proposition
2.2]{MR2516305},
\begin{equation}
  \Phi_{a} =
  \left\{
    \begin{array}{ll}
      \pm \epsilon_{i} + \frac{1}{2} \mathbb{Z}, & 1 \leq i \leq m, \\
      \pm 2 \epsilon_{i} + \frac{1}{2} + \mathbb{Z}, & 1 \leq i \leq m, \\
      \pm \epsilon_{i} \pm \epsilon_{j} + \frac{1}{2} \mathbb{Z}, & 1
                                                                   \leq
                                                                   i <
                                                                   j
                                                                   \leq m
    \end{array}
  \right\}.
\end{equation}
So the affine hyperplanes associated to $\Phi_{a}$ can be viewed as the
zero loci of the affine functions
\begin{equation}
  \label{eq:11}
  \left\{
    \begin{array}{l}
    \pm 2\epsilon_{i} + \frac{1}{2} \mathbb{Z}, \\
    \pm \epsilon_{i} \pm \epsilon_{j} + \frac{1}{2} \mathbb{Z},
    \end{array}
  \right.
\end{equation}
which can be viewed as an affine root system of type $C_{m}$. Let
$W_{0} = N(L)/T(L)$ be the Weyl group, which is isomorphic to
$\mathfrak{S}_{m} \rtimes \{\pm 1\}^{m}$ in the spirit of
\eqref{eq:11}. The affine Weyl group
$W_{a} = X_{*}(T^{\mathrm{sc}})_{\Gamma} \rtimes W_{0} \cong
\mathbb{Z}^{m} \rtimes W_{0}$,
where $T^{\mathrm{sc}}$ is $T \cap \SU(V,\phi)$, and the Iwahori-Weyl
group $\tilde{W} =X_{*}(T)_{\Gamma} \rtimes W_{0}$ is isomorphic to
$W_{a} \rtimes \pi_{1}(G)_{\Gamma}$, where $\pi_{1}(G)_{\Gamma}$ is
isomorphic to
$X_{*}(T)_{\Gamma}/X_{*}(T^{\mathrm{sc}})_{\Gamma} = \mathbb{Z}$.

Following \cite[1.8]{MR546588}, we choose a basis of $\Phi_{a}$
\begin{equation}
  \label{eq:4}
  \left\{
    \begin{array}{ll}
      \alpha_{i} = \epsilon_{m+1-i} - \epsilon_{m-i}, & 1 \leq i \leq m-1,
      \\
      \alpha_{m} = 2 \epsilon_{1}, & \\
      \alpha_{0} = \frac{1}{2} - \epsilon_{m}, & \\
    \end{array}
  \right.
\end{equation}
then we get the local Dynkin diagram of type $C\text{-}BC_{m}$.
\begin{equation}
  \label{eq:244}
  \xymatrix{
    \underset{\alpha_{0}}{\circ} & \underset{\alpha_{1}}{\circ} \ar@{=}[l] |-{\SelectTips{eu}{11}\object@\dir{>}} & \underset{\alpha_{2}}{\circ} \ar@{-}[l] &
    {\circ} \ar@{-}[l] &  \ar@{.}[l] &{\circ} \ar@{.}[l] &{\circ} \ar@{-}[l] &
    \underset{\alpha_{m}}{\circ} \ar@{=}[l] |-{\SelectTips{eu}{11}\object@\dir{>}}
  }.
\end{equation}
Note that $\alpha_{0}$ and $\alpha_{m}$ are special vertices.

\paragraph{(b) even case.} We write $n = 2m$. Following the same
procedure as in the odd case, the root system
$(X^{*}, X_{*}, \Phi, \Phi^{\vee})$ can be computed similarly.

Similarly, $X_{*}(T) = \mathbb{Z}^{n} \times \mathbb{Z}$, so
$X_{*}(T)_{\Gamma}$ can be identified with $\mathbb{Z}^{m} \times
\mathbb{Z}$ under the natural projection $X_{*}(T) \to
X_{*}(T)_{\Gamma}$. Then $X_{*}(S)$ consists of those $(x_{1}, \ldots,
x_{n}; y) \in X_{*}(T)$ satisfying $x_{1} + x_{n} = \dots = x_{m} +
x_{m+1} = y$, we identify $X_{*}(S) \otimes \mathbb{R} \cong
X_{*}(T)_{\Gamma} \otimes \mathbb{R} = X_{*}$. Furthermore, $X^{*}(S)$
can be identified with $\mathbb{Z}^{m} \times \mathbb{Z}$ under the
natural projection $X^{*}(T) \to X^{*}(T)_{\Gamma} / \text{torsion} =
X^{*}(S)$. So the relative roots are
\begin{equation}
  \Phi =
  \left\{
    \begin{array}{ll}
      \pm 2 \epsilon_{i}, & 1 \leq i \leq m ,\\
      \pm \epsilon_{i} \pm \epsilon_{j}, & i \neq j
    \end{array}
  \right\}.
\end{equation}
Then the affine roots are
\begin{equation}
  \Phi_{a} =
  \left\{
    \begin{array}{ll}
      \pm 2 \epsilon_{i} + \mathbb{Z}, & 1 \leq i \leq m, \\
      \pm \epsilon_{i} \pm \epsilon_{j} + \frac{1}{2} \mathbb{Z}, & i \neq j
    \end{array}
  \right\}.
\end{equation}
The affine hyperplanes can be viewed as zero loci of the affine
functions
\begin{equation}
  \left\{
     \begin{array}{ll}
       \pm \epsilon_{i} + \frac{1}{2} \mathbb{Z}, & 1 \leq i \leq m ,\\
       \pm \epsilon_{i} \pm \epsilon_{j} + \frac{1}{2} \mathbb{Z}, & i \neq j,
     \end{array}
  \right.
\end{equation}
which are, as affine root hyperplanes, of type $B_{m}$. So the Weyl
group is $W_{0} = \mathfrak{S}_{m} \rtimes \{ \pm 1 \}^{m}$, the
Iwahori-Weyl group is $\tilde{W} = X_{*}(T)_{\Gamma} \rtimes W_{0}$,
and the affine Weyl group is $W_{a} = X_{*}(T^{\mathrm{sc}})_{\Gamma}
\rtimes W_{0}$, $\pi_{1}(G)_{\Gamma} = \mathbb{Z}/2\mathbb{Z} \times
\mathbb{Z}$.

We choose a basis of $\Phi_{a}$
\begin{equation}
  \left\{
    \begin{array}{ll}
      \alpha_{i} = \epsilon_{i} - \epsilon_{i+1}, & 1 \leq i \leq m-1,
      \\
      \alpha_{m} = 2 \epsilon_{m}, & \\
      \alpha_{0} = \epsilon_{1} + \epsilon_{2} - \frac{1}{2} ,& \\
    \end{array}
  \right.
\end{equation}
then we get the local Dynkin diagram of type $B\text{-}C_{m}$.
\begin{equation}
  \label{eq:245}
  \xymatrix{
    \underset{\alpha_{0}}{\circ} \ar@{-}[rd] & & & & & & & \\
    & \underset{\alpha_{2}}{\circ} & \underset{\alpha_{3}}{\circ} \ar@{-}[l] &
    {\circ} \ar@{-}[l] &  \ar@{.}[l] &{\circ} \ar@{.}[l] &{\circ} \ar@{-}[l] &
    \underset{\alpha_{m}}{\circ} \ar@{=}[l]
    |-{\SelectTips{eu}{11}\object@\dir{>}} \\
   \underset{\alpha_{1}}{\circ} \ar@{-}[ur] & & & & & & &
  }
\end{equation}
Note that $\alpha_{0}$ and $\alpha_{1}$ are special vertices.

\subsubsection{$\mu$-admissible set}
\label{sec:bruh-decomp-mu}

Let $\mu \in X_{*}(T)$ be a minuscule cocharacter, $\lambda$ its image
in $X_{*}(T)_{\Gamma}$, in \cite{MR2141705} the \emph{admissible
  subset} of $\tilde{W}$ is defined as
\begin{equation}
  \Adm(\mu) = \{ w \in \tilde{W} : w \leq t^{w_{0}(\lambda)} \text{ for
    some } w_{0} \in W_{0} \}.
\end{equation}
In the spirit of the Bruhat decomposition, we are interested in the image,
denoted by $\Adm^{0}(\mu)$, of $W_{0}\cdot\Adm(\mu)$ in $W_{0} \backslash
\tilde{W} / W_{0}$. Note that all elements in $\Adm(\mu)$ have the same image in
$\pi_{1}(G)_{\Gamma}$. Because once a special vertex is chosen, we may write $
\tilde{W} = X_{*}(T)_{\Gamma}\rtimes W_{0}$, $\Adm^{0}(\mu)$ is completely
determined by the dominance order on $X_{*}(T)_{\Gamma}$ induced by the Bruhat
order on $\tilde{W}$.

From now on, $\mu = (1,(0)^{n-1};1) \in X_{*}(T) = \mathbb{Z}^{n}
\times \mathbb{Z}$; for $s = 0,1$, $\lambda_{s} = (1^{s},0^{m-s};1)
\in X_{*}(T)_{\Gamma} = \mathbb{Z}^{m} \times \mathbb{Z}$ in both odd
and even cases. Then, as in \cite[2.4.1 \& 2.4.2]{MR2516305},
\begin{equation}
  \label{eq:6}
  \Adm^{0}(\mu) =
  \left\{
    \begin{array}{ll}
      \{\lambda_{1}, \lambda_{0} \} & n \text{ odd,} \\
      \{\lambda_{1}\} & n \text{ even.}
    \end{array}
  \right.
\end{equation}
For convenience of computation, we choose representative(s) $\mu_{1}$
(and $\mu_{0}$ in the odd case) of $\Adm^{0}(\mu)$ in $T(L)$ under the
Kottwitz map of $T$ as follows
\begin{equation}
  \label{eq:9}
  \mu_{1} = \mathrm{diag}(\pi^{2}, \pi, \ldots, \pi, -1)
\end{equation}
in both odd and even case, and
\begin{equation}
  \label{eq:12}
  \mu_{0} = \mathrm{diag}(\pi, \ldots, \pi)
\end{equation}
in odd case. Then, if $K$ is the special parahoric subgroup of $G(L)$
corresponding to the $0$-th vertex of the local Dynkin
diagram~\eqref{eq:244} and \eqref{eq:245},
\begin{equation}
\label{eq:13}
 \bigcup_{w \in \Adm(\mu)} K w K =
 \left\{
   \begin{array}{ll}
     K \mu_{1} K \cup K \mu_{0} K & \text{ odd case,} \\
     K \mu_{1} K & \text{ even case.}
   \end{array}
 \right.
\end{equation}

\subsubsection{Lattice models for Bruhat-Tits buildings and parahoric subgroups}
\label{sec:lattice-model-bruhat}

Recall that $G = \GU(V, \phi)$, now we describe parahoric subgroups of
$G(L)$ in terms of lattices, following \cite{MR2435422}
\cite{MR2516305}.

For $i = 0, \ldots, n-1$, let
\begin{equation}
\label{eq:14}
\Lambda_{i} = \mathrm{span}_{\mathcal{O}_{L}} \{\pi^{-1} e_{1},
\ldots, \pi^{-1} e_{i}, e_{i+1}, \ldots, e_{n} \}.
\end{equation}
More generally, for $j = k n + i$,
$\Lambda_{j} \coloneqq \pi^{-k} \Lambda_{i}$. Let $\mathcal{L}_{I}$ be
the lattice chain $\{\Lambda_{j} : j \in n \mathbb{Z} \pm I\}$ for any
non-empty subset $I \subset \{0,1, \ldots, m\}$. For simplicity, we
write $\mathcal{L}_{i} \coloneqq \mathcal{L}_{\{i\}}$. Note that for
each minimal lattice chain $\mathcal{L}_{i}$, there exists a unique
lattice $M \in \mathcal{L}_{i}$ such that
$M \subset M^{\vee} \subset \pi^{-1} M$, such $M$ is called the
\emph{standard representative} of $\mathcal{L}_{i}$ (see
\cite[6.1]{MR1895706}). It is easy to see that $\Lambda_{i}^{\vee}$ is
the standard representative of $\mathcal{L}_{i}$. Let $P_{I}$ be the
stabilizer of of $\mathcal{L}_{I}$.

\paragraph{(a) odd case.}

In this case, the Kottwitz map is given by
\begin{align}
\label{eq:15}
\kappa_{G} \colon G(L) & \to \pi_{1}(G)_{\Gamma} = \mathbb{Z}, \\
g & \mapsto \val(c(g)). \nonumber
\end{align}
It is easy to see that each element $g \in P_{I}$ has trivial Kottwitz
invariant, then, as described in \cite[1.2.3.(a)]{MR2516305}, the
subgroup $P_{I}$ is a parahoric subgroup of $G(L)$ and each parahoric
subgroup of $G(L)$ is conjugate to $P_{I}$ for some $I$. Note that for
maximal parahoric subgroups, we have
\begin{equation}
  \label{eq:18}
  P_{\{i\}} = \Stab_{G(L)}(M \subset M^{\vee} \subset \pi^{-1} M),
\end{equation}
where $M$ is the standard representative of $\mathcal{L}_{i}$.

\begin{remark}
  \label{rmk:20}
  The maximal parahoric subgroup $P_{\{ i \}}$ for some
  $i \in \{ 0, 1, 2, \ldots, m \}$ corresponds to the $(m-i)$-th
  vertex of the local Dynkin diagram~\eqref{eq:244}, in particular the
  special parahoric subgroup $P_{ \{ m \}}$ corresponds to the $0$-th
  vertex. 
\end{remark}

\paragraph{(b) even case.}

In this case, the Kottwitz map is given by
\begin{align}
  \label{eq:19}
  \kappa_{G} \colon G(L) & \to \pi_{1}(G)_{\Gamma} = \mathbb{Z} \times
                           \{ \pm 1 \},\\
  g & \mapsto (\val(c(g)), (-1)^{\val(b)}), \nonumber
\end{align}
where $b \in \breve{F}^{\times}$ such that $b/\bar{b} = \det(g)\cdot
c(g)^{-m}$ by Hilbert's Satz 90.

Let $P_{I}^{0} \coloneqq P_{I} \cap \ker(\kappa_{G})$, then, as
described in \cite[1.2.3(b)]{MR2516305}, the subgroup $P_{I}^{0}$ is a
parahoric subgroup of $G(L)$, and each parahoric subgroup of $G(L)$ is
conjugate to $P_{I}^{0}$ for a unique subset $I$ satisfying that if
$m - 1 \in I$, then $m \in I$. Note that if $m \in I$, then
$P_{I}^{0} = P_{I}$.

\begin{remark}
  \label{rmk:21}
  Similar to the odd case, the special parahoric subgroup
  $P_{\{ m \}}$ corresponds to the $0$-th vertex.
\end{remark}

\section{Rapoport-Zink spaces and affine Deligne-Lusztig varieties}
\label{sec:rapoport-zink-spaces}

\subsection{Unitary $p$-divisible groups}
\label{sec:unitary-p-divisible}

Let $\Nilp_{\mathcal{O}_{\breve{F}}}$ be the category of
$\mathcal{O}_{\breve{F}}$-schemes $S$ such that $\pi$ is locally
nilpotent on $S$. For $S \in \Nilp_{\mathcal{O}_{\breve{F}}}$, a
\emph{unitary $p$-divisible group of signature $(1,n-1)$} over $S$,
following \cite[3.1]{RSZ}, consists of the following data:
\begin{enumerate}
\item a $p$-divisible group $X$ over $S$,
\item an $\mathcal{O}_{F}$-action $\iota_{X} \colon \mathcal{O}_{F}
  \to \End_{S}(X),$
\item a polarization $\lambda_{X} \colon X \to X^{\vee}$ such that the
  Rosati involution on $\End_{S}(X)$ attached to $\lambda_{X}$ induces
  the non-trivial automorphism on $\mathcal{O}_{F}$ over
  $\mathbb{Q}_{p}$,
\end{enumerate}
satisfying the following conditions:
\begin{enumerate}
\item \emph{Kottwitz condition}:
  \begin{equation}
    \mathrm{charpol}(\iota_{X}(\pi)|\Lie(X)) = (T-\pi)(T+\pi)^{n-1}
    \in \mathcal{O}_{S}[T],
  \end{equation}
\item \emph{Wedge condition}:
  \begin{align}
    \bigwedge^{n} (\iota(\pi) - \pi | \Lie(X)) & =0 \label{eq:22} ,\\
    \bigwedge^{2} (\iota(\pi) + \pi | \Lie(X)) & =0 \text{ if } n \geq
                                                 3 ,\label{eq:23}
  \end{align}
\item when $n$ is even, the extra \emph{Spin condition} is assumed:
  $\iota_{X}(\pi) | \Lie(X_{s})$ non-vanishing for any $s \in S$,
\item \emph{Periodicity condition}: if $n$ is even,
  $\ker(\lambda_{X}) = X[\iota_{X}(\pi)]$; if $n$ is odd,
  $\ker(\lambda_{X}) \subset X[\iota_{X}(\pi)]$ is of height $n-1$.
\end{enumerate}

\begin{remark}
  \label{rmk:22}
  Our definition of unitary $p$-divisible groups is slightly different
  from the one in \cite[3.1]{RSZ}. In our context, the Spin condition
  is not assumed in the odd case, because later we will see that the
  corresponding Rapoport-Zink space is has the same underlying
  topological space with the honest Rapoport-Zink space (see
  Proposition~\ref{pr:9}) which is enough for our purposes because in
  the $p$-adic uniformization theorem (see
  Theorem~\ref{thr:uniformization}), the underlying reduced scheme
  structure is required.
\end{remark}

\subsection{Moduli space of $p$-divisible groups}
\label{sec:moduli-space-p}

From now on, the sign $\mathbb{F}$ denotes the algebraic closure
$\overline{\mathbb{F}}_{p}$. To define the \emph{Rapoport-Zink space},
we fix a supersingular unitary $p$-divisible group $(\mathbb{X},
\iota_{\mathbb{X}}, \lambda_{\mathbb{X}})$ of signature $(1,n-1)$ over
$\mathbb{F}$ as the framing object henceforth. Note that
\cite[Proposition 3.1]{RSZ} shows that such a framing object exists
and is unique up to a quasi-isogeny.

Let $\mathcal{M}^{\mathrm{naive}}$ be the \emph{naive Rapoport-Zink
  space}, the formal scheme $\mathcal{M}^{\mathrm{naive}}$ is formally
locally of finite type over $\Spf \mathcal{O}_{\breve{F}}$ (cf.
\cite[Theorem 3.25]{MR1393439}). Unfortunately,
$\mathcal{M}^{\mathrm{naive}}$ is not flat over
$\mathcal{O}_{\breve{F}}$ (cf. \cite[Proposition 3.8]{MR1752014})
because its corresponding local model is not flat. Let
$\mathbf{M}^{\mathrm{loc}}$ be the honest local model, $\mathcal{M}$
its corresponding Rapoport-Zink space via the local model diagram,
then $\mathcal{M}$ is called the \emph{honest Rapoport-Zink space}.
However, it's not clear whether $\mathcal{M}$ has a moduli
description. Another way to define Rapoport-Zink spaces is to add some
extra conditions on the $p$-divisible groups and get a moduli space of
$p$-divisible groups with extra conditions.

Now we associate to
$(\mathbb{X}, \iota_{\mathbb{X}}, \lambda_{\mathbb{X}})$ a set-valued
functor $\mathcal{M}^{e}$ on the category
$\Nilp_{\mathcal{O}_{\breve{F}}}$. The superscript $e$ stands for
``exotic''.

\begin{definition}
  \label{df:14}
  For any $S \in \Nilp_{\mathbb{F}}$, $\mathcal{M}^{e}(S)$ is the set
  of isomorphism classes of $(X, \iota_{X}, \lambda_{X}, \rho_{X})$,
  where
  \begin{itemize}
  \item $(X, \iota_{X}, \lambda_{X})$ is a unitary $p$-divisible group
    of signature $(1,n-1)$ over $S$;
  \item
    $\rho_{X} \colon X \times_{S} \overline{S} \to \mathbb{X}
    \times_{\mathbb{F}} \overline{S}$
    is an $\mathcal{O}_{F}$-linear quasi-isogeny (of any height) such
    that $\rho^{*}(\lambda_{\mathbb{X}})$ and $\lambda_{X}$ differ
    locally on $\overline{S}$ by a scalar in
    $\mathbb{Q}_{p}^{\times}$.
  \end{itemize}
  Two quadruples $(X, \iota_{X}, \lambda_{X}, \rho_{X})$ and $(Y,
  \iota_{Y}, \lambda_{Y}, \rho_{Y})$ are isomorphic if there exists an
  $\mathcal{O}_{F}$-linear isomorphism of $p$-divisible groups $\alpha
  \colon X \to Y$ such that $\rho_{Y} \circ \alpha = \rho_{X}$ and
  $\alpha^{*}(\lambda_{Y})$ and $\lambda_{X}$ differ locally on
  $\overline{S}$ by a scalar in $\mathbb{Q}_{p}^{\times}$.
\end{definition}

\begin{prop}[Smithling]
  \label{pr:9}
  The functor $\mathcal{M}^{e}$ is represented by a separated formal
  scheme over $\Spf(\mathcal{O}_{\breve{F}})$, which is locally
  formally of finite type, and of relative formal dimension $n-1$ over
  $\mathcal{O}_{\breve{F}}$, and has the same underlying topological
  space with $\mathcal{M}$. Furthermore, if $n$ is even,
  $\mathcal{M}^{e}$ is flat over $\mathcal{O}_{\breve{F}}$.
\end{prop}

\begin{proof}
  Let $\mathbf{M}^{\mathrm{naive}}$, $\mathbf{M}^{e}$,
  $\mathbf{M}^{\mathrm{loc}}$ be the corresponding local models of
  $\mathcal{M}^{\mathrm{naive}}$, $\mathcal{M}^{e}$, $\mathcal{M}$. By
  \cite[Corollary 5.6.3]{MR2764885} and \cite[Theorem 1.3]{MR3177281},
  $\mathbf{M}^{e}$ is topologically flat, hence $\mathcal{M}^{e}$ has
  the same underlying topological space with $\mathcal{M}$. The
  flatness of $\mathbf{M}^{e}$ in the even case follows from
  \cite[Proposition 3.10]{RSZ}, so by the local model diagram,
  $\mathcal{M}^{e}$ is flat, i.e. $\mathcal{M} = \mathcal{M}^{e}$.
\end{proof}

Note that the framing object $\mathbb{X}$ can be defined over
$\mathbb{F}_{p}$. Let $(\mathbb{M}_{0}, \mathcal{F}_{0},
\mathcal{V}_{0})$ be the Dieudonn\'{e} module of $\mathbb{X}$ over
$\mathbb{F}_{p}$, which is a free $\mathbb{Z}_{p}$-module of rank $2
n$, $N_{0}$ its isocrystal with Frobenius $\mathcal{F}_{0}$ and
Verschiebung $\mathcal{V}_{0}$. The action of $F$ on $\mathbb{X}$
makes $N_{0}$ an $F$-vector space. The polarization
$\lambda_{\mathbb{X}}$ induces an alternating
$\mathbb{Q}_{p}$-bilinear non-degenerate form on $N_{0}$
\begin{equation}
\label{eq:25}
\langle \thinspace, \thinspace \rangle \colon N_{0} \times N_{0}
\to \mathbb{Q}_{p},
\end{equation}
such that
\begin{equation}
  \label{eq:26}
\langle \pi x, y\rangle = \langle x, \bar{\pi} y\rangle.
\end{equation}
This is equivalent to giving a
hermitian form $\varphi$ on $N_{0}$ such that
\begin{equation}
  \label{eq:27}
\langle x, y\rangle = \frac{1}{2} \Tr_{{F}/\mathbb{Q}_{p}}
(\pi^{-1} \varphi(x,y)).
\end{equation}
Then the periodicity condition for $\mathbb{X}$ means that
$\mathbb{M}_{0}$ is a nearly $\pi$-modular lattice in the odd case,
and a $\pi$-modular lattice in the even case. Note that by
\eqref{eq:27}, the dual of $\mathbb{M}_{0}$ with respect to $\varphi$
is the same as the dual with respect to $\langle \thinspace,
\thinspace\rangle$. So by Lemma~\ref{split-lemma}, we can choose a
${F}$-basis $\{e_{1}, \ldots, e_{n}\}$ of $N_{0}$ such that
$\varphi(e_{i},e_{j}) = \delta_{i,n+1-j}$. We borrow the notation from
\eqref{eq:14} to denote the ``standard'' lattices, so $\mathbb{M}_{0}
= \Lambda_{m}^{\vee}$. Now Let $(\mathbb{M}, \mathcal{F}, \mathcal{V})
= (\mathbb{M}_{0}, \mathcal{F}_{0}, \mathcal{V}_{0}) \otimes L$, $N =
N_{0} \otimes L$. Then for any $x, y \in N$, we have
\begin{equation}
  \label{eq:20}
  \langle \mathcal{F} x, y\rangle =
  \langle x, \mathcal{V} y\rangle^{\sigma}.
\end{equation}

Let $\mathcal{N}, \mathcal{N}^{e}, \mathcal{N}^{\mathrm{naive}}$ be
the special fibers of $\mathcal{M}, \mathcal{M}^{e},
\mathcal{M}^{\mathrm{naive}}$ respectively, we are interested in the
geometric structure of $\mathcal{N}$. Because $\mathcal{M}^{e}$ has
the same underlying topological space with $\mathcal{M}$, we have
$\mathcal{N}_{\mathrm{red}} = \mathcal{N}^{e}_{\mathrm{red}}$. The
$\mathbb{F}$-valued points of $\mathcal{N}^{e}$ have a simple
description: the Kottwitz condition means that $X$ is of dimension $n$
and of height $2 n$; the Wedge condition \eqref{eq:22} is trivial,
\eqref{eq:23} means
\begin{equation}
\label{eq:24}
\bigwedge^{2}(\iota(\pi) | \Lie(X)) = 0,
\end{equation}
i.e. the rank of the operator $\iota(\pi) | \Lie(X)$ is less than or
equal to $1$; in the even case the spin condition means the rank of
the operator $\iota(\pi)| \Lie(X)$ is $1$.

\begin{prop}
  \label{pr:2}
  Via Dieudonn\'{e} theory,
  $\mathcal{N}(\mathbb{F}) = \mathcal{N}^{e}(\mathbb{F})$ can be
  identified with the set of $\mathcal{O}_{\breve{F}}$-lattices $M$ in
  $N$ satisfying the following conditions:
  \begin{enumerate}
  \item $M$ is stable under $\mathcal{F}$ and $\mathcal{V}$;
  \item
    $M \stackrel{n-1}{\subset} p^{h} M^{\vee} \stackrel{1}{\subset}
    \pi^{-1} M$
    if $n$ is odd, and $p^{h} M^{\vee} = \pi^{-1} M$ if $n$ is even
    for some $h \in \mathbb{Z}$;
  \item $p M \stackrel{n}{\subset} \mathcal{V} M \stackrel{n}{\subset}
    M$; \label{item:3}
  \item $\mathcal{V} M \stackrel{\leq 1}{\subset} \mathcal{V} M + \pi
    M$; \label{item:4}
  \item if $n$ is even,
    $\mathcal{V} M \stackrel{1}{\subset} \mathcal{V} M + \pi
    M$. \label{item:5}
  \end{enumerate}
\end{prop}

\begin{proof}
  Via Dieudonn\'{e} theory, condition~\ref{item:3} is just the
  Kottwitz condition, condition~\ref{item:4} is the wedge condition
  and condition~\ref{item:5} is the extra Spin condition.
\end{proof}

\subsection{Local PEL datum}
\label{sec:local-pel-datum}

Let $G$ be the algebraic group $\GU(N_{0},\varphi)$ which is defined
over $\mathbb{Q}_{p}$. We write $\mathcal{F} = b \cdot \id_{F} \otimes
\sigma$ for some $b \in \GL_{\breve{F}}(N)$ in terms of the basis
$\{e_{1}, \ldots, e_{n}\}$, by \eqref{eq:20} and \eqref{eq:27}, we
have
\begin{equation}
  \label{eq:28}
  \varphi(\mathcal{F}x, \mathcal{F}y) = p \cdot \varphi(x^{\sigma}
  ,y^{\sigma})
\end{equation}
which implies that $b \in G(L)$ with $\val(c(b)) = 1$. Let
$[b] \in B(G)$ be the $G(L)$-conjugacy classes of $b$, i.e. the set
$\{g^{-1} b \sigma(g) : g \in G(L)\}$. We use the notation from
Section~\ref{sec:group-data}, i.e. $S$ is the maximal $L$-split torus
of $G$, $T$ is the centralizer of $S$ with
$X_{*}(T) \simeq \mathbb{Z}^{n} \times \mathbb{Z}$, and
$\mu \in X_{*}(T)$ is the geometric minuscule cocharacter
$(1, 0^{n-1}; 1)$. By the assumption of supersingularity and Kottwitz
condition on $\mathbb{X}$, we have
\begin{equation}
  \label{eq:36}
  [b] \in B(G,\{ \mu \})_{b},
\end{equation}
where $\{ \mu \}$ is the geometric conjugacy classes of $\mu$ and the
well-known set $B(G, \{ \mu \})$ is the subset of $B(G)$ consisting of
neutral acceptable elements (cf. \cite[6.2]{MR1485921}
\cite[Definition 2.3]{MR3271247}).

In summary, $({F}, F, N_{0}, \varphi, \langle \thinspace, \thinspace
\rangle, \thinspace {\bar{}} \thinspace, \{ \mu \}, [b], \pi,
\mathbb{M}_{0})$ forms a \emph{simple integral Rapoport-Zink
  PEL-datum} in the sense of \cite[4.1]{MR3271247} (cf.
\cite[Definition 3.18]{MR1393439}).

Another important group is the algebraic group $J$ consisting of
automorphisms of the unitary isocrystal $N$, i.e.
\begin{equation}
\label{eq:37}
J(R) = \left\{ g \in \GL_{\breve{F} \otimes R} (N
  \otimes_{\mathbb{Q}_{p}}R)
  \left|
    \begin{array}{l}
      g \mathcal{F} = \mathcal{F} g, \varphi(g x, g y) = c(g) \varphi(x,y) \\
      \text{for some } c(g) \in
      (L \otimes R)^{\times}
    \end{array}
  \right.
\right\}   
\end{equation}
for any $\mathbb{Q}_{p}$-algebra $R$. The group $J$ acts on
$\mathcal{N}^{e}$: for $g \in J$, the action is given by sending
$(X, \iota_{X}, \lambda_{X}, \rho_{X}) \in \mathcal{N}^{e}$ to
$(X, \iota_{X}, \lambda_{X}, g \circ \lambda_{X})$. By
\cite[5.2]{MR809866}, $J$ is an inner form of $G$ because $[b]$ is
basic.

The group $J$ is closely related to a hermitian space, namely $C$,
with respect to $F/\mathbb{Q}_{p}$ as discussed in
\cite{MR3175176}. Recall that $\pi^{2} = \varpi = \epsilon p$, let
$\eta, \delta \in \mathcal{O}_{L}^{\times}$ such that
$ \eta^{2} = \epsilon^{-1}$ and $\delta^{\sigma} = -\delta$
respectively. Then all slopes of the $\id \otimes \sigma$-linear
operator $\chi \coloneqq \eta \pi \mathcal{V}^{-1} \colon N \to N$ are
zero. Let $C$ be the set of points in $N$ fixed by $\chi$, then $C$ is
a vector space over $F$ and the isomorphism
\begin{equation}
\label{eq:38}
C \otimes_{\mathbb{Q}_{p}} L \simeq N
\end{equation}
identifies $\id_{C} \otimes \sigma$ with $\chi$. Let $\psi(x,y) \coloneqq
\delta \varphi(x,y)$ for $x,y \in C$, then by \eqref{eq:28}, we have
\begin{equation}
\label{eq:40}
\psi(x,y) = \psi(x,y)^{\sigma}.
\end{equation}
So $\psi$ takes values in $F$ and hence $(C, \psi)$ becomes a
hermitian space with respect to $F/\mathbb{Q}_{p}$. By \cite[Lemma
2.3]{MR3175176}, The group $J$ is isomorphic to the general unitary
group $\GU(C,\psi)$. By \cite[Lemma 3.3]{RSZ}, the hermitian space
$(C, \psi)$ is split if $n$ is odd, non-split if $n$ is even.

\begin{remark}
  \label{rmk:3} In \cite{smithling-moduli}, when $n$ is odd, the
  moduli description of $\mathbf{M}^{\mathrm{loc}}$, hence of
  $\mathcal{N}$, is formulated by proposing a further refinement of
  the spin condition, which is unfortunately very complicated. For the
  purposes of studying basic loci of Shimura varieties in this paper,
  for us it is enough to work with $\mathcal{N}^{e}$ since
  $\mathcal{N}^{e}_{\mathrm{red}} = \mathcal{N}_{\mathrm{red}}$.
\end{remark}

\subsection{Kottwitz invariants of quasi-isogenies}
\label{sec:kottw-invar-quasi}

In this section, we will define a morphism
\begin{equation}
\label{eq:39}
\kappa \colon \mathcal{N}^{e} \to \pi_{1}(G)_{\Gamma}.
\end{equation}

Let $X \in \mathcal{N}^{e}(\mathbb{F})$ be a unitary $p$-divisible
group with a quasi-isogeny $\rho \colon X \to \mathbb{X}$, let $M$ be
its corresponding Dieudonn\'{e} lattice in $N$, recall that the height
of $\rho$ is defined as
\begin{equation}
\label{eq:30}
\height(\rho) \coloneqq \height(p^{s} \rho) - \height(p^{s}),
\end{equation}
where $s$ is an integer such that $p^{s} \rho$ is an honest
isogeny.

\begin{prop}
  \label{pr:3}
  If $M$ satisfies $M \subset p^{h} M^{\vee} \subset \pi^{-1}M$ for
  some integer $h$, then $\height(\rho) = n h$.
\end{prop}

\begin{proof}
  The proof is pretty easy, so we leave it to the reader.
\end{proof}

For a unitary $p$-divisible group $X \in \mathcal{N}^{e}(S)$ with a
quasi-isogeny $\rho$, the height is locally constant on $S$, so we get
a morphism
\begin{align}
  \label{eq:33}
  \kappa_{1} \colon \mathcal{N}^{e} & \longrightarrow \mathbb{Z},\\
  \label{eq:34}
  (X,\rho) & \longmapsto \frac{1}{n} \height(\rho). \nonumber
\end{align}
By abuse of notation, we denote by $\kappa_{1}$ the composite morphism
$\mathcal{N} \subset \mathcal{N}^{e} \to \mathbb{Z}$.  Let
$\mathcal{N}_{h}$ (resp. $\mathcal{N}^{e}_{h}$) be the fiber
$\kappa_{1}^{-1}(h)$ for $h\in \mathbb{Z}$, then $\mathcal{N}_{h}$
(resp. $\mathcal{N}^{e}_{h}$) is an open and closed subscheme of
$\mathcal{N}$ (resp. $\mathcal{N}^{e}$), we have a decomposition
\begin{equation}
  \label{eq:35}
  \mathcal{N} = \coprod_{h \in \mathbb{Z}} \mathcal{N}_{h} \thinspace
  \text{ and  } \thinspace \mathcal{N}^{e} =
  \coprod_{h \in \mathbb{Z}} \mathcal{N}^{e}_{h}.
\end{equation}

In the even case, $\mathcal{N} = \mathcal{N}^{e}$, there is an extra
invariant of quasi-isogenies, which has been discussed in \cite[Lemma
3.2]{RSZ}. Let $(\tilde{X}, \rho_{1}, \rho_{2})$ be the \emph{minimal
  cover} of the quasi-isogeny $\rho$ in the following sense:
$\tilde{X}$ is a $p$-divisible group together with isogenies
$\rho_{1}, \rho_{2}$ making $\rho \circ \rho_{1} = \rho_{2}$, such
that for any $p$-divisible group $Y$ with isogenies $\alpha_{1},
\alpha_{2}$ satisfying $\rho \circ \alpha_{1} = \alpha_{2}$, there
exists a unique isogeny $\beta \colon Y \to \tilde{X}$ making the
following diagram commutative
\begin{equation}
\label{eq:45}
\xymatrix{
  & Y \ar@/_/[ddl]_{\alpha_{1}} \ar@/^/[ddr]^{\alpha_{2}} \ar@{-->}[d] & \\
  & \tilde{X} \ar[dl]^{\rho_{1}} \ar[dr]_{\rho_{2}} & \\
  X \ar[rr]_{\rho} & & \mathbb{X}.
}
\end{equation}
Note that, via Dieudonn\'{e} theory, $\tilde{X}$ corresponds to the
lattice $M \cap \mathbb{M}$. By \eqref{eq:30} and
Proposition~\ref{pr:3}, we have
\begin{equation}
\label{eq:46}
n h = \height(\rho) = \height(\rho_{2}) - \height(\rho_{1}).
\end{equation}
Because $n$ is even,
\begin{equation}
\label{eq:47}
\height(\rho_{2}) \equiv \height(\rho_{1}) \mod{2}.
\end{equation}
Hence we get a morphism
\begin{align}
  \label{eq:48}
  (\kappa_{1}, \kappa_{2}) \colon \mathcal{N} & \longrightarrow
                                                \mathbb{Z} \times \mathbb{Z} / 2
                                                \mathbb{Z} ,\\
  (X, \rho) & \longmapsto  (\frac{1}{n} \height(\rho), \thinspace
              \height(\rho_{1}) \mod{2}). \nonumber
\end{align}

In summary, we have the \emph{Kottwitz morphism}
\begin{equation}
\label{eq:49}
\kappa \colon \mathcal{N} \longrightarrow \pi_{1}(G)_{\Gamma},
\end{equation}
when $n$ is odd, $\kappa = \kappa_{1}$, when $n$ is even, $\kappa =
(\kappa_{1}, \kappa_{2})$. We have the decomposition
\begin{equation}
\label{eq:50}
\mathcal{N} = \coprod_{\kappa \in \pi_{1}(G)_{\Gamma}} \mathcal{N}_{(\kappa)},
\end{equation}
where $\mathcal{N}_{(\kappa)}$ consists of those quasi-isogenies with
Kottwitz invariants $\kappa \in \pi_{1}(G)_{\Gamma}$.
For any $\kappa, \kappa'$, let $g \in J$ such that $\kappa(g) =
\kappa' -\kappa$, then $g$ defines an isomorphism
\begin{align}
\label{eq:41}
  \mathcal{N}_{(\kappa)} & \longrightarrow \mathcal{N}_{(\kappa')} ,\\
  (X, \rho) & \longmapsto (X, g \circ \rho). \nonumber
\end{align}
Via Dieudonn\'{e} theory, we have, in the odd case
\begin{equation}
\label{eq:42}
\mathcal{N}_{(\kappa)}(\mathbb{F}) = \{ M \in \mathcal{N}(\mathbb{F}) :
M \subset p^{\kappa} M^{\vee} \subset \pi^{-1} M \}.
\end{equation}
In the even case for $\kappa = (\kappa_{1}, \kappa_{2})$
\begin{equation}
\label{eq:43}
\mathcal{N}_{(\kappa)} = \{M \in \mathcal{N}(\mathbb{F})
\left|
  \begin{array}{c}
    p^{\kappa_{1}} M^{\vee} = \pi^{-1} M, \\
    \dim_{\mathbb{F}} (M + \mathbb{M} / M) \equiv \kappa_{2} \mod{2}
  \end{array}
\right\}.
\end{equation}

\subsection{Affine Deligne-Lusztig varieties}
\label{sec:affine-delig-luszt}

Recall that for the local PEL-datum $({F}, F, N_{0}, \varphi, \langle
\thinspace, \thinspace \rangle, \thinspace {\bar{}} \thinspace, \{ \mu
\}, [b], \pi, \mathbb{M}_{0})$ in Section \ref{sec:local-pel-datum},
we may associate to it the \emph{generalized affine Deligne-Lusztig
  variety} (cf. \cite[Definition 4.1]{MR2141705})
\begin{equation}
  \label{eq:51}
  X(\mu, b)_{K} \coloneqq \{ g \in G(L)/K : g^{-1} b \sigma(g) \in
  \bigcup_{w \in \Adm(\mu)} K w K \},
\end{equation}
where $K = \Stab_{G(L)}(\mathbb{M} \subset \mathbb{M}^{\vee} \subset
\pi^{-1} \mathbb{M})$ which is the special parahoric subgroup $P_{\{m
  \}}$ of $G(L)$ corresponding to the $0$-th vertex of the local
Dynkin diagram in both odd and even cases because $\mathbb{M} =
\Lambda^{\vee}_{m}$. Note that the group $J$ also acts on $X(\mu,
b)_{K}$ because $J$ is just the $\sigma$-centralizer of $b$ in $G(L)$.

\begin{prop}
  \label{pr:6}
  The map
  \begin{align}
    \label{eq:52}
    \Phi \colon X(\mu, b)_{K} & \longrightarrow \mathcal{N}(\mathbb{F}),
    \\
    g  & \longmapsto g \mathbb{M}, \nonumber
  \end{align}
  is bijective.
\end{prop}

\begin{proof}
  For $g \in X(\mu, b)_{K}$, we need to check that
  $g \mathbb{M} \in \mathcal{N}(\mathbb{F})$, i.e. it satisfies the
  conditions in Proposition~\ref{pr:2}. Recall that we choose
  representative(s) $\mu_{1}$ (and $\mu_{0}$ in odd case) of
  $\Adm^{0}(\mu)$ in $T(L)$ in the subsection
  \ref{sec:bruh-decomp-mu}. The condition $g \mathbb{M}$ is stable under
    $\mathcal{F} = b \cdot \id \otimes \sigma$ is equivalent to the
    condition
    \begin{equation}
      \label{eq:53}
      g^{-1} b \sigma (g) \mathbb{M} \subset \mathbb{M},
    \end{equation}
    so by \eqref{eq:13}, it is enough to check
    $\mu_{1} \mathbb{M} \subset \mathbb{M}$ (and
    $\mu_{0} \mathbb{M} \subset \mathbb{M}$ in odd case). By the
    choice of $\mu_{1}$ (and $\mu_{0}$) in \eqref{eq:9} and
    \eqref{eq:12}, it is easy to see $g \mathbb{M}$ is
    $\mathcal{F}$-stable, and similarly, $\mathcal{V}$-stable. For the
    rest of conditions, we leave them to the reader.

    It's very easy to see that $\Phi$ is injective. For the
    surjectivity of $\Phi$, we use the \emph{G\"{o}rtz local model
      diagram} (cf. \cite[5.2]{MR2602029}):
    \begin{equation}
      \label{eq:81}
      \mathcal{N}^{e}(\mathbb{F}) \xhookrightarrow[ ]{\quad }
      G(L)/K \xleftarrow[]{\;\mathrm{pr} \;}
      G(L) \xrightarrow[ ]{\; \mathrm{pr}^{\sigma} \;}
      G(L)/K \xhookleftarrow[ ]{\quad}
      \mathbf{M}^{e}_{\mathbb{F}}
    \end{equation}
    where $\mathrm{pr}$ is the natural projection,
    $\mathrm{pr}^{\sigma}$ is the composite of the Lang map $g \mapsto
    g^{-1} b \sigma(g)$, with the projection $\mathrm{pr}$. Then
    $\mathrm{pr}^{-1} (\mathcal{N}^{e}(\mathbb{F})) =
    (\mathrm{pr}^{\sigma})^{-1}(\mathbf{M}^{e}_{\mathbb{F}})$. Note
    that
  \begin{equation}
    \mathbf{M}^{e}_{\mathbb{F}} = \bigcup_{w \in \Adm(\mu)} K w K / K
    \label{eq:82}
  \end{equation}
  in $G(L)/K$ by \cite[Corollary 5.6.2]{MR2764885} and \cite[Theorem
  1.4]{MR3177281}. So
  \begin{equation}
    \mathrm{pr}^{-1}(X(\mu, b)_{K}) =
    (\mathrm{pr}^{\sigma})^{-1}
    (\mathbf{M}^{e}_{\mathbb{F}}) \label{eq:83},
  \end{equation}
  and the injectivity of $\Phi$ implies that $\Phi$ is surjective.
\end{proof}

\begin{remark}
  \label{rmk:2}
  For $g \in X(\mu,b)_{K}$, the Kottwitz invariant is well defined by
  the definition of parahoric subgroups and compatible with the
  Kottwitz map for $\mathcal{N}^{e}(\mathbb{F})$ via the map $\Phi$.
  So $X(\mu,b)_{K}$ can be decomposed into a disjoint union of some
  subsets indexed by Kottwitz invariants. In the odd case, for any
  $\kappa \in \pi_{1}(G)_{\Gamma} \simeq \mathbb{Z}$,
  $\mathcal{N}_{(\kappa)} (\mathbb{F})$ can be identified with a
  generalized affine Deligne-Lusztig variety
  $X(\tilde{\mu},\tilde{b})_{K^{\prime}}$ associated to the derived
  group of $G$, i.e. the special unitary group $\SU(V, \varphi)$,
  where $K^{\prime} = K \cap \SU(V, \varphi)$. Because in this case,
  there exists a central element $\zeta$ such that $\mu = \zeta
  \tilde{\mu}$ and $b = \tilde{b} \zeta$. Then the map
  \begin{align}
    \label{eq:59}
    X(\tilde{\mu}, \tilde{b})_{K^{\prime}} & \longrightarrow
     \mathcal{N}_{(1)}(\mathbb{F}) ,\\
    \label{eq:60}
    g & \longmapsto g \zeta \mathbb{M}, \nonumber
  \end{align}
  gives the desired identification. However, this is no longer true in
  the even case, because
  $\mu = (1, 0^{(m-1)}; 1) \in X_{*}(T)_{\Gamma}$ and
  $\tilde{\mu} = (2,0^{(m-1)}) \in X_{*}(T^{\mathrm{sc}})$ differ in a
  non-central element in $\Omega$. We will work with the corresponding
  \emph{semisimple group of adjoint type} $G_{\mathrm{ad}}$, i.e. the
  quotient of $G$ by its center.
\end{remark}

Let $b_{\ad}, \mu_{\ad}, K_{\ad}$ be the images in $G_{\ad}(L)$ of
$b, \mu, K$ respectively. Similarly to \eqref{eq:51}, we define
\begin{equation}
  \label{eq:63}
  X(\mu_{\ad}, b_{\ad})_{K_{\ad}} \coloneqq
  \{ g \in G_{\ad}(L) / K_{\ad} : g^{-1} b_{\ad}
  \sigma(g) \in \bigcup_{w \in \Adm(\mu_{\ad})} K_{\ad}
  w K_{\ad} \},
\end{equation}
where $\Adm(\mu_{ad})$ is the $\mu_{\ad}$-admissible subset of the
Iwahori-Weyl group $\tilde{W}_{\ad}$ of $G_{\ad}$, which is bijective
to $\Adm(\mu)$ under the canonical map
$\tilde{W} \to \tilde{W}_{\ad}$. However, there is no reasonable map
from $X(\mu_{\ad}, b_{\ad})_{K_{\ad}}$ to $\mathcal{N}(\mathbb{F})$,
because for $g \in X(\mu_{\ad}, b_{\ad})_{K_{\ad}}$,
$M \in \mathcal{N}(\mathbb{F})$, the notation $g M$ doesn't make
sense. $g M$ is no longer a lattice, but a homothety class of
lattices. However, by \cite[6.a]{MR2435422}, (see also
\cite[2.2]{MR3377055}) the natural map
\begin{equation}
  \label{eq:84}
  G(L)/K \longrightarrow G_{\ad}(L) / K_{\ad}
\end{equation}
induces a bijection
\begin{equation}
  \label{eq:85}
  (G(L)/K)_{\kappa} \stackrel{\cong}{\longrightarrow}
  (G_{\ad}(L)/K_{\ad})_{\kappa_{\ad}},
\end{equation}
where $\kappa \in \pi_{1}(G)_{\Gamma}$, $\kappa_{\ad}$ is the image of
$\kappa$ in $\pi_{1}(G_{\ad})_{\Gamma}$, the notation
$(\thinspace \cdot \thinspace)_{\kappa}$ stands for the fiber of
corresponding Kottwitz maps. Immediately, we have
\begin{equation}
  \label{eq:86}
  (X(\mu, b)_{K})_{\kappa} \stackrel{\cong}{\longrightarrow}
  (X(\mu_{\ad}, b_{\ad})_{K_{\ad}})_{\kappa_{\ad}}.
\end{equation}
In particular, $X(\mu_{\ad}, b_{\ad})_{K_{\ad}}$ can be identified
with the moduli space of quasi-isogenies of height $0$, i.e.
\begin{align}
  \label{eq:87}
  \Phi_{\mathrm{ad}} \colon X(\mu_{\ad}, b_{\ad})_{K_{\ad}}
  & \stackrel{\cong}{\longrightarrow} \mathcal{N}_{0}(\mathbb{F}), \\
  g K_{\ad} & \longmapsto \pi^{- \val(\dot{g})} \dot{g} \mathbb{M},
              \nonumber
\end{align}
where $\dot{g} K$ is a lifting of $g K_{\ad}$ under the map
\eqref{eq:84}.

Let $G_{\mathrm{ad}}(L)'$ be the subgroup of $G_{\mathrm{ad}}(L)$
generated by all the parahoric subgroups of $G_{\mathrm{ad}}(L)$. Let
\begin{equation}
  \label{eq:116}
  X(\mu_{\mathrm{ad}}, b_{\mathrm{ad}})_{K_{\mathrm{ad}}}' \coloneqq
  \{ g \in G_{\mathrm{ad}}(L)' / K_{\mathrm{ad}} :
  g^{-1} b_{\mathrm{ad}} \sigma(g) \in \bigcup_{w \in \Adm(\mu_{\ad})} K_{\ad}
  w K_{\ad} \}.
\end{equation}
Note that $G_{\mathrm{ad}}(L)' = G_{\mathrm{ad}}(L)_{1}$ by
\cite[Lemma 17]{Haines-Rapoport}, in other words, the kernel of
Kottwitz map is generated by all the parahoric subgroups. When $n$ is
odd, nothing is new because
$X(\mu_{\mathrm{ad}}, b_{\mathrm{ad}})_{K_{\mathrm{ad}}}' =
X(\mu_{\mathrm{ad}}, b_{\mathrm{ad}})_{K_{\mathrm{ad}}}$; when $n$ is
even, the map $\Phi_{\mathrm{ad}}$ in \eqref{eq:87} induces the
following isomorphism:
\begin{equation}
  \label{eq:120}
  \Phi_{\mathrm{ad}} \colon X(\mu_{\mathrm{ad}},
  b_{\mathrm{ad}})_{K_{\mathrm{ad}}}'
  \stackrel{\cong}{\longrightarrow} \mathcal{N}_{(0,0)}(\mathbb{F}).
\end{equation}
From now on, let
\begin{equation}
  \label{eq:121}
  \mathcal{S} \coloneqq 
  \left\{
    \begin{array}{ll}
      \mathcal{N}_{0} & \text{ if } n \text{ is odd,} \\
      \mathcal{N}_{(0,0)} & \text{ if } n \text{ is even.}
    \end{array}
  \right.
\end{equation}

\begin{remark}
  \label{rmk:9}
  When $n$ is even, by the definition of $\mathcal{S}$,
  $\mathcal{S}(\mathbb{F})$ is a single $G_{\ad}(L)'$-orbit of
  $\mathbb{M}$.
\end{remark}

\begin{remark}
  An equivalent way to identify $X(\mu_{\ad}, b_{\ad})_{K_{\ad}}$ with
  a reasonable Rapoport-Zink space is to define the \emph{adjoint
    Rapoport-Zink space}
  $$\mathcal{N}_{\ad} \coloneqq \mathcal{N} / \mathbb{G}_{m}(F).$$
  Because in both odd and even cases, the action of $\pi$ on
  $\mathcal{N}$ via $\iota_{\mathbb{X}} \colon F \to \End(\mathbb{X})$
  gives an isomorphism
  \begin{equation}
    \label{eq:61}
    \mathcal{N}_{h} \longrightarrow \mathcal{N}_{h+1}.
  \end{equation}
  The set $\mathcal{N}_{\ad}(\mathbb{F})$ can be described as a set of
  homothety classes of lattices satisfying the Kottwitz, wedge and the
  extra spin conditions, so that for
  $g \in X(\mu_{\ad}, b_{\ad})_{K_{\ad}}$, the notation $g \mathbb{M}$
  makes sense as a homothety class of lattices.
\end{remark}

\section{Set structure of $\mathcal{N}$}
\label{sec:set-structure}

\subsection{Deligne-Lusztig varieties}
\label{sec:deligne-luszt-vari}

We need some results about classical Deligne-Lusztig varieties. Let
$H_{0}$ be a reductive group over $\mathbb{F}_{q}$. We fix a maximal
torus $T_{0}$ and Borel subgroup $B_{0}$ over $\mathbb{F}_{q}$. Let
$H$ be the reductive group $H_{0} \otimes \overline{\mathbb{F}}_{q}$
over $\overline{\mathbb{F}}_{q}$,
$B \coloneqq B_{0} \otimes \overline{\mathbb{F}}_{q}$, with a
Frobenius action $\sigma$. We fix a $\sigma$-stable maximal torus $T$
and Borel subgroup $B$. Let $W = W_{H}$ be the Weyl group of $H$. The
\emph{(classical) Deligne-Lusztig variety} (cf. \cite[Definition
1.4]{MR0393266}) $X(w)$ is defined as
\begin{equation}
  \label{eq:132}
  X(w) \coloneqq \{ g \in H/B : g^{-1} \sigma(g) \in B w B \},
\end{equation}
for each $w \in W_{H}$. We also say that $g$ and $h$ are in
\emph{relative position} $w$ if $g^{-1} h \in B w B$ for
$g, h \in G/B$ and $w \in W_{H}$.

\begin{prop}[{cf. \cite{MR0393266}, see also \cite[Proposition
    4.4]{MR3013026}}]
  \label{pr:14}
  For $w \in W_{H}$.
  \begin{enumerate}
  \item The Deligne-Lusztig variety $X(w)$ is smooth and of pure
    dimension $\ell(w)$, where $\ell(w)$ is the length of $w$.
  \item The flag variety $H/B$ is the disjoint union of all
    Deligne-Lusztig varieties, indexed by the Weyl group $W_{H}$. The
    closure $\overline{X(w)}$ of $X(w)$ in the flag variety $H/B$ is
    normal, and
    \begin{equation}
      \label{eq:157}
      \overline{X(w)} = \bigcup_{w' \leq w} X(w'),
    \end{equation}
    where $\leq$ denotes the Bruhat order in $W_{H}$. Furthermore, if
    $w$ is a Coxeter element, $\overline{X(w)}$ is smooth.
  \item The Deligne-Lusztig variety $X(w)$ is irreducible if and only
    if $w$ is not contained in any $\sigma$-stable standard parabolic
    subgroup of $W_{H}$.
  \end{enumerate}
\end{prop}

\begin{eg}[The split odd orthogonal group]
  \label{eg:2}
  Let $V$ be an $l$-dimensional vector space over $\mathbb{F}_{q}$,
  where $l = 2 d + 1$ is odd, together with ``the'' split
  non-degenerate symmetric form
  $\langle \thinspace, \thinspace \rangle$. Let $\mathrm{SO}(V)_{0}$
  be the (split) special orthogonal group over $\mathbb{F}_{q}$. We
  fix a Borel subgroup $B_{0}$ over $\mathbb{F}_{q}$. Let
  $\mathrm{SO}(V) \coloneqq \mathrm{SO}(V)_{0} \otimes
  \overline{\mathbb{F}}_{q}$ and
  $B \coloneqq B_{0} \otimes \overline{\mathbb{F}}_{q}$. Note that a
  Borel subgroup of $\mathrm{SO}(V)$ can be described as the
  stabilizer of a complete isotropic flag:
  \begin{equation}
    \label{eq:133}
    0 \subset V_{1} \subset \cdots \subset V_{d} \subset V_{d}^{\bot}
    \subset \cdots \subset V_{1}^{\bot} \subset V_{\mathbb{F}}.
  \end{equation}
  The (absolute) Weyl group $W$ can be identified with a subgroup of
  $S_{l}$:
  \begin{equation}
    \label{eq:150}
    W = \{ w \in S_{l} : w(i) + w(l+1-i) = l+1 \}.
  \end{equation}
  Let $\mathbb{S}$ be the set of simple reflections $\{ s_{i}, 1 \leq
  i \leq d \}$, where
  \begin{equation}
    \label{eq:151}
    s_{i} = 
    \left\{
      \begin{array}{ll}
        (i, i+1) (l-i, l-i+1), & \text{ if } 1 \leq i \leq d-1, \\
        (d, d+2), & \text{ if } i = d.
      \end{array}
    \right.
  \end{equation}
  The Dynkin diagram of type $B_{d}$ is
  \begin{equation}
    \xymatrix{
       \underset{s_{1}}{\circ} & \underset{s_{2}}{\circ} \ar@{-}[l] &
      {\circ} \ar@{-}[l] &  \ar@{.}[l] &{\circ} \ar@{.}[l] &{\circ} \ar@{-}[l] &
      \underset{s_{d}}{\circ} \ar@{=}[l] |-{\SelectTips{eu}{11}\object@\dir{<}}
    }
  \end{equation}

  Let $w = s_{d} s_{d-1} \cdots s_{1} = (d, d+2) (d, d-1, \ldots, 1)
  (l, l-1, \ldots, d+2)$, then, similarly as \cite[2.2]{MR0393266},
  for a flag $V_{\LargerCdot}$, the $V_{\LargerCdot}$ and
  $\sigma(V_{\LargerCdot})$ are in relative position $w$ if and only
  if the flag $V_{\LargerCdot}$ is of the form:
  \begin{equation}
    \label{eq:153}
    V_{d-i} = V_{d} \cap \sigma(V_{d}) \cap \cdots \cap \sigma^{i}(V_{d}),
  \end{equation}
  for $1 \leq i \leq d$. Let $P$ be the standard parabolic subgroup
  corresponding to $\mathbb{S} - \{ s_{d} \}$, then the natural map
  $\phi \colon X(w) \to \mathrm{SO}(V)/B \to \mathrm{SO}(V)/P$ sending
  $\phi \colon V_{\LargerCdot} \mapsto V_{d}$ is injective. And
  $\mathrm{im}(\phi)$ is the subvariety of $\mathrm{SO}(V)/P$
  parameterizing all the $d$-dimensional isotropic subspaces $V_{d}$
  such that for $1 \leq i \leq d$, we have
  \begin{equation}
    \label{eq:167}
    \dim(V_{d} \cap \sigma(V_{d}) \cap \cdots \cap \sigma^{i}(V_{d}))
    = d-i.
  \end{equation}
\end{eg}

\begin{eg}[The non-split even orthogonal group]
  \label{eg:3}
  Let $V$ be a vector space of even dimension $l = 2 d$ over
  $\mathbb{F}_{q}$, together with ``the'' non-split non-degenerate
  symmetric form $\langle \thinspace, \thinspace \rangle$. We assume
  that $d \geq 2$. Let $\mathrm{SO}(V)_{0}$ be the special orthogonal
  group which is a quasi-split but non-split reductive group over
  $\mathbb{F}_{q}$. We fix a Borel subgroup $B_{0}$ over
  $\mathbb{F}_{q}$.

  Now consider the group
  $\mathrm{SO}(V) \coloneqq \mathrm{SO}(V)_{0} \otimes
  \overline{\mathbb{F}}_{q}$. Let
  $B \coloneqq B_{0} \otimes \overline{\mathbb{F}}_{q}$. Note that
  a Borel subgroup can be described as the stabilizer of some flag
  $W_{\LargerCdot}$ of the form:
  \begin{equation}
    \label{eq:155}
    W_{\LargerCdot} \colon 0 \subset W_{1} \subset \cdots \subset W_{d-2} \subset
    (W_{d} \text{ and } W_{d'}) \subset W_{d-2}^{\bot} \subset \cdots
    \subset W_{1}^{\bot} \subset V,
  \end{equation}
  where $\dim(W_{i}) = i$ for $1 \leq i \leq d-2$, $\dim(W_{d}) =
  \dim(W_{d'}) = d$ and $\dim(W_{d} \cap W_{d'}) = d-1$.

  The absolute Weyl group $W$ can be identified as the subgroup of
  $S_{l}$:
  \begin{equation}
    \label{eq:156}
    W = 
    \left\{
      w \in S_{l} \left|
      \begin{array}{ll}
        w(i) + w(l+1-i) = l+1, \\
        \#\{i , 1 \leq i \leq d : w(i) > d \} \text{ is even}
      \end{array}
    \right. \right\}.
  \end{equation}
  Let $\mathbb{S} = \{ s_{i} : 1 \leq i \leq d \}$ be the set of
  simple reflections, where
  \begin{equation}
    \label{eq:158}
    s_{i} = 
    \left\{
      \begin{array}{ll}
        (i, i+1) (l-i, l-i+1), & 1 \leq i \leq d-1 ,\\
        (d, d+2) (d-1, d+1), & i = d.
      \end{array}
    \right.
  \end{equation}
  The Dynkin diagram of type $D_{d}$ is
  \begin{equation}
    \xymatrix{
      & & & & & & \underset{s_{d-1}}{\circ} \\
      \underset{s_{1}}{\circ} & \underset{s_{2}}{\circ} \ar@{-}[l] &
      {\circ} \ar@{-}[l] &  \ar@{.}[l] &{\circ} \ar@{.}[l] &{\circ}
      \ar@{-}[l] \ar@{-}[ru] \ar@{-}[rd] & \\
      & & & & & & \underset{s_{d}}{\circ} 
    }
  \end{equation}

  Let $w_{1} = s_{d-1}\cdots s_{2} s_{1} =
  (d,\ldots,1)(d+1,\ldots,l)$, then a flag $U_{\LargerCdot}$ lying in
  the Deligne-Lusztig variety $X(w_{1})$ is of the form:
  \begin{align}
    \label{eq:187}
    U_{d-i} & = U_{d'} \cap \delta(U_{d'}) \cap \cdots \cap
              \delta^{i}(U_{d'}), 
  \end{align}
  for $1 \leq i \leq d$. Therefore the natural map $\phi \colon
  X(w_{1}) \to \mathrm{SO}(V)/B \to \mathrm{SO}(V)/P_{d'}$ sending
  $U_{\LargerCdot} \to U_{d}$, where $P_{d'}$ is the standard
  parabolic subgroup corresponding to $\mathbb{S} - \{s_{d-1} \}$, is
  an injection. And $\mathrm{im}(\phi)$ is the subvariety of
  $\mathrm{SO}(V)/P_{d'}$ parameterizing all the maximal isotropic
  subspaces $U$ of $V_{\mathbb{F}}$ such that $U$ lies in the
  $\mathrm{SO}(V)$-orbit of the $d$-dimensional subspace fixing by
  $P_{d'}$, and for $1 \leq i \leq d$, we have
  \begin{equation}
    \label{eq:206}
    \dim( U \cap \delta(U) \cap \cdots \cap \delta^{i}(U)) = d-i.
  \end{equation}

  For $w_{2}= s_{d} s_{d-2} \ldots s_{2} s_{1}$, we have similar
  result.
\end{eg}

Returning to the general case, let $P_{I}$ be the standard parabolic
subgroup of $H$, where $I$ is a subset of the set of simple
reflections $\mathbb{S}$ of $W_{H}$. Let $W_{I}$ be the subgroup of
$W_{H}$ generated by simple reflections in $I$, $W^{I}$
(resp. $\prescript{I}{}{W}$) the set of minimal length representatives
of the cosets in $W_{H}/W_{I}$ (resp.  $W_{I} \backslash W_{H}$). Let
$\prescript{I}{}{W}^{J}$ denote $\prescript{I}{}{W} \cap W^{J}$. Then
we can define the generalized Deligne-Lusztig varieties.

\begin{definition} 
  For each $w \in W_{H}$, the \emph{generalized Deligne-Lusztig
    variety} $X_{P_{I}}(w)$ is defined as
  \begin{equation}
    \label{eq:169}
    X_{P_{I}}(w) \coloneqq \{ g \in H/P_{I} : g^{-1} \sigma(g) \in P_{I}
    w P_{\sigma(I)} \}.
  \end{equation}
\end{definition}

\begin{prop}[{\cite[Lemma 2.1.3]{hoeve-phd}}]
  \label{pr:18}
  For $w \in \prescript{I}{}{W}^{\sigma(I)}$, the Deligne-Lusztig
  variety $X_{P_{I}}(w)$ is smooth of dimension $\ell(w) +
  \ell(W_{\sigma(I)}) - \ell(W_{I \cap \prescript{w}{}{\sigma(I)}})$,
where $\ell(W_{J})$ denotes the maximal length of elements in $W_{J}$
for $J \subset \mathbb{S}$.
\end{prop}

The partial flag variety $H/P_{I}$ can be written as the disjoint
union of all such Deligne-Lusztig varieties indexed by the set
$\prescript{I}{}{W}^{\sigma(I)}$.

\begin{eg}[The odd orthogonal group]
  \label{eg:7}
  Notations are the same as in Example~\ref{eg:2}. We will show that
  the closure of $\mathcal{P}_{w} \coloneqq \phi(X(w))$ is normal
  using the same method as \cite[Proposition 7.3.2]{GH}. Consider the
  inclusion of closures
  \begin{equation}
    \label{eq:216}
    \overline{\mathcal{P}_{w}}  \subset \overline{X_{P}(w)}.
  \end{equation}
  The variety $X_{P}(w)$ is irreducible, so the inclusion
  \eqref{eq:216} is an equality if and only if
  $\ell(w) = \dim(X_{P}(w))$. In this case, let
  $I = \{ s_{1}, s_{2}, \ldots, s_{d-1} \}$ be the type of $P$,
  $w_{\mathrm{min}} = s_{d}$ the minimal representative of $w$ in
  $\prescript{I}{}{W}^{\sigma(I)}$. Note that $\sigma$ acts on the
  Dynkin diagram trivially. Then
  $I \cap \prescript{w_{\mathrm{min}}}{}{\sigma(I)} = \{ s_{1}, s_{2},
  \ldots, s_{d-2} \}$. Therefore
  \begin{equation}
    \label{eq:217}
    \dim(X_{P}(w)) = 1 + \frac{d (d-1)}{2} - 
    \frac{(d-2)(d-1)}{2} = \ell(w).
  \end{equation}
  Then the closure $\overline{\mathcal{P}_{w}}$ is normal, and has
  isolated singularities by \cite[Proposition 7.3.2]{GH}. The closure
  $\overline{\mathcal{P}_{w}}$ can be described as the subvariety of
  $\mathrm{SO}(V) / P$ parameterizing all the $d$-dimensional
  isotropic subspaces $V_{d}$ such that
  \begin{equation}
    \label{eq:278}
    \mathrm{dim}(V_{d} \cap \sigma(V_{d})) \geq d-1.
  \end{equation}
\end{eg}

\begin{eg}
  \label{eg:8}
  Notations are the same as Example~\ref{eg:3}. Let
  $\mathcal{P}_{w_{1}} \coloneqq \phi(X(w_{1}))$. Let
  $I = \{s_{1}, s_{2}, \ldots, s_{d-2}, s_{d} \}$ be the type of
  $P_{w_{1}}$. Note that the Frobenius $\delta$ exchanges $s_{d-1}$
  and $s_{d}$, and fixes all the other $s_{i}$'s. Let
  $w_{\mathrm{min}} = 1$ the minimal representative of $w_{1}$ in
  $\prescript{I}{}{W}^{\delta(I)}$. Then
  $I \cap \prescript{w_{\mathrm{min}}}{}{\delta(W)} = \{ s_{1}, s_{2},
  \ldots, s_{d-2} \}$. We have
  \begin{equation}
    \label{eq:218}
    \dim(X_{P_{d'}}(w_{1})) = \frac{d(d-1)}{2} - \frac{(d-2)(d-1)}{2} =
    \ell(w_{1}). 
  \end{equation}
  Therefore the closure $\overline{\mathcal{P}_{w_{1}}}$ is
  normal. Furthermore, $\overline{\mathcal{P}_{w_{1}}}$ is smooth
  because $\delta(P_{d'}) = P_{d}$. The closure
  $\overline{\mathcal{P}_{w_{1}}} = X_{P_{d'}}(\mathrm{id})$ can be
  described as the subvariety of $\mathrm{SO}(V) / P_{d'}$
  parameterizing all the maximal isotropic subspaces $U$ of $V$ such
  that $U$ lies in the $\mathrm{SO}(V)$-orbit of $W_{d}$ and
  \begin{equation}
    \label{eq:279}
    \mathrm{dim}(U \cap \delta(U)) = d-1.
  \end{equation}
\end{eg}

\begin{remark}
  \label{rmk:15}
  The odd orthogonal case has been listed in \cite[Proposition
  7.3.2]{GH}, which corresponds to the triple
  $(\tilde{C}_{d}, \omega^{\vee}_{1}, \mathbb{S},
  \mathrm{id})$. However, the even orthogonal case is not in the list,
  which corresponds to the triple
  $(\tilde{B}_{m}, \omega_{1}^{\vee}, \mathbb{S},
  \mathrm{id})$. Furthermore, by the same procedure, it is easy to
  check that the triples
  $(\tilde{D}_{l}, \omega^{\vee}_{1}, \mathbb{S}, \mathrm{id})$ and
  $(\tilde{D}_{l}, \omega^{\vee}_{1}, \mathbb{S}, \sigma_{0})$ should
  be also included in the list of \cite[Proposition 7.3.2]{GH}. The
  extra three smooth cases make the list complete.
\end{remark}

\subsection{The group-theoretic approach}
\label{sec:gortz-hes-result}

In this section, we apply the group-theoretic results in \cite{GH} to
our case. There is no harm to look at only one connected component of
the Rapoport-Zink space, because all connected components are
isomorphic to each other by \eqref{eq:41}. By \eqref{eq:120}, we may
work with the group $G_{\mathrm{ad}}(L)'$ instead of
$G_{\mathrm{ad}}(L)$. All notations are the same as in previous
sections.

Let $\tilde{\mathbb{S}} = \{ s_{0}, s_{1}, \ldots, s_{m}\}$ be the set
of affine simple reflections in $\tilde{W}_{\ad}$.  For
$Z \subset \tilde{\mathbb{S}}$, we denote $P_{Z}$ the corresponding
standard parahoric subgroup of $G_{\ad}(L)$.  We will write
$\mathrm{EO}_{\mathrm{cox}}$ instead of
$\mathrm{EO}_{\sigma, \mathrm{cox}}^{\mathbb{S}}(\mu_{\ad})$ in
\cite[5.1]{GH} to lighten the notations because in our case $\sigma$
acts on the affine Dynkin diagram trivially. Let
$\mathbb{J} = J_{\mathrm{ad}}(\mathbb{Q}_{p})'$.

Let us compute the $\mathrm{EO}$ set explicitly first.

\subsubsection{$\mathrm{EO}$ set}
\label{sec:eo-set}

Let $\mathrm{EO}(\mu)$ be the set
$\mathrm{Adm}^{\circ}(\mu) \cap \prescript{J}{}{\tilde{W}_{\ad}}$. Let
$\mathrm{EO}_{\mathrm{cox}}$ be the subset of $\mathrm{EO}(\mu)$
consisting of those $w$ such that $\mathrm{supp}_{\sigma}(w)$ is a
proper subset of $\tilde{\mathbb{S}}$ and $w$ is a $\sigma$-Coxeter
element of $W_{\mathrm{supp}_{\sigma}(w)}$.

Let $\tau$ be the image of $b_{\mathrm{ad}}$ in $\Omega$. For
$v \in \tilde{\mathbb{S}}$ let $d(v)$ be the minimal distance between
the $\tau \sigma$-orbit containing $v$ and the vertex outside
$\mathbb{S}$. Let $\mathscr{J}$ be the set of subsets $\Sigma$ of
$\tilde{\mathbb{S}}$, that is $\tau \sigma$-stable and $d(v) = d(v')$
for any $v, v' \in \Sigma$. For $\Sigma \in \mathscr{J} $ let
$d(\Sigma) \coloneqq d(v)$ for some $v \in \Sigma$, $\Sigma^{\flat}$
the union of all the $\tau \sigma$-orbits $\Sigma'$ that is not
contained in $\Sigma$ and $d(\Sigma') \leq d(\Sigma)$,
$\Sigma^{\sharp}$ the union of all the $\tau \sigma$-orbits $\Sigma'$
such that $d(\Sigma') > d(\Sigma)$. By \cite[Proposition 7.1.1]{GH},
the map
\begin{align}
  \label{eq:62}
  \mathscr{J} & \longrightarrow \mathrm{EO}_{\mathrm{cox}} ,\\
  \Sigma & \longmapsto w_{\Sigma}, \nonumber
\end{align}
is bijective, where $w_{\Sigma}$ is the unique element in
$\mathrm{EO}_{\mathrm{cox}}$ such that
$\mathrm{supp}_{\sigma}(w_{\Sigma}) = \Sigma^{\flat}$. We have $\ell
(w_{\Sigma}) = d(\Sigma)$.

\paragraph{(a) odd case}

In this case, $\tau$ is identity.
\begin{align}
  \label{eq:117}
  \mathrm{EO}_{\mathrm{cox}} & = \{1, s_{0},
  s_{0}s_{1}, \ldots, s_{0}s_{1}\dots s_{m-1} \}. \\
  \label{eq:118}
  \mathscr{J} &= \{ \{ s_{0} \}, \{ s_{1}\} , \ldots, \{ s_{m}\} \}.
\end{align}
If $\Sigma = \{s_{i}\} \in \mathscr{J}$ for some $i$, then
$\Sigma^{\flat} = \{s_{0}, \ldots, s_{i-1}\}$ if $i>0$ or empty
otherwise; $\Sigma^{\sharp} = \{ s_{i+1}, \ldots, s_{m}\}$ if $i<m$ or
empty otherwise; $w_{\Sigma} = s_{0}s_{1}\dots s_{i-1}$ if $i > 0$ or
$1$ otherwise.

For example, if $m = 7$, $\Sigma = \{s_{4} \}$, then
$\Sigma^{\flat} = \{ s_{0}, s_{1}, s_{2}, s_{3}\}$ and
$\Sigma^{\sharp} = \{s_{5}, s_{6}, s_{7} \}$. See the
diagram~\eqref{eq:246}, where $\Sigma^{\flat}$ is surrounded by the
solid frame, and $\Sigma^{\sharp}$ is surrounded by the dashed frame.

\begin{equation}
  \label{eq:246}
  \xymatrix{
    \underset{{0}}{\circ} & \underset{{1}}{\circ} \ar@{=}[l] |-{\SelectTips{eu}{11}\object@\dir{>}} & \underset{{2}}{\circ} \ar@{-}[l] &
   \underset{{3}}{\circ} \ar@{-}[l] & \underset{{4}}{\circ} \ar@{-}[l] &\underset{{5}}{\circ} \ar@{-}[l] &\underset{{6}}{\circ} \ar@{-}[l] &
    \underset{{7}}{\circ} \ar@{=}[l]
    |-{\SelectTips{eu}{11}\object@\dir{>}} \\
    & &  \Sigma^{\flat} & & & & \Sigma^{\sharp}
    \save "1,1"."1,4"*\frm{-}
    \restore
    \save "1,6"."1,8"*\frm{--}
    \restore
  }
\end{equation}

\paragraph{(b) even case}

In this case, $\tau$ switches between $s_{0}$ and $s_{1}$, and fixes
all the other vertices.
\begin{align}
  \label{eq:66}
  \mathrm{EO}_{\mathrm{cox}} &= \{ \tau, s_{0} \tau, s_{0}s_{2}\tau,
  \ldots, s_{0}s_{2}\dots s_{m-1} \tau \}. \\
  \label{eq:119}
  \mathscr{J} &= \{ \{s_{0}, s_{1}\}, \{s_{2}\}, \ldots,
  \{s_{m}\} \}.
\end{align}
For $\Sigma \in \mathscr{J}$, if $\Sigma = \{ s_{0}, s_{1} \}$, then
$\Sigma^{\flat} = \emptyset$,
$\Sigma^{\sharp} = \{ s_{2}, \ldots, s_{m}\}$, $w_{\Sigma} = \tau$; if
$\Sigma = \{s_{i}\}$ for some $i > 1$, then
$\Sigma^{\flat} = \{ s_{0}, s_{1}, \ldots, s_{i-1}\}$,
$\Sigma^{\sharp} = \{ s_{i+1}, \ldots, s_{m}\}$ if $i < m$ or empty
otherwise $w_{\Sigma} = s_{0}s_{2}\dots s_{i-1}\tau$.

For example, if $m = 8$, $\Sigma = \{s_{5} \}$, then
$\Sigma^{\flat} = \{s_{0}, s_{1}, \ldots, s_{4}\}$ and
$\Sigma^{\sharp} = \{s_{6}, s_{7}, s_{8} \}$. See the
diagram~\eqref{eq:247}, where $\Sigma^{\flat}$ is surrounded by the
solid frame, and $\Sigma^{\sharp}$ is surrounded by the dashed frame.

\begin{equation}
  \label{eq:247}
  \xymatrix{
    \underset{{0}}{\circ} \ar@{-}[rd] & & & & & & & \\
    & \underset{{2}}{\circ} & \underset{{3}}{\circ} \ar@{-}[l] &
    \underset{4}{\circ} \ar@{-}[l] & \underset{5}{\circ} \ar@{-}[l]
    &\underset{6}{\circ} \ar@{-}[l] &\underset{7}{\circ} \ar@{-}[l] &
    \underset{{m}}{\circ} \ar@{=}[l]
    |-{\SelectTips{eu}{11}\object@\dir{>}} \\
    \stackrel{{1}}{\circ} \ar@{-}[ur] & & \Sigma^{\flat}
    & & & & \Sigma^{\sharp} &
    \save "1,1"."3,4"*\frm{-}
    \restore
    \save "2,6"."2,8"*\frm{--}
    \restore
  }
\end{equation}

\subsubsection{Lusztig's partition}
\label{sec:lusztigs-partition}

Let $P$ and $Q$ be two parahoric subgroups of $G_{\mathrm{ad}}(L)'$,
similarly to the parabolic case, we define the refinement of $P$ with
respect to $Q$ as
\begin{equation}
  \label{eq:109}
  P^{Q} \coloneqq (P \cap Q) \cdot U_{P},
\end{equation}
where $U_{P}$ is the pro-unipotent radical of $P$. The group $P^{Q}$
is a parahoric subgroup of $G_{\mathrm{ad}}(L)'$ again, and its
pro-unipotent radical is $(P \cap U_{Q}) \cdot U_{P}$, which can be
proved analogously. G\"{o}rtz-He generalize the Lusztig-B\'{e}dard
sequence to the case of affine Weyl groups.
\begin{definition}[Lusztig-B\'{e}dard sequence]
  \label{df:7}
  For $J \subset \tilde{\mathbb{S}}$, let
  $\mathcal{T}(J, \tau \sigma)$ be the set of sequences
  $(J_{i}, w_{i})_{i \geq 0}$ such that
  \begin{itemize}
  \item [(a)] $J_{0} = J$, $w_{0} \in \prescript{J}{}{W_{a}}^{\tau
      \sigma (J)}$,
  \item [(b)] $J_{i+1} = J_{i} \cap \prescript{w_{i}}{}{(\tau
      \sigma(J_{i}))}$, for $i \geq 0$,
  \item [(c)] $w_{i} \in \prescript{J_{i}}{}{W_{a}}^{\tau \sigma
      (J_{i})}$ and $w_{i+1} \in W_{J_{i+1}} w_{i} W_{\tau
      \sigma(J_{i})}$, for $i \geq 1$.
  \end{itemize}
\end{definition}
Let $J_{\infty} \coloneqq J_{i}$, $w_{\infty} \coloneqq w_{i}$, for
$i \gg 0$, then the map $(J_{i}, w_{i}) \mapsto w_{\infty}$ defines a
bijection $\mathcal{T} \to \prescript{J}{}{W_{a}}$. For each parahoric
subgroup $P$ of type $J$, we associate a sequence of parahoric
subgroups
\begin{equation}
  \label{eq:185}
  P^{0} \coloneqq P, \quad P^{i+1} \coloneqq (P^{i})^{b_{\ad} \sigma(P^{i})
    b_{\ad}^{-1}} \text{ for } i \geq 0,
\end{equation}
and a sequence $(J_{i}, w_{i})$, where $J_{i}$ is the type of $P^{i}$
and $w_{i}$ is the relative position of $P^{i}$ and
$b_{\ad} \sigma(P^{i}) b_{\ad}^{-1}$. Then $(J_{i}, w_{i})$ is a
Lusztig-B\'{e}dard sequence. Let $P^{\infty}$ denote $P^{m}$ for
$m \gg 0$. Then $P^{\infty}$ is of type $J_{\infty} \coloneqq J_{m}$
for $m \gg 0$.

For any \emph{affine Deligne-Lusztig variety}
\begin{equation}
  \label{eq:89}
  X_{w}(b) \coloneqq \{ g \in G_{\ad}(L)' / I : g^{-1} b \sigma(g) \in
  I w I \}
\end{equation}
(cf. \cite[Definition 4.1]{MR2141705}), its image
\begin{equation}
  \label{eq:90}
  X_{w}^{f}(b) \coloneqq \{ g \in G_{\mathrm{ad}}(L)'/K_{\mathrm{ad}} 
  : g^{-1} b \sigma(g) \in K_{\mathrm{ad}} \cdot_{\sigma} I w I \}
\end{equation}
under the map $G_{\ad}(L)'/I \rightarrow G_{\ad}(L)'/K_{\mathrm{ad}}$
is called the \emph{fine affine Deligne-Lusztig variety}, where the
superscript $f$ stands for ``fine''.

Now for each $g \in X^{f}_{w}(b_{\mathrm{ad}})$, let
$P \coloneqq \prescript{g}{}{K_{\mathrm{ad}}}$ be the corresponding
parahoric subgroup of $G_{\mathrm{ad}}(L)'$. We write $w = x \tau$ for
some $x \in \prescript{J}{}{W_{a}}$ and $\tau \in \Omega$. Recall that
$w = w_{\Sigma}$. Then the type of $P$ is $J_{0} = \mathbb{S}$. By the
same procedure as in previous paragraph, we get the Lusztig-B\'{e}dard
sequence $(J_{n}, x_{n})$, and by \cite[Lemma 1.4]{MR2299444},
$J_{\infty} = \Sigma^{\sharp}$. Then, by \cite[4.2 (c)(d)]{MR2336607},
the map
$\prescript{g}{}{K_{\mathrm{ad}}} \mapsto
({\prescript{g}{}{K_{\mathrm{ad}}}})^{\infty}$ gives the isomorphism:
\begin{equation}
  \label{eq:100}
  \theta \colon X^{f}_{w}(b_{\mathrm{ad}}) 
  \stackrel{\cong}{\longrightarrow}
  \{ g \in G_{\mathrm{ad}}(L)' / P_{\Sigma^{\sharp}} : g^{-1}
  b_{\mathrm{ad}} \sigma(g) \in P_{\Sigma^{\sharp}} w P_{\Sigma^{\sharp}}\},
\end{equation}
whose inverse map is the natural projection map, i.e. the following
diagram is commutative:
\begin{equation}
  \label{eq:111}
  \xymatrix{
    X_{w}(b_{\mathrm{ad}}) \ar@{->>}[d] \ar@{^{(}->}[r] &
    G_{\mathrm{ad}}(L)'/I \ar[d] \\
    \{ g \in G_{\ad}(L)' /
    P_{\Sigma^{\sharp}} : g^{-1} b_{\mathrm{ad}} \sigma(g) \in
    P_{\Sigma^{\sharp}} w P_{\Sigma^{\sharp}}\} \ar@{^{(}->}[r] \ar[d]^{\cong}
    & G_{\mathrm{ad}}(L)'/P_{\Sigma^{\sharp}} \ar[d] \\
    X^{f}_{w}(b_{\mathrm{ad}}) \ar@{^{(}->}[r] & 
    G_{\mathrm{ad}}(L)'/K_{\mathrm{ad}}.
  }
\end{equation}

\begin{remark}
  \label{rmk:7}
  Using the same trick as Example~\ref{eg:2} and Example~\ref{eg:3},
  we can describe the map $\theta$ in \eqref{eq:100} in terms of
  lattices. Let $g \in G_{\ad}(L)' / P_{\Sigma^{\sharp}}$ such that
  $g^{-1} b_{\ad} \sigma(g) \in P_{\Sigma^{\sharp}} w
  P_{\Sigma^{\sharp}}$. 
  \begin{enumerate}
  \item When $n = 2 m +1$ is odd, let $\Sigma = \{ s_{i} \}$ for some
    $0 \leq i \leq m$, then $P_{\Sigma^{\sharp}}$ is the stabilizer of
    the lattice chain
    \begin{equation}
      \label{eq:179}
      \Lambda_{m-i} \subset \cdots \subset \Lambda_{m-1} \subset
      \Lambda_{m},
    \end{equation}
    Let $\Delta_{j} \coloneqq g \cdot \Lambda_{j}$, then
    $\prescript{g}{}{P_{\Sigma^{\sharp}}}$ is the stabilizer of the
    lattice chain
    \begin{equation}
      \label{eq:180}
      \Delta_{m-i} \subset \cdots \subset \Delta_{m-1} \subset \Delta_{m}.
    \end{equation}
    The condition
    $g^{-1} b_{\ad} \sigma(g) \in P_{\Sigma^{\sharp}} w
    P_{\Sigma^{\sharp}}$ is equivalent to the condition that the pair
    $(\prescript{g}{}{P_{\Sigma^{\sharp}}}, b_{\ad} \sigma
    (\prescript{g}{}{P_{\Sigma^{\sharp}}}))$ lies in the
    $G_{\ad}(L)'$-orbit of
    $(P_{\Sigma^{\sharp}}, w P_{\Sigma^{\sharp}})$. So the lattice
    chain $\Delta_{\LargerCdot}$ is of the form:
    \begin{equation}
      \label{eq:181}
      \Delta_{m-j} = \Delta_{m} \cap (b_{\ad}\sigma)(\Delta_{m}) \cap
      \cdots \cap (b_{\ad} \sigma)^{j} (\Delta_{m}),
    \end{equation}
    for $1 \leq j \leq i$ and
    $\Delta_{m-i} = (b_{\ad} \sigma)(\Delta_{m-i})$. Therefore, the
    map $\theta$ is
    \begin{equation}
      \label{eq:182}
      \Delta_{m} \longmapsto \Delta_{\LargerCdot},
    \end{equation}
    such that $\Delta_{\LargerCdot}$ satisfies condition
    \eqref{eq:181}.
  \item When $n = 2 m$ is even, let $\Sigma = \{ s_{i} \}$ for some $2
    \leq i \leq m$, then $P_{\Sigma^{\sharp}}$ is the stabilizer of
    the lattices
    \begin{equation}
      \label{eq:184}
      \Lambda_{m-i} \subset \cdots \subset \Lambda_{m-2} \subset
      (\Lambda_{m} \text{ and } \Lambda_{m'}) \subset \Lambda_{m-2}^{\vee}
      \subset \cdots \subset \Lambda_{m-i}^{\vee}.
    \end{equation}
    Note that $\Lambda_{m} \cap \Lambda_{m'} = \Lambda_{m-1}$, and the
    lattice $\Lambda_{m'}$ is uniquely determined by $\Lambda_{m-1}$
    and $\Lambda_{m}$, because there are exactly two isotropic lines
    in the hyperbolic plane $\Lambda_{m-1}^{\vee} / \Lambda_{m-1}$. We
    also use the notation $\Lambda_{m}^{\dagger}$ to denote the unique
    isotropic line (hence the lattice) determined by $\Lambda_{m}$ and
    $\Lambda_{m-1}$. Let $\Delta_{j} \coloneqq g \cdot \Lambda_{j}$
    for $m-i \leq j \leq m$ and $j = m'$. Then the condition
    $g^{-1} b_{\ad} \sigma(g) \in P_{\Sigma^{\sharp}} w
    P_{\Sigma^{\sharp}}$ implies that $\Delta_{\LargerCdot}$ is of the
    form:
    \begin{align}
      \label{eq:189}
      \Delta_{m'} & = b_{\ad} \sigma (\Delta_{m}) \\
      \Delta_{m-j} & = \Delta_{m} \cap (b_{\ad} \sigma)(\Delta_{m})
                     \cap \cdots \cap (b_{\ad} \sigma)^{j}(\Delta_{m})
                     \nonumber
    \end{align}
    for $1 \leq j \leq i$ and $\Delta_{m-i} = (b_{\ad}
    \sigma)(\Delta_{m-i})$. Therefore, the map $\theta$ is 
    \begin{equation}
      \label{eq:190}
      \Delta_{m} \longmapsto \Delta_{\LargerCdot},
    \end{equation}
    such that $\Delta_{\LargerCdot}$ satisfies the condition described
    above.
    
    For $\Sigma = \{s_{0}, s_{1} \}$, $P_{\Sigma^{\sharp}}$ is the
    stabilizer of the lattice chain:
    \begin{equation}
      \label{eq:288}
      \Lambda_{m-1} \subset \Lambda_{m}.
    \end{equation}
    Note that $P_{\Sigma^{\sharp}}$ is also the stabilizer of the
    lattices
    \begin{equation}
      \label{eq:289}
      \Lambda_{m-1} \subset (\Lambda_{m} \text{ and } \Lambda_{m'})
      \subset \Lambda_{m-1}^{\vee}.
    \end{equation}
    Therefore, the map $\theta$ is
    \begin{equation}
      \label{eq:290}
      \Delta \longmapsto \Delta_{\LargerCdot},
    \end{equation}
    where $\Delta_{\LargerCdot} = \{ \Lambda_{m-1} \subset (\Lambda_{m}
    \text{ and } \Lambda_{m'}) \subset  \Lambda^{\vee}_{m-1}
    \}$ satisfying $\Lambda_{m'} = b_{\ad} \sigma(\Delta_{m})$.
  \end{enumerate}
\end{remark}

\subsubsection{Bruhat-Tits stratification}
\label{sec:bruh-tits-strat-1}

\begin{prop}
  \label{pr:ulrich-he}
  We have the following decomposition:
  \begin{align}
    \label{eq:58}
    X(\mu_{\ad}, b_{\ad})_{K_{\ad}}' & = \biguplus_{w \in
                                           \mathrm{EO}_{\mathrm{cox}}}
                                           X_{w}^{f} (b_{\ad}), \\
    \label{eq:295}
    & = \biguplus_{\Sigma \in \mathscr{J}}
    \biguplus_{j \in \mathbb{J} / \mathbb{J} \cap
      P_{\tilde{\mathbb{S}} \backslash \Sigma}}
    j\cdot Y_{\Sigma^{\sharp}}(w_{\Sigma}),
  \end{align}
  where
  \begin{equation}
    \label{eq:296}
    Y_{\Sigma^{\sharp}}(w) = \{ g \in P_{\tilde{\mathbb{S}}-\Sigma} /
    P_{\Sigma^{\sharp}} : g^{-1} b_{\ad} \sigma(g) \in
    P_{\Sigma^{\sharp}} w P_{\Sigma^{\sharp}} \}.
  \end{equation}
  Moreover, the natural projection
  $G_{\ad}(L)'/I \to G_{\ad}(L)'/K_{\mathrm{ad}}$ induces an
  isomorphism from the Deligne-Lusztig variety
  \begin{equation}
    \label{eq:297}
    Y(w) = \{ g \in P_{\Sigma^{\flat}} / I : g^{-1} b_{\ad} \sigma(g)
    \in I w I \}
  \end{equation}
  in $P_{\Sigma^{\flat}}/I$ to $Y_{\Sigma^{\sharp}}(w)$,
  i.e. $Y_{\Sigma^{\sharp}}(w)$ is a fine Deligne-Lusztig variety.
\end{prop}

\begin{proof}
  To prove \eqref{eq:58} is just to show that the triple $(G_{\ad},
  \mu_{\ad}, \mathbb{S})$ is of Coxeter type, which has been
  completely listed in \cite[Theorem 5.1.2]{GH}. For the second
  equality, similarly to \cite[Proposition 2.2.1]{GH}, the set
  \begin{equation}
    \label{eq:294}
    \{ g \in G_{\mathrm{ad}}(L)' / P_{\Sigma^{\sharp}} : g^{-1}
    b_{\mathrm{ad}} \sigma(g) \in P_{\Sigma^{\sharp}} w P_{\Sigma^{\sharp}}\}
  \end{equation}
  can be decomposed into a disjoint union of some classical
  Deligne-Lusztig varieties $Y_{\Sigma^{\sharp}}(w)$, then by
  \eqref{eq:100}, we get the desired decomposition. 
\end{proof}

\begin{remark}
  \label{rmk:54}
  Note that the reductive quotient $\bar{P}_{\Sigma^{\flat}}$ of
  $P_{\Sigma^{\flat}}$ has the Dynkin diagram $\Sigma^{\flat}$,
  i.e. when $n = 2 m + 1$ is odd, $\bar{P}_{\Sigma^{\flat}}$ is an
  algebraic group of type $B_{m}$; when $n = 2 m$,
  $\bar{P}_{\Sigma^{\flat}}$ is of type $D_{m}$.
\end{remark}

\begin{remark}
  \label{rmk:6}
  When $n = 2 m + 1$ is odd, $w_{\Sigma}$ is a Coxeter element of
  $W_{\Sigma^{\flat}}$. However, when $n = 2 m$, $w_{\Sigma}$ is not a
  Coxeter element, but a $\sigma$-Coxeter element. Let $w' \coloneqq w
  \tau^{-1}$, $\delta \coloneqq \tau \sigma \tau^{-1}$, then $\delta$
  is the twisted Frobenius on the reductive quotient
  $\bar{P}_{\Sigma^{\flat}}$, and $w'$ is a $\delta$-twisted Coxeter
  element. We have
  \begin{equation}
    \label{eq:178}
    Y(w) = \{ g \in \bar{P}_{\Sigma^{\flat}} / \bar{I} : g^{-1}
    \delta(g) \in \bar{I} w' \bar{I} \},
  \end{equation}
  where $\bar{I}$ is the image of $I$ in
  $\bar{P}_{\Sigma^{\flat}}$. Note that $\bar{I}$ is a $\delta$-stable
  Borel subgroup of $\bar{P}_{\Sigma^{\flat}}$.
\end{remark}

Now let us consider the closure of $Y(w)$ in the partial flag variety
$P_{\Sigma^{\flat}} / P_{\Sigma^{\flat} \cap \mathbb{S}}$.

\begin{prop}
  \label{pr:7}
  For each $w = w_{\Sigma} \in \mathrm{EO}_{\mathrm{cox}}$, we have
  \begin{equation}
    \label{eq:99}
    \overline{Y(w)} = \coprod_{(\Sigma')^{\flat} \subset
      \Sigma^{\flat}}
    \coprod_{\substack{(\mathbb{J} \cap P_{\tilde{\mathbb{S}} - \Sigma}) \cap
        \prescript{j}{}{(\mathbb{J} \cap P_{\tilde{\mathbb{S}}- \Sigma'})} \neq \emptyset ,\\
        j \in \mathbb{J} / (\mathbb{J} \cap P_{\tilde{\mathbb{S}} - \Sigma'})}}
    j Y(w_{\Sigma'}).
  \end{equation}
\end{prop}

\begin{proof}
  Let $Q = P_{\Sigma^{\flat} \cap \mathbb{S}}$,
  $\mathrm{pr} \colon P_{\Sigma^{\flat}} / I \to P_{\Sigma^{\flat}} /
  Q$ the natural projection which is proper. Then
  \begin{equation}
    \label{eq:186}
    \overline{Y(w)} = \bigcup_{v \leq w} \mathrm{pr}(Y(v)).
  \end{equation}
  The rest of the proof is exactly the same as \cite[Theorem
  7.2.1]{GH}, so we omit it.
\end{proof}

\begin{remark}
  \label{rmk:10}
  For $i \in \mathbb{J}$, $i \cdot Y_{\Sigma^{\sharp}}(w)$ is a
  Deligne-Lusztig variety in the partial flag variety
  $\prescript{i}{}{P_{\tilde{\mathbb{S}} - \Sigma}} /
  \prescript{i}{}{P_{\Sigma^{\sharp}}}$, more precisely,
  \begin{equation}
    \label{eq:203}
    i \cdot Y_{\Sigma^{\sharp}}(w) = \{ x \in 
    \prescript{i}{}{P_{\tilde{\mathbb{S}} - \Sigma}} /
    \prescript{i}{}{P_{\Sigma^{\sharp}}} : x^{-1} b_{\ad} \sigma (x)
    \in \prescript{i}{}{P_{\Sigma^{\sharp}}} i w \sigma(i)^{-1}
    (\prescript{\sigma(i)}{}{P_{\Sigma^{\sharp}}}),
  \end{equation}
  which is isomorphic to $Y_{\Sigma^{\sharp}}(w)$. Using the same
  method as the proof of Proposition~\ref{pr:7}, it is easy to show
  that 
  \begin{equation}
    \label{eq:204}
    \overline{i \cdot Y_{\Sigma^{\sharp}}(w)} 
    = \coprod_{(\Sigma')^{\flat} \subset
      \Sigma^{\flat}}
    \coprod_{\substack{i \cap j \neq \emptyset ,\\
        j \in \mathbb{J} / (\mathbb{J} \cap P_{\tilde{\mathbb{S}} - \Sigma'})}}
    j \cdot Y(w_{\Sigma'}).
  \end{equation}
\end{remark}

The closure relations can be described by the rational Bruhat-Tits
building of $\mathbb{J}$.

\begin{prop}[{\cite[Proposition 7.2.2]{GH}}]
  \label{pr:13}
  Let $i, j \in \mathbb{J}$, $\Sigma, \Sigma' \in \mathscr{J}$, the
  following are equivalent\footnote{The proposition in loc cit is not
    correct, we should assume $i$ and $j$ have the same Kottwitz
    invariants. However in our case this is true because we are
    working with the groups $G_{\ad}(L)'$ and
    $\mathbb{J} = J_{\ad}(\mathbb{Q}_{p})'.$}:
  \begin{enumerate}
  \item $i (\mathbb{J} \cap P_{\tilde{\mathbb{S}} - \Sigma}) \cap j
    (\mathbb{J} \cap P_{\tilde{\mathbb{S}} - \Sigma'}) \neq
    \emptyset$,
  \item $\prescript{i}{}{(\mathbb{J} \cap P_{\tilde{\mathbb{S}} -
        \Sigma})} \cap \prescript{j}{}{(\mathbb{J} \cap
      P_{\tilde{\mathbb{S}} - \Sigma'})}$ contains an Iwahori subgroup
    of $\mathbb{J}$,
  \item The faces in the building of $\mathbb{J}$ corresponding to
    $\prescript{i}{}{(\mathbb{J} \cap P_{\tilde{\mathbb{S}} -
        \Sigma})}$ and $\prescript{j}{}{(\mathbb{J} \cap
      P_{\tilde{\mathbb{S}} - \Sigma'})}$ are neighbors.
  \end{enumerate}
\end{prop}

\subsection{Crucial lemma}
\label{sec:essential-lemma}

Recall $\chi \coloneqq \eta \pi^{-1} \mathcal{F}$ in
section~\ref{sec:local-pel-datum}. For each
$M \in \mathcal{S}(\mathbb{F})$ and $r \in \mathbb{Z}_{\geq 1}$, we
define the lattices
\begin{equation}
  \label{eq:97}
  \Xi_{r}(M) \coloneqq M + \chi(M) + \cdots + \chi^{r}(M).
\end{equation}
By \cite[Proposition 2.17]{MR1393439}, $\Xi_{n-1}(M)$ is invariant
under $\chi$. Note that when $n$ is even, $M \stackrel{1}{\subset} M +
\chi(M)$ by Proposition~\ref{pr:2}.

\begin{lemma}
  \label{lm:3}
  Let $d$ be the minimal number such that $\Xi_{d}(M)$ is
  $\chi$-stable, then $0 \leq d \leq n/2$ and we have the following
  long lattice chain
  \begin{equation}
    \label{eq:98}
    M \stackrel{1}{\subset} \Xi_{1}(M) \stackrel{1}{\subset} \cdots
    \stackrel{1}{\subset} \Xi_{d}(M) \subset \Xi_{d}(M)^{\vee}
    \stackrel{1}{\subset} \cdots \stackrel{1}{\subset}
    \Xi_{1}(M)^{\vee} \subset M^{\vee}.
  \end{equation}
  Furthermore, if $n$ is even, $1 \leq d \leq n/2$.
\end{lemma}

\begin{proof}
  Recall that we have a bijection
  \begin{equation}
    \label{eq:112}
    \Phi_{\mathrm{ad}} \colon X(\mu_{\mathrm{ad}}, b_{\mathrm{ad}})_{K_{\mathrm{ad}}}'
    \longrightarrow \mathcal{S}(\mathbb{F}).
  \end{equation}
  Let $g K_{\ad}$ be the pre-image of $M$ for some
  $g \in G_{\ad}(L)'$, then $g K_{\ad}$ corresponds to the parahoric
  subgroup
  \begin{equation}
    \label{eq:113}
    \mathrm{Stab}_{G_{\mathrm{ad}}}(M \subset M^{\vee} \subset
    \pi^{-1}M), 
  \end{equation}
  which is equal to $\prescript{g}{}{K_{\ad}}$. Then by \eqref{eq:58},
  there exists a unique $w \in \mathrm{EO}_{\mathrm{cox}}$ such that
  $g K_{\mathrm{ad}} \in X^{f}_{w}(b_{\mathrm{ad}})$. And $w$ is of
  the form $w = w_{\Sigma}$ for some $\Sigma \in \mathscr{J}$. So
  $(\prescript{g}{}{K_{\mathrm{ad}}})^{\infty}$ is of type
  $\Sigma^{\sharp}$ and the natural projection
  $G_{\mathrm{ad}}(L)'/P_{\Sigma^{\sharp}} \to
  G_{\mathrm{ad}}(L)'/K_{\mathrm{ad}}$ sending
  $(\prescript{g}{}{K_{\mathrm{ad}}})^{\infty}$ to
  $\prescript{g}{}{K_{\mathrm{ad}}}$ by \eqref{eq:100}. In other
  words, the lattice $M$ sits inside a long lattice chain whose
  connected stabilizer is the parahoric subgroup
  $(\prescript{g}{}{K_{\mathrm{ad}}})^{\infty}$. The lattice chain
  corresponding to $(\prescript{g}{}{K_{\mathrm{ad}}})^{\infty}$ is
  \begin{equation}
    \label{eq:114}
    M \stackrel{1}{\subset} \Xi_{1}(M) \stackrel{1}{\subset} \cdots
    \stackrel{1}{\subset} \Xi_{d}(M) \subset \Xi_{d}(M)^{\vee}
    \stackrel{1}{\subset} \cdots \stackrel{1}{\subset}
    \Xi_{1}(M)^{\vee} \subset M^{\vee},
  \end{equation}
  where $\Xi_{d}(M)$ is $\chi$-stable and
  $\Xi_{d}(M) = \Xi_{d+1}(M) = \cdots$ by
  Remark~\ref{rmk:7}. Furthermore, we have
  \begin{equation}
    \label{eq:191}
    d = 
    \left\{
      \begin{array}{ll}
        \ell(w_{\Sigma}) = |\Sigma^{\flat}| ,
        & \text{ if } n \text{ is odd,} \\
        \ell(w_{\Sigma})-1 = |\Sigma^{\flat}|, 
        & \text{ if } n \text{ is even and } \Sigma = \{s_{i} \}
          \text{ for } 2 \leq i \leq m, \\
        \ell(w_{\Sigma}) = 1, & \text{ if } n \text{ is even and } \Sigma = \{ s_{0},
                                s_{1} \},
      \end{array}   
    \right.
  \end{equation}
  by the calculations of $\mathscr{J}$ in \eqref{eq:118} and
  \eqref{eq:119}.
\end{proof}

\begin{remark}
  \label{rmk:4}
  The terminology ``crucial lemma'' is inherited from \cite[Lemma
  2.1]{MR2666394}, and the lemma does play a ``crucial'' role in the
  theory of Bruhat-Tits stratification. Since the work of Vollaard
  \cite{MR2666394} and Vollaard-Wedhorn \cite{MR2800696},
  Rapoport-Terstiege-Wilson \cite{MR3175176} and Howard-Pappas
  \cite{MR3272278} adopt almost the same approach to the Bruhat-Tits
  stratification, i.e. proving some variant of ``crucial lemma'', see
  \cite[Proposition 4.1]{MR3175176} and \cite[Proposition
  2.19]{MR3272278}. However, the proof of crucial lemmas in all the
  mentioned literature is elementary and not conceptual so that one
  can only prove them case by case. Thanks to Lusztig's work in
  \cite{MR2336607}, we give the ``crucial lemma'' a conceptual proof
  using a group-theoretic method.
\end{remark}

Let $\Xi_{\infty}(M) \coloneqq \Xi_{m}(M)$ for $m \gg 0$. Via the
identification \eqref{eq:38}, the $\chi$-invariant lattice
$\Xi_{\infty}(M)$ can be viewed as an $\mathcal{O}_{F}$-lattice in the
vector space $C$. And we have
$\pi \cdot \Xi_{\infty}(M)^{\vee} \subset \Xi_{\infty}(M) \subset
\Xi_{\infty}(M)^{\vee}$.

\begin{definition}
  \label{df:10}
  An $\mathcal{O}_{F}$-lattice $\Lambda$ in $C$ is called a vertex
  lattice if
  $\Lambda \subset \Lambda^{\sharp} \subset \pi^{-1} \Lambda$, where
  $\Lambda^{\sharp}$ is the dual of $\Lambda$ with respect to the
  hermitian form $\psi$ in Section \ref{sec:local-pel-datum}. The
  dimension of the $\mathbb{F}_{p}$-vector space
  $\Lambda / \pi \Lambda^{\sharp}$ is called the type of the lattice,
  denoted by $t(\Lambda)$.
\end{definition}

For $M \in \mathcal{S}(\mathbb{F})$, it is easy to see that the
lattice $\Xi_{\infty}(M)$ is a vertex lattice and its type $t$ is
\begin{equation}
  \label{eq:192}
  t = 
  \left\{
    \begin{array}{ll}
      2 d + 1, & \text{ if } n \text{ is odd,} \\
      2 d, & \text{ if } n \text{ is even.}
    \end{array}
  \right.
\end{equation}

\begin{remark}
  \label{rmk:8}
  Our definition of vertex lattices is slightly different from the one
  in \cite[Definition 3.1]{MR3175176}, an $\mathcal{O}_{F}$-lattice
  $\Delta$ is a vertex lattice in loc. cit. if and only if
  $\Delta^{\sharp}$ is a vertex lattice in our sense.
\end{remark}

\begin{prop}[Properties of vertex lattices]
  \label{pr:10}
  Let $\Lambda, \Lambda'$ be two vertex lattices.
  \begin{enumerate}
  \item The type of $\Lambda$ has the same parity as $n$.
  \item The inclusion $\Lambda \subset \Lambda'$ implies
    $t(\Lambda) \leq t(\Lambda')$, and in this case, the equality
    holds if and only if $\Lambda = \Lambda'$.
  \item If $t(\Lambda) = t(\Lambda')$, then either
    $\Lambda = \Lambda'$ or $\Lambda \not \subset \Lambda'$ and
    $\Lambda' \not \subset \Lambda$.
  \item The intersection $\Lambda \cap \Lambda'$ is a vertex lattice
    if and only if $\Lambda^{\sharp} \subset \pi^{-1} \Lambda'$.
  \item When $n$ is odd, for each odd number $t$ satisfying
    $1 \leq t \leq n$, there exists a vertex lattice of type $t$.
  \item When $n$ is even, for each even number $t$ satisfying
    $2 \leq t \leq n$, there exists a vertex lattice of type $t$, but
    there is no vertex lattice of type $0$.
  \end{enumerate}
\end{prop}

\begin{proof}
  We will prove $(6)$ and leave the rest to the reader. When $n = 2 m$
  is even, the hermitian space $(C, \psi)$ is non-split. So there is
  no vertex lattice of type $0$, because a lattice is of type $0$ if
  and only if it is a $\pi$-modular lattice, which exists if and only
  if $C$ is split by Lemma~\ref{split-lemma}. We may assume $C$ is the
  direct product of an $(n-2)$ dimensional split hermitian space with
  the unique non-split $2$ dimensional hermitian space. Note that
  every lattice in the non-split $2$ dimensional hermitian space is
  self-dual. Then similarly to the construction in the odd case, for
  each even number $t$ satisfying $2 \leq t \leq n$, there exists a
  lattice $\Delta$ of type $t$.
\end{proof}

Let $\mathcal{B}$ be the set of vertex lattices in $C$. Two vertex
lattices $\Lambda$ and $\Lambda'$ are called neighbors if $\Lambda
\subset \Lambda'$ or $\Lambda \subset \Lambda'$. A $d$-simplex is a
vertex lattice chain:
\begin{equation}
  \label{eq:137}
  \Lambda_{0} \subset \Lambda_{1} \subset \cdots \subset \Lambda_{d}
  \subset \pi^{-1} \Lambda_{0}.
\end{equation}
Then $\mathcal{B}$ forms a simplicial complex which is connected and
isomorphic to the (rational) Bruhat-Tits building of $\mathbb{J}$ by
\cite[Proposition 3.4]{MR3175176}.

\subsection{The set structure of Bruhat-Tits stratification}
\label{sec:set-structure-bruhat}

\begin{definition}
  \label{df:1}
  For each vertex lattice $\Lambda$,
  \begin{equation}
    \label{eq:123}
    \mathcal{S}_{\Lambda}(\mathbb{F}) \coloneqq
    \{ M \in \mathcal{S}(\mathbb{F}) : M \subset \Lambda \}.
  \end{equation}
\end{definition}

\begin{prop}
  \label{pr:8}
  \leavevmode
  \begin{enumerate}
  \item $\mathcal{S}(\mathbb{F}) = \bigcup_{\Lambda \in \mathcal{B}}
    \mathcal{S}_{\Lambda}(\mathbb{F})$.
  \item Let $\Lambda, \Lambda'$ be two vertex lattices, then the
    inclusion $\Lambda \subset \Lambda'$ implies that
    $\mathcal{S}_{\Lambda}(\mathbb{F}) \subset
    \mathcal{S}_{\Lambda'}(\mathbb{F})$.
  \item Let $\Lambda, \Lambda'$ be two vertex lattices, then
    \begin{equation}
      \label{eq:124}
      \mathcal{S}_{\Lambda}(\mathbb{F}) \cap
      \mathcal{S}_{\Lambda'}(\mathbb{F}) =
      \left\{
        \begin{array}{ll}
          \mathcal{S}_{\Lambda \cap \Lambda'}(\mathbb{F}), &
          \text{ if } \Lambda \cap \Lambda' 
          \text{ is a vertex lattice,} \\
          \emptyset, & \text{otherwise.}
        \end{array}
      \right.
    \end{equation}
  \end{enumerate}
\end{prop}

\begin{proof}
  We will prove $(3)$ and leave the rest to the reader. If $\Lambda
  \cap \Lambda'$ is a vertex lattice, then $M \in
  \mathcal{S}_{\Lambda}(\mathbb{F}) \cap
  \mathcal{S}_{\Lambda'}(\mathbb{F})$ implies that $M \subset \Lambda
  \cap \Lambda'$, in other words, $M \in \mathcal{S}_{\Lambda \cap
    \Lambda'}(\mathbb{F})$. If $\Lambda \cap \Lambda'$ is not a vertex
  lattice, and $\mathcal{S}_{\Lambda}(\mathbb{F}) \cap
  \mathcal{S}_{\Lambda'}(\mathbb{F})$ is non-empty, we take $M \in
  \mathcal{S}_{\Lambda}(\mathbb{F}) \cap
  \mathcal{S}_{\Lambda'}(\mathbb{F})$, then we have
  \begin{align}
    \label{eq:195}
    \pi (\Lambda)^{\sharp} \subset & M \subset 
                                     \Lambda \subset \Lambda^{\sharp}, \\
    \pi (\Lambda')^{\sharp} \subset & M \subset 
                                      \Lambda' \subset (\Lambda')^{\sharp}.
  \end{align}
  In particular, we have
  \begin{align}
    \label{eq:196}
    \pi (\Lambda^{\sharp} + (\Lambda')^{\sharp}) \subset M \subset
    \Lambda \cap \Lambda' \subset \Lambda^{\sharp} + (\Lambda')^{\sharp},
  \end{align}
  which implies that $\Lambda \cap \Lambda'$ is a vertex lattice,
  contrary to the assumption.
\end{proof}

\begin{definition}
  \label{df:3}
  For each vertex lattice $\Lambda$,
  \begin{equation}
    \label{eq:126}
    \mathcal{S}_{\Lambda}^{\circ}(\mathbb{F}) \coloneqq \{ M \in
    \mathcal{S}(\mathbb{F}) : \Xi_{\infty}(M) = \Lambda \}.
  \end{equation}
\end{definition}

\begin{prop}
  \label{pr:11}
  \leavevmode
  \begin{enumerate}
  \item $\mathcal{S}_{\Lambda}^{\circ} (\mathbb{F}) =
    \mathcal{S}_{\Lambda}(\mathbb{F}) \backslash \bigcup_{\Lambda'
      \subsetneq \Lambda} \mathcal{S}_{\Lambda'}(\mathbb{F})$.
  \item
    $\mathcal{S}(\mathbb{F}) = \biguplus_{\Lambda \in \mathcal{B}}
    \mathcal{S}_{\Lambda}^{\circ}(\mathbb{F})$ and
    $\mathcal{S}_{\Lambda}(\mathbb{F}) = \biguplus_{\Lambda' \subset
      \Lambda} \mathcal{S}_{\Lambda'}^{\circ}(\mathbb{F})$.
  \end{enumerate}
\end{prop}

We leave the proof to the reader.

For a vertex lattice $\Lambda$, let $\mathbb{B}_{\Lambda}$ be the
$\mathbb{F}_{p}$-vector space $\Lambda / \pi \Lambda^{\sharp}$ of
dimension $t(\Lambda)$. The form $\psi$ in \eqref{eq:40} induces a
$\mathbb{F}_{p}$-valued bilinear symmetric form $\bar{\psi}$ on
$\mathbb{B}_{\Lambda}$ (because $\Lambda \subset \Lambda^{\sharp}$)
defined by
\begin{equation}
  \label{eq:197}
  \bar{\psi}(x, y) \coloneqq \overline{ \psi(x, y)},
\end{equation}
where the overline denotes the reduction modulo $\pi$. Let
$\mathbb{B}_{\Lambda}$ denote the orthogonal space
$(\Lambda / \pi \Lambda^{\sharp}, \bar{\psi})$ by abuse of notation.

Via the identification \eqref{eq:38}, we identify $\psi \otimes L$
with the twisted form $\delta \varphi$ in
section~\ref{sec:local-pel-datum}, viewing $\Lambda$ as a lattice in
$N$, then $\Lambda^{\sharp} = \Lambda^{\vee}$.

\begin{lemma}
  \label{lm:11}
  The symmetric form $\bar{\psi}$ is non-degenerate.
\end{lemma}

The proof is trivial.

Let $\mathrm{SO}(\mathbb{B}_{\Lambda})$ be the special orthogonal
group with respect to the orthogonal space $\mathbb{B}_{\Lambda,
  \mathbb{F}} \coloneqq \mathbb{B}_{\Lambda} \otimes \mathbb{F}$
defined over $\mathbb{F}_{p}$. Recall that via the identification
\eqref{eq:38}, $\chi = \mathrm{id} \otimes \mathrm{Frob}_{\mathbb{F} /
  \mathbb{F}_{p}}$. Let $B$ be a fixed $\chi$-stable Borel subgroup of
$\mathrm{SO}(\mathbb{B}_{\Lambda})$.

For each $M \in \mathcal{S}(\mathbb{F})$, let $\Lambda =
\Xi_{\infty}(M)$ and $\bar{M} \coloneqq M / \pi \Lambda^{\sharp}$,
then $\bar{M}^{\bot} = \pi M^{\vee} / \pi \Lambda^{\sharp}$ and thus
$\bar{M}^{\bot}$ is a maximal isotropic subspace in
$\mathbb{B}_{\Lambda, \mathbb{F}}$ of dimension
$[\frac{t(\Lambda)}{2}]$. Note that by Remark~\ref{rmk:9}, every $M$
lies in the same $G_{\ad}(L)'$-orbit, and hence every $\bar{M}$ lies
in the same $\mathrm{SO}(\mathbb{B}_{\Lambda})$-orbit. Let $Q'$ be the
standard maximal parabolic subgroup corresponding to the
$\mathrm{SO}(\mathbb{B}_{\Lambda})$-orbit of some (or equivalently
any) $M \in \mathcal{S}(\mathbb{F})$.

\begin{lemma}
  \label{lm:12}
  The map
  \begin{align}
    \label{eq:199}
    \mathcal{S}_{\Lambda}^{\circ}(\mathbb{F}) & \longrightarrow
    \mathrm{SO}(\mathbb{B}_{\Lambda}) / Q', \\
    M & \longmapsto \mathrm{Stab}{(\bar{M})} ,
        \nonumber
  \end{align}
  is injective.
\end{lemma}

We leave the proof to the reader.

\begin{remark}
  \label{rmk:13}
  The proof of Lemma~\ref{lm:12} also shows that the map 
  \begin{align}
    \label{eq:19999}
    \mathcal{S}_{\Lambda}(\mathbb{F}) & \longrightarrow
    \mathrm{SO}(\mathbb{B}_{\Lambda}) / Q', \\
    M & \longmapsto \mathrm{Stab}{(\bar{M})} ,
        \nonumber
  \end{align}
  is injective.
\end{remark}

\begin{prop}
  \label{pr:12}
  The map $\Phi_{\ad}$ induces a bijection
  \begin{equation}
    \label{eq:127}
    j \cdot Y_{\Sigma^{\sharp}}(w) \longrightarrow 
    \mathcal{S}_{\Lambda}^{\circ}(\mathbb{F}),
  \end{equation}
  for each $j \in \mathbb{J}$ and
  $w = w_{\Sigma} \in \mathrm{EO}_{\mathrm{cox}}$, where
  $\Lambda = j \cdot \Xi_{\infty}(M)$ for some (or equivalently any)
  $M \in Y_{\Sigma^{\sharp}}(w)$.
\end{prop}

\begin{proof}
  Assume $j = 1$ firstly.  Let $g \in Y_{\Sigma^{\sharp}}(w)$,
  $\dot{g}$ a lifting of $g$ in $G(L)'$, $M = \dot{g} \mathbb{M}$. Let
  $\Lambda \coloneqq \Xi_{\infty}(M)$.
  
  By Proposition~\ref{pr:ulrich-he}, we have the following diagram
  \begin{equation}
    \label{eq:194}
    \xymatrix{
      Y_{\Sigma^{\sharp}}(w) \ar@{^{(}->}[r] 
      \ar@{^{(}->}[rd]_{\phi}
      & P_{\tilde{\mathbb{S}} - \Sigma} /P_{\Sigma^{\sharp}}
      \ar@{^{(}->}[r] \ar@{->>}[d]^{\mathrm{pr}}
      & G_{\ad}(L)' / P_{\Sigma^{\sharp}} \ar[d] \\
      & P_{\Sigma^{\flat}} / P_{\Sigma^{\flat} \cap \mathbb{S}}
      \ar@{^{(}->}[r]
      & G_{\ad}(L)' / K_{\ad}.
    }
  \end{equation}
  The condition $g \in \phi(Y_{\Sigma^{\sharp}}(w)) $ implies that
  $(\prescript{g}{}{K_{\ad}})^{\infty} \in P_{\tilde{\mathbb{S}} -
    \Sigma} / P_{\Sigma^{\sharp}}$, which implies that
  $\Xi_{\infty}(M) = \Lambda_{m-i}^{\vee}$ if $\Sigma = \{ s_{i} \}$
  by Remark~\ref{rmk:7} and Lemma~\ref{lm:3}, where $\Lambda_{m-i}$ is
  the $(m-i)$-th standard lattice in the
  subsection~\ref{sec:lattice-model-bruhat}. In other words, for any
  $g_{1}, g_{2} \in Y_{\Sigma^{\sharp}(w)}$, we get the same vertex
  lattice
  $\Xi_{\infty}(g_{1} \mathbb{M}) = \Xi_{\infty}(g_{2}
  \mathbb{M})$. Therefore the map $\Phi_{\ad}$ takes
  $Y_{\Sigma^{\sharp}}(w)$ into
  $\mathcal{S}_{\Lambda}^{\circ}(\mathbb{F})$.

  We write $Q \coloneqq P_{\Sigma^{\flat} \cap \mathbb{S}}$. Then
  $\Sigma^{\flat} \cap \mathbb{S} = \Sigma^{\flat} - \{ s_{0} \}$ and
  the image $\bar{Q}$ of $Q$ in the reductive quotient
  $\bar{P}_{\Sigma^{\flat}}$ is a maximal parahoric subgroup if
  $\Sigma^{\flat}$ is non-empty, otherwise
  $P_{\Sigma^{\flat}} = Q = I$. By Remark~\ref{rmk:54}, the reductive
  quotient $\bar{P}_{\Sigma^{\flat}}$ has the Dynkin diagram
  $\Sigma^{\flat}$ which is the same as
  $\mathrm{SO}(\mathbb{B}_{\Lambda})$ in both odd and even cases. So
  we have the same (partial) flag varieties
  \begin{equation}
    \label{eq:198}
    P_{\Sigma^{\flat}} / I = \mathrm{SO}(\mathbb{B}_{\Lambda})
    / B, \quad P_{\Sigma^{\flat}} / Q =
    \mathrm{SO}(\mathbb{B}_{\Lambda}) / Q'.
  \end{equation}
  The map $\Phi_{\ad} \colon Y_{\Sigma^{\sharp}}(w) \to
  \mathcal{S}^{\circ}_{\Lambda}(\mathbb{F})$ is compatible with their
  embeddings into $P_{\Sigma^{\flat}} / Q$.
  For $M \in \mathcal{S}_{\Lambda}^{\circ}(\mathbb{F})$, by
  Lemma~\ref{lm:3}, its image $\bar{M}$ in the partial flag variety
  $P_{\Sigma^{\flat}} / Q$ satisfies
  \begin{equation}
    \label{eq:202}
    \bar{M} \stackrel{1}{\subset} \bar{M} + \chi(\bar{M}) \stackrel{1}{\subset} \cdots 
    \stackrel{1}{\subset} \bar{M} + \chi(\bar{M}) + \cdots \chi^{i}(\bar{M}) =
    \mathbb{B}_{\Lambda, \mathbb{F}}.
  \end{equation}
  By the description of the fine Deligne-Lusztig varieties, i.e. the
  image of $\phi$, in Example~\ref{eg:2} and \ref{eg:3} and taking
  dual of \eqref{eq:202}, we can see that $\bar{M}$ lies in
  $\mathrm{im}(\phi)$. Hence $\Phi_{\ad}$ is bijective.
  
  For general $j$, if $g \in j \cdot Y_{\Sigma^{\sharp}}(w)$, then
  $j^{-1} g \in Y_{\Sigma^{\sharp}}(w)$. Let $\Lambda'$ be the lattice
  such that
  $Y_{\Sigma^{\sharp}} (w) \cong \mathcal{S}_{\Lambda'}^{\circ}$. Let
  $\Lambda \coloneqq j \cdot \Lambda'$, then
  $j \cdot Y_{\Sigma^{\sharp}}(w) \cong \mathcal{S}_{\Lambda}^{\circ}$
  because $ j \cdot \Xi_{\infty}(M) = \Xi_{\infty}(j \cdot M)$.
\end{proof}

\begin{cor}
  \label{cr:1}
  The map $\Phi_{\ad}$ induces a bijection
  \begin{equation}
    \label{eq:138}
    \overline{j \cdot Y_{\Sigma^{\sharp}}(w)} \longrightarrow
    \mathcal{S}_{\Lambda}(\mathbb{F}),
  \end{equation}
  for each $j \in \mathbb{J}$,
  $w=w_{\Sigma} \in \mathrm{EO}_{\mathrm{cox}}$ and the vertex lattice
  $\Lambda$ corresponding to $j \cdot Y_{\Sigma^{\sharp}}(w)$ via
  Proposition~\ref{pr:12}.
\end{cor}

\begin{proof}
  Let $\Lambda$ be the vertex lattice such that
  $\Phi_{\ad} (j \cdot Y_{\Sigma^{\sharp}}(w) ) =
  \mathcal{S}_{\Lambda}^{\circ}(\mathbb{F})$. Then by
  Proposition~\ref{pr:13},
  $i(\mathbb{J} \cap P_{\tilde{\mathbb{S}}-\Sigma'}) \cap j(\mathbb{J}
  \cap P_{\tilde{\mathbb{S}}-\Sigma}) \neq \emptyset$ if and only if
  $\Lambda$ and $\Lambda'$ are neighbors, where $\Lambda'$ is the
  vertex lattice corresponding to
  $i \cdot Y_{(\Sigma')^{\sharp}}(w_{\Sigma'})$ via
  Proposition~\ref{pr:12}. And
  $(\Sigma')^{\flat} \subset \Sigma^{\flat}$ if and only if
  $(\Sigma')^{\sharp} \supset \Sigma^{\sharp}$, if and only if
  $\Lambda' \subset \Lambda$. So we have
  \begin{equation}
    \label{eq:139}
    \Phi_{\ad}(\overline{j \cdot Y_{\Sigma^{\sharp}}(w)}) =
    \bigcup_{\Lambda' \subset \Lambda} \mathcal{S}_{\Lambda'}^{\circ}(\mathbb{F}).
  \end{equation}
  Then by Proposition~\ref{pr:11} we get the desired result.
\end{proof}

\begin{cor}
  \label{cr:2}
  Let $\Lambda, \Lambda'$ be two vertex lattices, then $\Lambda
  \subset \Lambda'$ if and only if $\mathcal{S}_{\Lambda}(\mathbb{F})
  \subset \mathcal{S}_{\Lambda'}(\mathbb{F})$.
\end{cor}

\begin{proof}
  If
  $\mathcal{S}_{\Lambda}(\mathbb{F}) \subset
  \mathcal{S}_{\Lambda'}(\mathbb{F})$, then
  $\mathcal{S}_{\Lambda}^{\circ}(\mathbb{F}) \subset
  \mathcal{S}_{\Lambda'}(\mathbb{F})$. By Proposition~\ref{pr:12},
  there is a bijection between
  $\mathcal{S}_{\Lambda}^{\circ}(\mathbb{F})$ and a Deligne-Lusztig
  variety, in particular, $\mathcal{S}_{\Lambda}^{\circ}(\mathbb{F})$
  is non-empty. Take
  $M \in \mathcal{S}_{\Lambda}^{\circ}(\mathbb{F}) \subset
  \mathcal{S}_{\Lambda'}(\mathbb{F})$, then $M \subset \Lambda'$
  and by Lemma~\ref{lm:11} we have $\Lambda = \Xi_{\infty}(M) \subset
  \Lambda'$. 
\end{proof}

\begin{remark}
  \label{rmk:11}
  The notations $\mathcal{S}_{\Lambda}(\mathbb{F})$ and
  $\mathcal{S}_{\Lambda}^{\circ}(\mathbb{F})$ imply that they are the
  $\mathbb{F}$-points of the schemes $\mathcal{S}_{\Lambda}$ and
  $\mathcal{S}_{\Lambda}^{\circ}$ which will be defined in
  section~\ref{sec:scheme-theor-struct}.
\end{remark}

\begin{remark}
  \label{rmk:18}
  For each algebraically closed field extension $k$ of $\mathbb{F}$,
  replacing $\mathbb{F}$ by $k$, all results in Section 3 and Section
  4 are true because by the set-up of \cite{GH}, we may work with any
  algebraically closed field extension $k$ of $\mathbb{F}$.
\end{remark}

\section{Scheme-theoretic structure of $\mathcal{N}$}
\label{sec:scheme-theor-struct}

\subsection{The closed and open Bruhat-Tits strata}
\label{sec:closed-open-bruhat}

Let $\Lambda$ be a vertex lattice, we define
\begin{equation}
  \label{eq:141}
  \Lambda^{+}  \coloneqq \Lambda, \quad 
  \Lambda^{-}  \coloneqq \pi \Lambda^{\vee}.
\end{equation}
It is easy to see $\Lambda^{\pm}$ are Dieudonn\'{e} modules in $N$
(recall that $N$ is the rational Dieudonn\'{e} module of
$\mathbb{X}$). Let $X_{\Lambda^{\pm}}$ be the $p$-divisible
$\mathcal{O}_{\breve{F}}$-modules over $\mathbb{F}$ corresponding to
$\Lambda^{\pm}$, together with $\mathcal{O}_{F}$-linear
quasi-isogenies
$\rho_{\Lambda^{\pm}} \colon X_{\Lambda^{\pm}} \to \mathbb{X}$ and
polarizations $\lambda_{\Lambda^{\pm}}$. Note that the form
$\pi^{-1} \langle \thinspace, \thinspace \rangle$ induces a perfect
paring between $\Lambda^{+}$ and $\Lambda^{-}$.

For any $\mathbb{F}$-scheme $S$ and any unitary $p$-divisible group
$(X, \rho_{X}) \in \mathcal{S}(S)$, we define quasi-isogenies:
\begin{align}
  \label{eq:142}
  \rho_{X, \Lambda^{+}} & \colon X_{\bar{S}}
                          \stackrel{\rho_{X}}{\longrightarrow}
                          \mathbb{X}_{\bar{S}}
                          \stackrel{\rho^{-1}_{\Lambda^{+}}}{\longrightarrow}
                          (X_{\Lambda^{+}})_{\bar{S}}, \\
  \rho_{\Lambda^{-}, X} & \colon (X_{\Lambda^{-}})_{\bar{S}}
                          \stackrel{\rho_{\Lambda^{-}}}{\longrightarrow}
                          \mathbb{X}_{\bar{S}}
                          \stackrel{\rho_{X}^{-1}}{\longrightarrow}
                          X_{\bar{S}}. 
\end{align}
By the same reasoning as in Proposition~\ref{pr:3}, we have
\begin{equation}
  \label{eq:145}
  \height(\rho_{X, \Lambda^{+}})  = [\frac{t(\Lambda)}{2}],
  \quad
  \height(\rho_{\Lambda^{-}, X})  = 
  [\frac{t(\Lambda)+1}{2}].
\end{equation}

\begin{definition}
  \label{df:4}
  The subfunctor $\tilde{\mathcal{S}}_{\Lambda}$ is defined as
  \begin{equation}
    \label{eq:143}
    \tilde{\mathcal{S}}_{\Lambda}(S) \coloneqq \{ (X, \rho_{X}) \in
    \mathcal{S}(S) : \rho_{X, \Lambda^{+}} \text{ is an isogeny} \},
  \end{equation}
  for each vertex lattice $\Lambda$ and $\mathbb{F}$-scheme $S$.
\end{definition}

Note that $\rho_{X, \Lambda^{+}}$ is an isogeny if and only if
$\rho_{\Lambda^{-1}, X}$ is an isogeny.

\begin{lemma}
  \label{lm:10}
  The subfunctor $\tilde{\mathcal{S}}_{\Lambda}$ is represented by a
  projective scheme over $\mathbb{F}$ and the monomorphism
  $\tilde{\mathcal{S}}_{\Lambda} \hookrightarrow \mathcal{S}$ is a
  closed immersion.
\end{lemma}

\begin{proof}
  The proof is exactly the same as \cite[Lemma 3.2]{MR2800696}.
\end{proof}

\begin{definition}
  \label{df:8}
  Let
  $\mathcal{S}_{\Lambda} \coloneqq
  (\tilde{\mathcal{S}}_{\Lambda})_{\mathrm{red}}$, we call
  $\mathcal{S}_{\Lambda}$ the \emph{closed Bruhat-Tits stratum}
  associated to $\Lambda$.
\end{definition}

\begin{remark}
  \label{rmk:12}
  The definition of $\mathcal{S}_{\Lambda}$ coincides with
  Definition~\ref{df:1} on $k$-points, in the spirit of
  Remark~\ref{rmk:18}, for any algebraically closed field extension
  $k$ of $\mathbb{F}$.
\end{remark}

If $\Lambda, \Lambda'$ are two vertex lattices such that $\Lambda
\subset \Lambda'$, by Dieudonn\'{e} theory, the corresponding
quasi-isogeny $X_{\Lambda} \to X_{\Lambda'}$ is an isogeny, so we have
$\mathcal{S}_{\Lambda} \subset \mathcal{S}_{\Lambda'}$.

\begin{definition}
  \label{df:5}
  The locally closed subscheme $\mathcal{S}_{\Lambda}^{\circ}$ is
  defined as
  \begin{equation}
    \label{eq:144}
    \mathcal{S}_{\Lambda}^{\circ} \coloneqq
    \mathcal{S}_{\Lambda} \backslash \bigcup_{\Lambda' \subsetneq
      \Lambda} \mathcal{S}_{\Lambda'},
  \end{equation}
  for each vertex lattice $\Lambda$. Then
  $\mathcal{S}_{\Lambda}^{\circ}$ is an open subscheme of
  $\mathcal{S}_{\Lambda}$. We call $\mathcal{S}_{\Lambda}^{\circ}$ the
  \emph{open Bruhat-Tits stratum} associated to $\Lambda$.
\end{definition}

By definition, we have
\begin{equation}
  \label{eq:221}
  \mathcal{S}_{\Lambda} = \biguplus_{\Lambda' \subset \Lambda} 
  \mathcal{S}_{\Lambda'}^{\circ}.
\end{equation}

\begin{remark}
  \label{rmk:1222}
  The definition of $\mathcal{S}_{\Lambda}^{\circ}$ coincides with
  Definition~\ref{df:3} on $k$-valued points, in the spirit of
  Remark~\ref{rmk:18}, for any algebraically closed field extension
  $k$ of $\mathbb{F}$ by Proposition~\ref{pr:11}.
\end{remark}

\subsection{An $A$-windows-theory interlude}
\label{sec:wind-theory-interl}

We need some results about Zink's windows theory for formal
$p$-divisible groups. The reader is referred to \cite{MR1827031} for
all relevant concepts.
 
Let $k$ be a field of characteristic $p$, $A$ the Cohen subring of
$W(k)$ (cf. \cite[IX $\mathsection 2$ Definition 2]{MR2284892})

\begin{lemma}
  \label{lm:19}
  Let $Y$ be a $p$-divisible group over $k$ of height $2 d$ and
  dimension $d$, $(M_{Y}, M_{Y, 1}, \Upsilon_{Y})$ its
  $A$-window. Then giving a $p$-divisible group $X$ over $k$ of height
  $2 d$ and dimension $d$, together with an isogeny
  $\rho \colon X \to Y$, is equivalent to giving an $A$-submodule $M$
  of $M_{Y}$ such that $M$ is $\Upsilon_{Y}$-stable and
  $p M \stackrel{2 d}{\subset} M$.
\end{lemma}

\begin{proof}
  The proof is straightforward, the reader is referred to \cite[Lemma
  5.2.11]{wu-thesis}.
\end{proof} 

Now let us consider the $A$-windows associated to unitary
$p$-divisible groups. Let $k \supset \mathbb{F}$ be a field extension,
$A$ the Cohen subring of $W(k)$ which is also an
$\mathcal{O}_{L}$-algebra. Then
$(\mathcal{O}_{L}, p \mathcal{O}_{L}, \sigma)$ is a frame over
$\mathbb{F}$ and $(A, p A, \sigma_{A})$ is a frame over $k$. The
inclusion $\mathcal{O}_{L} \subset A$ induces a morphism of frames
\begin{equation}
  \label{eq:265}
  (\mathcal{O}_{L}, p \mathcal{O}_{L}, \sigma) \longrightarrow (A, p
  A, \sigma_{A}).
\end{equation}
By abuse of notation, let $\sigma$ denote $\sigma_{A}$. Recall that we
fix a supersingular unitary $p$-divisible group
$(\mathbb{X}, \iota_{\mathbb{X}}, \lambda_{\mathbb{X}})$ over
$\mathbb{F}$ of signature $(1, n-1)$ in
Section~\ref{sec:moduli-space-p}, then the $A$-window of the
underlying $p$-divisible group $\mathbb{X}$ over $k$ is the base
change of the Dieudonn\'{e} module
$(\mathbb{M}, \mathcal{F}, \mathcal{V})$ via the morphism of
frames~\eqref{eq:265}. More precisely, let
$(\mathbb{M}_{A}, \mathbb{M}_{A, 1}, \Upsilon)$ be the $A$-window of
$\mathbb{X} \otimes k$, then by \cite[Theorem 4]{MR1827031}, we have
$\mathbb{M}_{A} = \mathbb{M} \otimes_{\mathcal{O}_{L}} A$,
$\mathbb{M}_{A, 1}$ is the submodule of $\mathbb{M}_{A}$ generated by
$\mathcal{V} \mathbb{M} \otimes_{\mathcal{O}_{L}} A$ and
$\mathbb{M} \otimes_{\mathcal{O}_{L}} p A$, and
$\Upsilon = \mathcal{F} \otimes \sigma_{A}$. Note that
$\mathbb{M} \otimes_{\mathcal{O}_{L}} p A = p \mathbb{M}
\otimes_{\mathcal{O}_{L}} A \subset \mathcal{V} \mathbb{M}
\otimes_{\mathcal{O}_{L}} A$, so we have
$\mathbb{M}_{A, 1} = \mathcal{V} \mathbb{M} \otimes_{\mathcal{O}_{L}}
A$. Let $(N_{A}, \Upsilon)$ be the rational $A$-window, i.e.
$N_{A} = \mathbb{M}_{A} \otimes_{\mathcal{O}_{L}} \mathrm{Frac}(A)$,
together with the $\mathcal{O}_{F}$-action $\iota_{\mathbb{X}}$ and
the non-degenerate alternating form
$\langle \thinspace, \thinspace \rangle$ induced by the polarization
$\lambda_{\mathbb{X}} \otimes k$. For any $x, y \in N_{A}$, we have
\begin{equation}
  \label{eq:266}
  \langle \Upsilon(x), \Upsilon(y) \rangle = \langle x, y
  \rangle ^{\sigma},
\end{equation}
and
\begin{equation}
  \label{eq:267}
  \langle \iota(\pi) x, y \rangle =
  \langle x, \iota(\bar{\pi}) y \rangle.
\end{equation}
Henceforth, we write $\pi$ instead of $\iota(\pi)$ to lighten the
notations. The $\pi$-action defines an $A[\pi]$-module structure on
$\mathbb{M} \otimes_{\mathcal{O}_{L}} A$. Let $\Lambda$ be a vertex
lattice. Then
$(\Lambda^{\pm} \otimes A, \mathcal{V} \Lambda^{\pm} \otimes A,
\Upsilon)$ are the $A$-windows of the $p$-divisible groups
$X_{\Lambda^{\pm}} \otimes k$. We will write $\Lambda^{\pm}_{A}$
instead of $\Lambda^{\pm} \otimes A$, and write
$\mathcal{V} \Lambda^{\pm}_{A}$ instead of
$\mathcal{V} \Lambda^{\pm} \otimes A$ for short.

By Lemma~\ref{lm:19}, we have the following windows description of
$\mathcal{S}_{\Lambda}(k)$.

\begin{prop}
  \label{pr:24}
  Via $A$-windows theory, $\mathcal{S}_{\Lambda}(k)$ can be identified
  with the set of $A[\pi]$-lattices $M$ in $N_{A}$ satisfying the
  following conditions:
  \begin{enumerate}
  \item $M$ is $\Upsilon$-stable;
  \item
    $M \stackrel{n-1}{ \subset } M^{\vee} \stackrel{1}{\subset}
    \pi^{-1} M$ if $n$ is odd, and $M^{\vee} = \pi^{-1} M$ if $n$ is
    even;
  \item $p M \stackrel{n}{\subset } M_{1} \stackrel{n}{\subset } M$;
  \item $M_{1} \stackrel{\leq 1}{\subset } M_{1} + \pi M$;
  \item if $n$ is even, $M_{1} \stackrel{1}{\subset } M_{1} + \pi M$;
  \item $M \subset \Lambda_{A}^{+} $;
  \end{enumerate}
  where
  $M_{1} \coloneqq \ker(M \to \Lambda^{+} / \mathcal{V} \Lambda^{+}_{A}
  )$.
\end{prop}

\subsection{The Bruhat-Tits strata as Deligne-Lusztig varieties}
\label{sec:bruhat-tits-strata}

Let $T$ be a scheme over $\mathbb{F}$, $(X, \rho) \in \mathcal{S}(T)$
a unitary $p$-divisible group. Let $D(X)$ be the Lie algebra of the
universal vector extension of $X$ (cf. \cite[Chapter IV, Definition
1.12]{MR0347836}), then the functor
\begin{align}
  \label{eq:207}
  (p\text{-divisible groups over } T) &
  \longrightarrow (\text{locally free } 
  \mathcal{O}_{T} \text{-modules}), \\
  X & \longmapsto D(X),
      \nonumber
\end{align}
commutes with an arbitrary base change $T' \to T$. When $T =
\mathrm{Spec}(k)$ for an algebraically closed field extension $k$ of
$\mathbb{F}$, we have $D(X) \cong M(X) / p M(X)$ canonically, where
$M(X)$ is the Dieudonn\'{e} module of $X$.

\begin{lemma}
  \label{lm:13}
  Let $\rho_{i} \colon X \to Y_{i}$, for $i = 1, 2$, be two isogenies
  of naive unitary $p$-divisible groups (of any signature) over $T$,
  such that $\ker(\rho_{1}) \subset \ker(\rho_{2}) \subset X[\pi]$,
  then both $\ker(D(\rho_{1}))$ and $\ker(D(\rho_{2}))$ are locally
  free $\mathcal{O}_{T}$-modules and $\ker(D(\rho_{1}))$ is a locally
  direct summand of $\ker(D(\rho_{2}))$.
\end{lemma}

\begin{proof}
  Note that by definition $X$ is endowed with an
  $\mathcal{O}_{F}$-action, hence the proof is exactly the same as
  {\cite[Corollary 3.7]{MR2800696}} replacing $p$ by $\pi$.
\end{proof}

Using Lemma~\ref{lm:13}, we can construct a morphism from
$\tilde{\mathcal{S}}_{\Lambda}$ to the partial flag variety
$\mathrm{SO}(\mathbb{B}_{\Lambda}) / Q'$ defined in
section~\ref{sec:set-structure-bruhat}. Let
$(X, \rho) \in \tilde{\mathcal{S}}_{\Lambda}(R)$ for an
$\mathbb{F}$-algebra $R$ and a vertex lattice $\Lambda$, we have
isogenies
\begin{equation}
  \label{eq:208}
  (X_{\Lambda^{-}})_{\bar{R}}
  \stackrel{\rho_{\Lambda^{-}}}{\xrightarrow{\hspace*{1cm}}} 
  X_{\bar{R}}
  \stackrel{\rho_{\Lambda^{+}}}{\xrightarrow{\hspace*{1cm}}}
  (X_{\Lambda^{+}})_{\bar{R}}.
\end{equation}
where
$\rho_{\Lambda^{-}} = \rho_{\Lambda^{-}, X} \otimes \mathrm{id}_{R}$
by abuse of notation and similarly for $\rho_{\Lambda^{+}}$. The
composition
$\rho_{\Lambda} \coloneqq \rho_{\Lambda^{+}} \circ \rho_{\Lambda^{-}}$
corresponds to the isogeny
$(X_{\Lambda^{-}})_{\bar{R}} \to (X_{\Lambda^{+}})_{\bar{R}}$ induced
by the inclusion $\Lambda^{-} \subset \Lambda^{+}$. Then we have
$\ker(\rho_{\Lambda^{-}}) \subset \ker(\rho_{\Lambda}) \subset
X_{\Lambda^{-}}[\pi]$. Note that
$\ker(D(\rho_{\Lambda})) = \mathbb{B}_{\Lambda, R} \coloneqq
\mathbb{B}_{\Lambda} \otimes R$, and when $R = \mathrm{Spec}(k)$ for
an algebraically closed field $k$,
$\ker(D(\rho_{\Lambda^{-}})) = M(X) / \pi \Lambda^{\vee}$.

Recall that for any $\mathbb{F}$-algebra $R$, the partial flag variety
$\mathrm{SO}(\mathbb{B}_{\Lambda}) / Q'$ has the following description
as a functor
\begin{equation}
  \label{eq:209}
  (\mathrm{SO}(\mathbb{B}_{\Lambda}) / Q') (R) =  
  \left\{
    \begin{array}{ll}
      U \subset \mathbb{B}_{\Lambda, R} \\ 
      \text{a direct summand}
    \end{array}
    \left|
      \begin{array}{ll}
        U \subset U^{\bot}, \\
        \mathrm{rank}_{R}(U) = [\frac{t(\Lambda)}{2}], \\
        U \text{ lies in the } \mathrm{SO}(\mathbb{B}_{\Lambda})
        \text{-orbit}\\
        \text{corresponding to } Q'
      \end{array}
    \right.
  \right\}.
\end{equation}
For the orthogonal Grassmannian
$\mathrm{Grass}(\mathbb{B}_{\Lambda})$, we have
\begin{equation}
  \label{eq:20999}
  \mathrm{Grass}(\mathbb{B}_{\Lambda}) (R) =  
  \left\{
    \begin{array}{ll}
      U \subset \mathbb{B}_{\Lambda, R} \\ 
      \text{a direct summand}
    \end{array}
  \right.
  \left|
    \begin{array}{ll}
      U \subset U^{\bot}, \\
      \mathrm{rank}_{R}(U) = [\frac{t(\Lambda)}{2}] 
    \end{array}
  \right\}.
\end{equation}
Let $E(X) \coloneqq \ker(D(\rho_{\Lambda^{-}}))$ which is of rank
$\height(\rho_{\Lambda^{-}}) = [\frac{t(\Lambda)+1}{2}]$, then sending
$(X, \rho)$ to $E(X)^{\bot}$ defines a map
\begin{align}
  \label{eq:291}
  \tilde{f} \colon \tilde{\mathcal{S}}_{\Lambda}(R)
  & \longrightarrow
    \mathrm{Grass}(\mathbb{B}_{\Lambda})(R), \\
  (X, \rho)
  & \longmapsto E(X)^{\bot}.
    \nonumber
\end{align}
In summary we have a morphism
$\tilde{\mathcal{S}}_{\Lambda} \to
\mathrm{Grass}(\mathbb{B}_{\Lambda})$, which induces a morphism
\begin{equation}
  \label{eq:292}
  f \colon \mathcal{S}_{\Lambda} \longrightarrow
  \mathrm{Grass}(\mathbb{B}_{\Lambda}).
\end{equation}
Note that by
Remark~\ref{rmk:13}, for any algebraically field extension $k$ of
$\mathbb{F}$, we have
\begin{equation}
  \label{eq:210}
  \mathcal{S}_{\Lambda}(k)
  \hookrightarrow (\mathrm{SO}(\mathbb{B}_{\Lambda})/Q')(k) \subset
  \mathrm{Grass} (\mathbb{B}_{\Lambda})(k),
\end{equation}
i.e. the image of $\mathcal{S}_{\Lambda}$ lies in
$\mathrm{SO}(\mathbb{B}_{\Lambda}) / Q'$ because
$\mathcal{S}_{\Lambda}$ is reduced.

\begin{lemma}
  \label{lm:14}
  The morphism
  $f \colon \mathcal{S}_{\Lambda} \to
  \mathrm{SO}(\mathbb{B}_{\Lambda})/Q'$ is a closed immersion. In
  particular, taking closure of $\mathcal{S}_{\Lambda}^{\circ}$ in
  $\mathcal{S}$ is the same as taking closure in
  $\mathrm{SO}(\mathbb{B}_{\Lambda}) / Q'$.
\end{lemma}

The proof is trivial.
 
\begin{lemma}
  \label{lm:20}
  The morphism $f$ induces a morphism
  $f \colon \mathcal{S}_{\Lambda} \to \overline{j \cdot
    Y_{\Sigma^{\sharp}(w)}}$, where $j \in \mathbb{J}$ and
  $w = w_{\Sigma} \in \mathrm{EO}_{\mathrm{cox}}$ corresponding to
  $\Lambda$ via Proposition~\ref{pr:12}.
\end{lemma}

\begin{proof}
  For any algebraically closed field extension $k$ of $\mathbb{F}$, we
  have
  $f \colon \mathcal{S}_{\Lambda}(k) \to (\overline{j\cdot
    Y_{\Sigma^{\sharp}}(w)}) (k)$ by Corollary~\ref{cr:1}, since
  $\mathcal{S}_{\Lambda}$ is reduced, we prove the claim.
\end{proof}

\begin{lemma}
  \label{lm:16}
  Let $k$ be a field extension of $\mathbb{F}$ (not necessarily
  algebraically closed), then the morphism $f$ induces a bijection
  \begin{equation}
    \label{eq:253}
    \mathcal{S}_{\Lambda}(k) \longrightarrow
    \overline{j \cdot Y_{\Sigma^{\sharp}}(w)}(k).
  \end{equation}
\end{lemma}

\begin{proof}
  The injectivity of $f(k)$ follows from that $f(\bar{k})$ is
  bijective by Corollary~\ref{cr:1} and Remark~\ref{rmk:18}.

  Let us prove the surjectivity of $f(k)$. Let $A$ be the Cohen
  subring of $W(k)$, then $\mathcal{S}_{\Lambda}(k)$ can be described
  as the set of all the $A[\pi]$-lattices in $N_{A}$ satisfying all
  the conditions in Proposition~\ref{pr:24}. Let
  $U_{0} \in \overline{j \cdot Y_{\Sigma^{\sharp}}(w)} (k)$ be a
  maximal isotropic subspace of $\mathbb{B}_{\Lambda, k}$, then
  $U_{0}^{\bot}$ gives rise to an $A[\pi]$-module $M$ such that
  \begin{equation}
    \label{eq:268}
    \Lambda^{-}_{A} \subset \pi M^{\vee} \stackrel{}{\subset }
    M \subset  \Lambda^{+}_{A} .
  \end{equation}
  To prove the surjectivity of $f(k)$, it only needs to show that the
  $A[\pi]$-module $M$ lies in $\mathcal{S}_{\Lambda}(k)$, i.e. $M$
  satisfies all the conditions in Proposition~\ref{pr:24}. We leave it
  to the reader.
\end{proof}

\begin{prop}
  \label{pr:17}
  The morphism $f \colon \mathcal{S}_{\Lambda} \to \overline{j \cdot
    Y_{\Sigma^{\sharp}}(w)}$ is an isomorphism.
\end{prop}

\begin{proof}
  The proof is exactly the same as \cite[Theorem 4.8]{MR2800696},
  except that we didn't compute the dimension of the tangent space of
  $\mathcal{S}_{\Lambda}$ at every $k$-valued point, which seems not
  necessary.
\end{proof}

\begin{cor}
  \label{cr:6}
  The morphism $f$ induces an isomorphism
  $\mathcal{S}_{\Lambda}^{\circ} \to j \cdot
  Y_{\Sigma^{\sharp}}(w)$. In particular, the locally closed subscheme
  $\mathcal{S}^{\circ}_{\Lambda}$ is smooth of dimension $\ell(w) =
  [\frac{t(\Lambda) - 1}{2}]$.
\end{cor}

\begin{cor}
  \label{cr:5}
  The closure $\overline{\mathcal{S}_{\Lambda}^{\circ}}$ of
  $\mathcal{S}_{\Lambda}^{\circ}$ in $\mathcal{S}$ is
  $\mathcal{S}_{\Lambda}$.
\end{cor}
 
By Example~\ref{eg:7} and \ref{eg:8}, we have the following corollary.

\begin{cor}
  \label{cr:3}
  The closed subscheme $\mathcal{S}_{\Lambda}$ of $\mathcal{S}$ is
  projective and normal of dimension $\ell(w_{\Sigma}) =
  [\frac{t(\Lambda) - 1}{2}] $. When $n$ is odd,
  $\mathcal{S}_{\Lambda}$ has isolated singularities; when $n$ is
  even, $\mathcal{S}_{\Lambda}$ is smooth.
\end{cor}

\subsection{The Bruhat-Tits stratification}
\label{sec:bruh-tits-strat}

\begin{theorem}
  \label{thr:1}
  Let $\Lambda$ and $\Lambda'$ be two vertex lattices.
  \begin{enumerate}
  \item We have $\Lambda \subset \Lambda'$ if and only if
    $\mathcal{S}_{\Lambda} \subset \mathcal{S}_{\Lambda'}$.
  \item We have
    \begin{equation}
      \label{eq:222}
      \mathcal{S}_{\Lambda} \cap \mathcal{S}_{\Lambda'}
      = 
      \left\{
        \begin{array}{ll}
          \mathcal{S}_{\Lambda \cap \Lambda'}, & \text{ if } \Lambda
                                                \cap \Lambda'
                                                \text{ is a vertex lattice again,} \\
\emptyset, & \text{ otherwise.}
        \end{array}
      \right.
    \end{equation}
  \item Recall that $\mathcal{B}$ is the set of vertex lattices, then
    we have
    \begin{equation}
      \label{eq:223}
      \mathcal{S} = \bigcup_{\Lambda \in \mathcal{B}} \mathcal{S}_{\Lambda},
    \end{equation}
    and each closed Bruhat-Tits stratum $\mathcal{S}_{\Lambda}$ is
    projective and normal of dimension $[\frac{t(\Lambda) - 1}{2}]$,
    with isolated singularities when $n$ is odd, is smooth when $n$ is
    even.
  \end{enumerate}
\end{theorem}

\begin{proof}
  It follows from Proposition~\ref{pr:8} and Corollary~\ref{cr:2} that
  part $1$ and $2$ are true. Part $3$ follows from
  Corollary~\ref{cr:3}. 
\end{proof}

\begin{theorem}
  \label{thr:bt}
  \leavevmode
  \begin{enumerate}
  \item There is a stratification, which is called the Bruhat-Tits
    stratification, of $\mathcal{S}$ by locally closed subschemes
    \begin{equation}
      \label{eq:224}
      \mathcal{S} = \biguplus_{\Lambda \in \mathcal{B}} 
      \mathcal{S}_{\Lambda}^{\circ},
    \end{equation}
    and each stratum is isomorphic to the Deligne-Lusztig variety
    associated to the orthogonal group
    $\mathrm{SO}(\mathbb{B}_{\Lambda})$ and a $\sigma$-Coxeter
    element. The closure of each stratum
    $\mathcal{S}_{\Lambda}^{\circ}$ in $\mathcal{S}$ is given by
    \begin{equation}
      \label{eq:225}
      \overline{\mathcal{S}_{\Lambda}^{\circ}} = \biguplus_{\Lambda'
        \subset \Lambda} \mathcal{S}_{\Lambda'}^{\circ} 
      = \mathcal{S}_{\Lambda}.
    \end{equation}
  \item The scheme $\mathcal{S}$ is geometrically connected of pure
    dimension $[\frac{n-1}{2}]$. The irreducible components of
    $\mathcal{S}$ are those $\mathcal{S}_{\Lambda}$ with
    $t(\Lambda) = n$.
  \end{enumerate}
\end{theorem}

\begin{proof}
  \leavevmode
  \begin{enumerate}
  \item The stratification follows from \eqref{eq:221} and
    part $3$ of Theorem~\ref{thr:1}.
  \item For a vertex lattice, the form $\psi$ in
    section~\ref{sec:local-pel-datum} defines a non-degenerate
    symplectic form on the quotient space $\Lambda^{\sharp} / \Lambda$
    (cf. \cite[Lemma 6.4]{MR3175176}). Then a vertex lattices
    $\Lambda'$ such that $\Lambda' \supset \Lambda$ corresponds to an
    isotropic subspace of $\Lambda^{\sharp} / \Lambda$. In particular,
    $\Lambda$ is contained in a maximal type vertex lattice. By the
    part $1$ of Theorem~\ref{thr:1}, $\mathcal{S}_{\Lambda}$ is an
    irreducible component of $\mathcal{S}$ if $t(\Lambda) = n$. The
    simplicial complex $\mathcal{B}$ is connected, hence $\mathcal{S}$
    is connected of pure dimension $[\frac{n-1}{2}]$.
  \end{enumerate}
\end{proof}

\section{The supersingular locus of the unitary Shimura varieties}
\label{sec:supers-locus-unit}

\subsection{The integral model}
\label{sec:integral-model}

We start with the ramified unitary PEL datum of signature $(1, n-1)$
(cf. \cite[1.1]{MR2516305} or \cite[5.1]{MR3449174}). For the
definition of the general PEL datum, we refer to
\cite[2.1]{MR3449174}.

Let $E$ be an imaginary quadratic field extension of $\mathbb{Q}$ with
a fixed embedding $\gamma_{0} \colon E \hookrightarrow
\mathbb{C}$. Let $\bar{} \in \mathrm{Gal}(E/\mathbb{Q})$ be the unique
non-trivial automorphism. Then $\gamma_{0}$ and
$\gamma_{1} \coloneqq \gamma_{0} \circ \bar{}$ give rise to all the
embeddings of $E$ into $\mathbb{C}$. Let $W = E^{n}$ be an
$n$-dimensional vector space over $E$, where $n \geq 3$, together with
a hermitian form $\varphi$. We fix an element $\epsilon \in E$ such
that $\bar{\epsilon} = - \epsilon$, then the form
$\epsilon \cdot \varphi$ is a skew hermitian form on $W$. Furthermore,
we assume that the hermitian form $\varphi$ is of signature $(1, n-1)$
in the following sense: there exists a $\mathbb{C}$-basis of $W
\otimes_{E, \gamma_{0}} \mathbb{C}$ such that the matrix of $\varphi$ is
\begin{equation}
  \label{eq:226}
  H \coloneqq \mathrm{diag}(-1, 1, \ldots, 1).
\end{equation}
Note that we have an $\mathbb{R}$-isomorphism
$W \otimes_{\mathbb{Q}} \mathbb{R} \cong W \otimes_{E, \gamma_{0}}
\mathbb{C}$. Therefore the matrix $\sqrt{-1} \cdot H$ defines an
$\mathbb{R}$-endomorphism of $W \otimes \mathbb{R}$ satisfying
$(\sqrt{-1} \cdot H)^{2} = -\mathrm{id}$ and hence a complex structure
of $W \otimes \mathbb{R}$.

The hermitian form $\varphi$ defines a $\mathbb{Q}$-linear symplectic
form
$\langle \thinspace, \thinspace \rangle \colon W \times W \to
\mathbb{Q}$ by
$\langle \thinspace, \thinspace \rangle \coloneqq
\mathrm{Tr}_{E/\mathbb{Q}} (\epsilon \cdot \varphi(\thinspace,
\thinspace)) $. The form $ \langle v, \sqrt{-1} \cdot H w \rangle \ $,
for $v, w \in W \otimes \mathbb{R}$, is $\mathbb{R}$-symmetric, and if
it is not positive definite, we replace $\epsilon$ by $-\epsilon$
which will guarantee the positive definiteness.

Let $p$ be an odd prime which ramifies in $E$. Let $v$ be the place
above $p$, $E_{v}$ the completion of $E$ at $v$ with the ring of
integers $\mathcal{O}_{v}$. Let $\pi$ be a uniformizer of $E_{v}$ such
that $\bar{\pi} = -\pi$. We assume that the hermitian space
$(W \otimes_{E} E_{v}, \varphi)$ is split. We can define standard
lattices $\{ \Lambda_{i} \}_{i \in \mathbb{Z}}$ in the same manner as
in \ref{sec:lattice-model-bruhat}.

For any $\mathbb{Q}$-algebra $R$, let
\begin{equation}
  \label{eq:227}
  \mathbb{G}(R) \coloneqq \left\{ g \in \mathrm{GL}_{E \otimes R} (W
  \otimes R) 
  \left|
    \begin{array}{ll}
      \exists c = c(g) \text{ such that } \forall v, w \in W \\
      \langle g(v), g(w) \rangle = c \langle v, w \rangle
    \end{array}
    \right.
  \right\}.
\end{equation}
Then $\mathbb{G}$ is a reductive group over $\mathbb{Q}$. Sending
$\sqrt{-1}$ to $\sqrt{-1} \cdot H$ defines a homomorphism
\begin{equation}
  \label{eq:228}
  h \colon \mathrm{Res}_{\mathbb{C} / \mathbb{R}} (\mathbb{G}_{m,
    \mathbb{C}}) \longrightarrow \mathbb{G}_{\mathbb{R}}.
\end{equation}
Then the $\mathbb{Q}$-reductive group $\mathbb{G}$ and the
$\mathbb{G}(\mathbb{R})$-conjugacy class $X$ of $h$ define a Shimura
datum, hence the Shimura variety $\mathrm{Sh}(\mathbb{G}, h)$ over the
reflex field $E$. Let $C = \prod_{w} C_{w}$ be an open compact
subgroup of $\mathbb{G}(\mathbb{A}_{f})$ with
$C_{w} \subset \mathbb{G}(\mathbb{Q}_{w})$. Then the Shimura variety
$\mathrm{Sh}_{C}(\mathbb{G}, h)$ is a quasi-projective variety over
$E$ whose $\mathbb{C}$-valued points can be identified with
\begin{equation}
  \label{eq:229}
  \mathbb{G}(\mathbb{Q}) \backslash 
  ( X \times (\mathbb{G}(\mathbb{A}_{f}) /C )).
\end{equation}
We assume that the subgroup
$C^{p} \coloneqq \prod_{w \neq p} C_{w} \subset
\mathbb{G}(\mathbb{A}_{f}^{p})$ is sufficiently small, i.e. the
subgroup $C^{p}$ is contained in the principal congruence subgroup of
level $N \geq 3$, where $N$ is coprime to the discriminant of $E$. We
also assume that $C_{p}$ is the parahoric subgroup of
$\mathbb{G}(\mathbb{Q}_{p})$ stabilizing the lattice $\Lambda_{m}$.

Now we define the integral model of $\mathrm{Sh}_{C}(\mathbb{G}, h)$
over $E_{v}$ following \cite[Chapter 6]{MR1393439}. For a fixed base
scheme $S$, let $\mathrm{AV}(S)$ be the category of abelian
$\mathcal{O}_{E}$-varieties up to isogeny of order prime to $p$ over
$S$ (cf. \cite[6.3]{MR1393439}).

\begin{definition}
  \label{df:11}
  The naive moduli functor $\mathcal{A}_{C^{p}}^{\mathrm{naive}}$ over
  $\mathcal{O}_{E_{v}}$ is a set-valued functor:
  \begin{align}
    \label{eq:232}
    (\mathcal{O}_{E_{v}}) \text{-schemes} & 
                                            \longrightarrow (\text{Sets}), \\
    S & \longmapsto \text{isomorphism classes of } (A, \iota,
        \bar{\lambda}, \bar{\eta}),
        \nonumber
  \end{align}
  where $(A, \iota) \in \mathrm{AV}(S)$, $\bar{\lambda}$ is a
  $\mathbb{Q}$-homogeneous polarization of $(A, \iota)$ which contains
  a polarization $\lambda \colon A \to A^{\vee}$ such that
  \begin{itemize}
  \item if $n$ is odd, $\ker(\lambda) \subset A[\iota(\pi)]$ is of
    height $n-1$,
  \item if $n$ is even, $\ker(\lambda) = A[\iota(\pi)]$;
  \end{itemize}
  and $\bar{\eta}$ is a $C^{p}$-level structure
  \begin{equation}
    \label{eq:233}
    \bar{\eta} \colon \mathrm{H}_{1}(A, \mathbb{A}^{p}_{f}) \cong W
    \otimes \mathbb{A}^{p}_{f} \mod C^{p}.
  \end{equation}
  Furthermore, the pair $(A, \iota)$ is required to satisfy the
  determinant condition:
  \begin{equation}
    \label{eq:110}
    \mathrm{det}_{\mathcal{O}_{S}} ( \iota(a) | \mathrm{Lie}_{S}(A))
    = (T_{0} + T_{1} \pi) ( T_{0} + T_{1} \bar{\pi})^{n-1} \in
    \mathcal{O}_{E_{v}}[T_{0}, T_{1}]
  \end{equation}
  for all $a \in \mathcal{O}_{E_{v}}$. Then the functor
  $\mathcal{A}_{C^{p}}^{\mathrm{naive}}$ is represented by a
  quasi-projective scheme over $\mathcal{O}_{E_{v}}$, which is denoted
  by $\mathcal{A}^{\mathrm{naive}}_{C^{p}}$ by abuse of notation.
\end{definition}

The scheme $\mathcal{A}^{\mathrm{naive}}_{C^{p}}$ is not flat by
\cite[Proposition 3.8]{MR1752014}. 

\begin{definition}
  \label{df:12}
  The subfunctor $\mathcal{A}^{e}$ of
  $\mathcal{A}^{\mathrm{naive}}_{C^{p}}$ is defined by requiring that
  the quadruple
  $(A, \iota, \bar{\lambda}, \bar{\eta}) \in \mathcal{A}^{e}(S)$
  satisfy the following condition(s):
  \begin{enumerate}
  \item (Wedge condition.) For each $a \in \mathcal{O}_{E_{v}}$, the
    homomorphisms
    \begin{align}
      \label{eq:115}
      \wedge^{n} (\iota(a) - a) \colon &
                                         \wedge^{n} \mathrm{Lie}(A) \longrightarrow \wedge^{n}
                                         \mathrm{Lie}(A), \\
      \wedge^{2} (\iota(a) - \bar{a}) \colon &
                                               \wedge^{2} \mathrm{Lie}(A) \longrightarrow \wedge^{2}
                                               \mathrm{Lie}(A),
    \end{align}
    are both equal to zero.
  \item When $n$ is even, the extra Spin condition is assumed:
    $\iota(\pi) | \mathrm{Lie}(A_{s})$ non-vanishing for all
    $s \in S$.
  \end{enumerate}
\end{definition}

\begin{definition}
  \label{df:13}
  The \emph{honest integral model} $\mathcal{A}$ is defined as the
  flat closure of $\mathcal{A}^{\mathrm{naive}}_{C^{p}}$ in its
  generic fiber.
\end{definition}

\begin{prop}[Smithling]
  \label{pr:20}
  The functor $\mathcal{A}^{e}$ is represented by a closed subscheme
  of $\mathcal{A}^{\mathrm{naive}}_{C^{p}}$ over
  $\mathcal{O}_{E_{v}}$, which is topologically flat and of dimension
  $n-1$. Furthermore, when $n$ is even, $\mathcal{A}^{e}$ is flat over
  $\mathcal{O}_{E_{v}}$, in other words,
  $\mathcal{A}^{e} = \mathcal{A}$.
\end{prop}

\begin{proof}
  Note that $\mathcal{A}^{\mathrm{naive}}_{C^{p}}$, $\mathcal{A}^{e}$
  and $\mathcal{A}$ sit inside the usual local model diagram by
  \cite[Theorem 2.2]{MR1752014}. Then similar to the proof of
  Proposition~\ref{pr:9}, the proposition follows from the property of
  the local model.
\end{proof}

Let $\mathcal{A}_{\mathbb{F}}$ (resp. $\mathcal{A}^{e}_{\mathbb{F}}$)
be the special fiber of $\mathcal{A}$ (resp. $\mathcal{A}^{e}$), then
by Proposition~\ref{pr:20}, we have $\mathcal{A}_{\mathbb{F},
  \mathrm{red}} = \mathcal{A}^{e}_{\mathbb{F}, red}$.

\subsection{The supersingular locus}
\label{sec:supersingular-locus}

Let $\mathcal{A}^{e, \mathrm{ss}}_{\mathbb{F}}$
(resp. $\mathcal{A}^{\mathrm{ss}}_{\mathbb{F}}$) be the supersingular
locus of $\mathcal{A}^{e}_{\mathbb{F}}$
(resp. $\mathcal{A}^{\mathrm{ss}}_{\mathbb{F}}$), then
$\mathcal{A}^{\mathrm{ss}}_{\mathbb{F}} = \mathcal{A}^{e,
  \mathrm{ss}}_{\mathbb{F}} $ because by definition the supersingular
locus is endowed with the closed reduced subscheme structure
(cf. \cite[Theorem 6.27]{MR1393439}). 

Similarly to the naive case, we have the $p$-adic uniformization
theorem.

\begin{theorem}[{\cite[Theorem 6.30]{MR1393439} \&
    \cite[6.4]{MR2800696}}]
  \label{thr:uniformization}
  Let
  $(A_{0}, \iota_{0}, \bar{\lambda}_{0}, \bar{\eta}_{0}) \in
  \mathcal{A}^{e}(\mathbb{F})$ be a supersingular abelian variety,
  together with its corresponding Rapoport-Zink space
  $\mathcal{N}^{e}$. Then the uniformization morphism given by
  $(A_{0}, \iota_{0}, \bar{\lambda}_{0}, \bar{\eta}_{0})$
  \begin{equation}
    \label{eq:239}
    \Theta \colon \mathbb{I}(\mathbb{Q}) \backslash
    \mathcal{N}^{e}_{\mathrm{red}}
    \times \mathbb{G}(\mathbb{A}^{p}_{f}) / C^{p}
    \longrightarrow \mathcal{A}^{\mathrm{ss}}_{\mathbb{F}}
  \end{equation}
  is an isomorphism, $\mathbb{I}$ is the group of
  $\mathcal{O}_{E_{v}}$-linear quasi-isogenies in $\mathrm{End}(A_{0})
  \otimes \mathbb{Q}$ which respect the polarizations
  $\bar{\lambda}_{0}$. And the source of the uniformization morphism
  is a finite disjoint sum
  \begin{equation}
    \label{eq:240}
    \coprod_{i = 1}^{m} \Gamma_{i} \backslash
    \mathcal{N}^{e}_{\mathrm{red}},
  \end{equation}
  where $\Gamma_{i} = \mathbb{I}(\mathbb{Q}) \cap g_{i} C^{p}
  g_{i}^{-1} \subset J(\mathbb{Q}_{p})$ which is discrete and
  cocompact modulo center, and $g_{1}, \ldots, g_{m}$ are
  representatives of the finitely many double cosets in $\mathbb{I}(Q)
  \backslash \mathbb{G}(\mathbb{A}^{p}_{f}) / C^{p}$. Furthermore, the
  induced surjective morphism
  \begin{equation}
    \label{eq:241}
    \tilde{\Theta} \colon \coprod_{i = 1}^{m}
    \mathcal{N}^{e}_{\mathrm{red}} \longrightarrow
    \mathcal{A}^{\mathrm{ss}}_{\mathbb{F}},
  \end{equation}
  is a local isomorphism and the restriction of $\tilde{\Theta}$ to
  any closed quasi-compact subscheme of
  $\mathcal{N}^{e}_{\mathrm{red}}$ is finite.
\end{theorem}

\begin{theorem}
  \label{thr:conclusion}
  The supersingular locus $\mathcal{A}^{\mathrm{ss}}_{\mathbb{F}}$ is
  of pure dimension $[\frac{n-1}{2}]$. We have natural bijections
  \begin{equation}
    \label{eq:242}
    \{ \text{irreducible components of } \mathcal{A}^{\mathrm{ss}}_{\mathbb{F}} \}
    \stackrel{1:1}{\longrightarrow} \mathbb{I}(\mathbb{Q}) \backslash
    (J(\mathbb{Q}_{p}) / K_{\mathrm{max}} \times
    \mathbb{G}(\mathbb{A}^{p}_{f}) / C^{p}),
  \end{equation}
  and 
  \begin{equation}
    \label{eq:243}
    \{ \text{connected components of } \mathcal{A}^{\mathrm{ss}}_{\mathbb{F}} \}
    \stackrel{1:1}{\longrightarrow} \mathbb{I}(\mathbb{Q}) \backslash
    (J(\mathbb{Q}_{p}) / J^{0} \times \mathbb{G}(\mathbb{A}^{p}_{f})/
    C^{p}). 
  \end{equation}
  where $J^{0}$ is the subgroup of $J(\mathbb{Q}_{p})$ consisting of
  those $j$ with trivial Kottwitz invariant and $K_{\mathrm{max}}$ is
  the stabilizer of some maximal-type vertex lattice in
  $J(\mathbb{Q}_{p})$.
\end{theorem}

\begin{proof}
  The proof is the same as \cite[6.5]{MR2800696}.
\end{proof}

\bibliographystyle{amsalpha}
\bibliography{/home/weare1217/ownCloud/code-tex/AG}

\end{document}